\documentclass[11pt,reqno]{amsart}
\usepackage{amsfonts}
\usepackage{dsfont}
\usepackage{amsthm}
\usepackage{enumitem}
\usepackage{bbm}
\usepackage{times}
\usepackage{tabularx}
\usepackage{amsbsy}
\usepackage{mathtools}
\usepackage{tikz}
\usepackage[all]{xy}

\usetikzlibrary{arrows,decorations.markings}
\usetikzlibrary{matrix}
\usetikzlibrary{graphs}
\usetikzlibrary{backgrounds}
\usepackage{setspace}     \spacing{1}
\usepackage{color}

\usepackage[colorlinks=true,linkcolor=blue,citecolor=blue]{hyperref}
\usepackage{amssymb,latexsym}
\usepackage{mathrsfs,hyperref}
\usepackage{amsmath,latexsym}
\usepackage{amscd,latexsym} %for communtative diagram
\usepackage{float}
\usepackage{indentfirst}
\usepackage{amsfonts,latexsym}
 \usepackage{graphicx}
 \usepackage{subfig}

\newtheorem{thm}{Theorem}
\newtheorem{prop}{Proposition}%[section]
\newtheorem{lem}[prop]{Lemma}
\newtheorem{cor}[prop]{Corollary}
\newtheorem{clm}[prop]{Claim}

\theoremstyle{definition}

\newtheorem{df}[prop]{Definition} %[subsection]

\theoremstyle{remark}

\newtheorem{rmk}[prop]{Remark} %[subsection]
\newtheorem{ex}[prop]{Example}

\newcommand{\R}{{\mathbb{R}}}
\newcommand{\Z}{{\mathbb{Z}}}
\newcommand{\C}{{\mathbb{C}}}

\newcommand{\N}{{\mathbb{N}}}

\newcommand{\hL}{\widehat{L}}
\newcommand{\trv}{\mathfrak{t}_\eta}

\newcommand{\J}{\mathcal{J}}

\newcommand{\M}{\mathcal{M}}

\newcommand{\ham}{Ham}

\newcommand{\cL}{\mathcal{L}}

\newcommand{\hc}{\mathcal{C}(L,H)}

\newcommand{\ho}{\mathcal{O}(L,H)}
\newcommand{\eho}{\mathcal{O}_\eta(L,H)}

\newcommand{\al}{\alpha}

\newcommand{\ga}{\gamma}

\newcommand{\T}{\mathbb{T}}

\newcommand{\cA}{\mathcal{A}_{L,H}}

\newcommand{\cD}{\mathcal{D}}

\newcommand{\cH}{\mathcal{H}}

\newcommand{\cJ}{\mathcal{J}}

\newcommand{\cP}{\mathcal{P}}

\newcommand{\cT}{\mathcal{T}}
\newcommand{\cM}{\mathcal{M}}

\newcommand{\brat}[1]{{\left< #1 \right>}}
\newcommand{\tens}{\otimes}

\DeclareMathOperator{\supp}{\mathrm{supp}}

\DeclareMathOperator{\im}{\mathrm{Im}}
\DeclareMathOperator{\spec}{\mathrm{Spec}}

\DeclareMathOperator{\crit}{{\mathrm{Crit}}}
\DeclareMathOperator{\ind}{{\mathrm{ind}}}
\DeclareMathOperator{\virdim}{\mathrm{vir-dim}}

\numberwithin{equation}{section}

\title[Lagrangian ]{The unbounded Lagrangian spectral norm and wrapped Floer cohomology}% On a variant of Viterbo's conjecture
\author{Wenmin Gong }

\thanks{Dedicated to Professor Leonid Polterovich on the occasion of his 60th birthday}

\address{School of Mathematical Sciences, Beijing Normal University,
	Beijing, 100875, China}

\email{ wmgong@bnu.edu.cn}
\begin{document}

\maketitle
%\centerline{Preliminary version}

\begin{abstract}
We investigate the question of whether the spectral metric on the orbit space of a fiber in the disk cotangent bundle of a closed manifold, under the action of the compactly supported Hamiltonian diffeomorphism group, is bounded. We utilize wrapped Floer cohomology to define the spectral invariant of an admissible Lagrangian submanifold within a Weinstein domain. We show that the pseudo-metric derived from this spectral invariant is a valid $Ham$-invariant metric. Furthermore, we establish that the spectral metric on the orbit space of an admissible Lagrangian is bounded if and only if the wrapped Floer cohomology vanishes. Consequently, we prove that the Lagrangian Hofer diameter of the orbit space for any fiber in the disk cotangent bundle of a closed manifold is infinite.

\end{abstract}

\tableofcontents

\maketitle

\section{Introduction}\label{sec:1}

\noindent

\subsection{Viterbo's conjecture and results}
A well-known conjecture of Viterbo~\cite{Vi1} from 2007 states that, given a Riemannian metric $g$ on the $n$-dimensional torus $M=\T^n$, the spectral norm of any exact Lagrangian $L$ which is Hamiltonian isotopic to the zero section of the unit disk cotangent bundle $D_g^*M$, admits a uniform bound. Recently, Shelukhin~\cite{Sh1} has proven this conjecture, extending its validity beyond the special case of tori. In fact, Shelukhin has verified Viterbo's conjecture for various classes of closed manifolds such as compact rank one symmetric spaces $S^n, \mathbb{R}P^n, \mathbb{C}P^n, \mathbb{H}P^n$~\cite{Sh2}, as well as string point-invertible manifolds~\cite{Sh1}. Biran and Cornea~\cite{BC} established the equivalence between the Viterbo conjecture and the boundedness of the boundary depth of the Floer complex of the pair $(L, T_q^*M)$, where $(M,g)$ is any closed Riemannian manifold, and $L\subseteq D_g^*M$ is Lagrangian exact isotopic within $T^*M$ to the zero section. Despite these advancements, this conjecture remains open for arbitrary closed smooth manifolds.

On the contact side, Dimitroglou Rizell demonstrated in~\cite{Dim} that the spectral norm of Legendrians within $D^*S^1\times\R$ belonging to the contactization $ J^1 S^1$  of the cotangent bundle $T^*S^1$ of the unit circle $S^1$, which are Legendrian isotopic to the zero section, does not possess a uniform bound. This indicates that a direct generalization of Viterbo's conjecture to the Legendrian case is not feasible. In the same work, the author constructed a Hamiltonian isotopy of a closed exact Lagrangian within the $2$-torus with an open ball removed (a Liouville domain), in which the spectral norm attains arbitrarily large values. Further examples illustrating this phenomenon can be found in~\cite{Za2,KS}.

Viterbo's conjecture motivates us to consider the following
\begin{description}
  \item[Question A] Whether or not the spectral norm $\ga(F_q,L)$ of every exact Lagrangian $L$ (with boundary) in the co-disk bundle $D_g^*M$ of a closed manifold $M$, which is isotopic to a fixed fiber $F_q\subset D_g^*M$ through a Hamiltonian isotopy with support in $int(D_g^*M)$, possesses a uniform bound?
\end{description}

Before answering this question, the first problem with which we are confronted is how to define the spectral norm of a Hamiltonian deformation $L$ of a fiber under a Hamiltonian flow with support inside the unit disk cotangent bundle $D_g^*M\subset T^*M$. Originally, Viterbo~\cite{Vi} defined the spectral norm $\ga(o_M,L)$ of an exact Lagrangian $L\subset T^*M$, which is Hamiltonian isotopic to the zero section $o_M\subset T^*M$, in terms of the difference between two homological minimax values of a generating function. Later, Oh~\cite{Oh2,Oh3} provided an alternative definition using Lagrangian Floer theory (see~\cite{Mi} for the proof of the equivalence of these two constructions), which offers a more modern perspective. For a comprehensive introduction to this subject, we refer to~\cite{Oh4}.

In our case, we are interested in working with a Liouville domain $(M,d\theta)$ equipped with an exact Lagrangian submanifold $L$ intersecting the boundary $\partial M$ transversally. Somewhat surprisingly, to the author's knowledge, the spectral invariant and the associated spectral norm for such $L$ have not yet been explored. We notice that the spectral norm for Hamiltonian deformations of a closed embedded weakly exact Lagrangian in a compact or convex-at-infinity symplectic manifold has been defined by Leclercq~\cite{Le} for over a decade. Previous work on this topic by Viterbo~\cite{Vi}, Oh~\cite{Oh2,Oh3}, and subsequent research conducted by Monzner, Vichery, and Zapolsky~\cite{MVZ}, Leclercq and Zapolsky~\cite{LZ}, Kati\'{c}, Milinkovi\'{c}, and Nikoli\'{c}~\cite{KMD1,KMD2}, and Fukaya, Oh, Ohta, and Ono~\cite{FOOO} are also worth mentioning.

Let $(M^{2n},\omega=d\theta)$ be a Weinstein domain and $L^n\subset M$ be an admissible Lagrangian (cf.  Definition~\ref{df:lag}). We denote $\mathcal{H}_c(M)$ as the set of Hamiltonians $H\in C^\infty([0,1]\times M)$ with support in $[0,1]\times int(M)$ where $int(M):=M\setminus\partial M$. The Hamiltonian flow $\varphi_H^t$ associated with $H\in \mathcal{H}_c(M)$ is defined as the flow obtained by integrating the time-dependent vector field $X_{H_t}$, where $H_t=H(t,\cdot)$ and $X_{H_t}$ is determined by $-dH_t=\omega(X_{H_t},\cdot)$. The group of time-one maps of Hamiltonian flows for $H\in \mathcal{H}_c(M)$ is denoted by $\ham_c(M,\omega)$.

Using wrapped Floer cohomology, we define the spectral invariant $\ell(H,\al)$ for each $H\in\cH_c(M)$ and each non-zero class $\al\in H^*(L)$ by
\[
\ell(H,\alpha)=\sup\big\{a\in\R\setminus\spec(L,H)\big|\pi_a\circ\psi_{pss}^H(\alpha)=0\big\},
\]
see (\ref{df:siv2}) or (\ref{df:siv1}) for an equivalent definition. In this paper, we establish some basic properties of this spectral invariant, see Proposition~\ref{prop:spectral}.

The spectral pseudo-norm $\ga$ is given for $H\in\cH_c(M)$ by
\[
\ga(L,H)=-\ell(H,\mathds{1}_L)-\ell(\overline{H},\mathds{1}_L)
\]
where $\mathds{1}_L\in H^0(L)$ is the fundamental class. We denote by $$\cL(L)=\{\varphi(L)|\varphi\in\ham_c(M,\omega)\}$$ the orbit space of $L$ under the action of the group $\ham_c(M,\omega)$.

The \emph{Lagrangian spectral pseudo-metric} is given by
$$\ga(L_1,L_2)=\inf\ga(L,H),\quad \hbox{for any}\;L_1=\varphi_F^1(L),L_2=\varphi_G^1(L)\in \cL(L),$$
where the infimum runs over ${H\in\cH_c(M)}$ such that $\varphi_H^1(L)=\varphi_F^{-1}\circ\varphi_G^1(L)$.

Based on the properties of $\ell$ (cf. Proposition~\ref{prop:spectral}), we show that the spectral pseudo-metric $\ga$ satisfies the following properties.

\begin{thm}\label{nontrivality}
The pseudo-metric $\ga$ on $\cL(L)$ satisfies
\begin{itemize}
\item[(a)] $\ga(L_1,L_2)=0$ if and only if $L_1=L_2$;
\item[(b)] $\ga(L_1,L_2)=\ga(L_2,L_1)$;
\item[(c)] $\ga(\varphi(L_1),\varphi(L_2))=\ga(L_1,L_2)$ for all $\varphi\in\ham_c(M,\omega)$;
\item[(d)] $\ga(L_1,L_2)\leq\ga(L_1,L_3)+\ga(L_2,L_3)$;
\item[(e)] $\ga(L_1,L_2)\leq \delta_H(L_1,L_2)$.
\end{itemize}

\end{thm}

From this theorem, one can see that the pseudo-metric $\ga$ on $\cL(L)$ is indeed a genuine $\ham$-invariant metric. Bi-invariant metrics on $\ham_c(M,\omega)$ that satisfy these properties have appeared in the work of Viterbo~\cite{Vi} for the standard symplectic vector space $(\R^{2n},\omega_0)$ and cotangent bundles of closed manifolds, Schwarz~\cite{Sch} for closed weakly exact manifolds, Frauenfelder and Schlenk~\cite{FS} for weakly exact convex symplectic manifolds and Oh~\cite{Oh5,Oh6} in general case.

%For a Liouville domain $W$,  Mailhot~\cite{Ma} shows that the spectral diameter of $\mathcal{H}_c(W)$ is bounded if and only if the symplectic cohomology $SH^*(W)$ vanishes. 

Our main result is an open string analogue of Mailhot's (cf.~\cite[Theorem~A1]{Ma}).

\begin{thm}\label{Mthm0}
Let $(M^{2n},d\theta)$ be a Weinstein domain and $L^n\subseteq M$ an admissible Lagrangian submanifold. Then the metric space $(\mathcal{L}(L),\ga)$ is bounded if and only if the wrapped Floer cohomology of $L$ vanishes.
\end{thm}
%with the pull back $1$-form $\theta|_L=0$
\subsection{Lagrangian Hofer geometry}
 In addition to Viterbo's conjecture, another motivation of the present paper arises from Hofer geometry, for which we direct the reader to the book~\cite{Po} by Polterovich for a fascinating introduction. For a symplectic manifold $(M,\omega)$, there is a pseudo-norm on $\ham_c(M,\omega)$ defined as follows:
 \[
\|\varphi\|=\inf\bigg\{\int^1_0\bigg(\sup_{x\in M}H(t,x)-\inf_{x\in M}H(t,x)\bigg)dt\bigg|\varphi=\varphi_H^1\bigg\}.
\]

A significant result by Hofer~\cite{HZ} establishes that the pseudo-norm $\|\cdot\|$ on $\ham_c(\mathbb{R}^{2n},\omega_0)$ is a valid norm, where $\omega_0=\sum_i dx_i\wedge dy_i$. This result was later extended to general symplectic manifolds in~\cite{Po0,LM}. The Hofer norm $\|\cdot\|$ gives rise to a bi-invariant metric $d_H$ on $\ham_c(M,\omega)$ through the formula $d_H(\varphi,\psi) = \|\varphi \circ \psi^{-1}\|$.

 In~\cite{Ch}, Chekanov introduced a Lagrangian version of the Hofer metric $d_H$, where the Lagrangian Hofer pseudo-metric $\delta_H$ on the orbit space $\mathcal{L}(L)$ of a fixed Lagrangian $L\subset M$ is defined as:
\begin{equation}\label{df:dH}
\delta_H(L_1,L_2)=\inf\big\{\|\varphi\|\big|\varphi(L_1)=L_2,\;\varphi\in \ham_c(M,\omega)\big\}.
\end{equation}
It can be showed that
if $L$ is a closed Lagrangian submanifold of a tame symplectic manifold $(M,\omega)$, then $(\cL(L),\delta_H)$ is a metric space, see~\cite[Theorem~1]{Ch}.

In~\cite{Su}, Sugimoto obtained a similar result on the Lagrangian Hofer metric by replacing the closedness condition with the weekly exact condition for a Lagrangian $L$ with boundary and its completion $\widehat{L}$. Specifically, \cite[Theorem~1.4~(1)]{Su} implies the following.
%Indeed, he used wrapped Floer theory to obtain a Lagrangian version of energy-capacity inequality for the case that $L$ is a weakly exact Lagrangian submanifold in a compact symplectic manifold with a contact-type boundary.

\begin{thm}\label{thm:Hmetr}
Let $L$ be an admissible Lagrangian submanifold of a Liouville domain $(M,d\theta)$. Then $(\cL(L),\delta_H)$ is a metric space.
\end{thm}

\begin{rmk}
The nontrivial part of the proof of Theorem~\ref{thm:Hmetr} is to show the non-degeneracy of the pseudo-metric $\delta_H$.  This result can be viewed as a simple application of properties~(a) and (e) in Theorem~\ref{nontrivality}.
\end{rmk}

Let $B\subset\R^{2n}$ be the open unit ball which is endowed with the symplectic structure $\omega=\sum_{i=1}^n\frac{1}{\pi}dx_i\wedge dy_i$. Denote by $L_0=\{(x_i,y_i)\in B| y_i=0\}$ the standard Lagrangian in the unit ball. In~\cite{Kh}, Khanevsky proved that the metric space $(\cL(L_0),\delta_H)$ is unbounded in two-dimensional case, and later Seyfaddini~\cite{Se} proved the unboundedness of $(\cL(L_0),\delta_H)$ in full generality.

In the field of Lagrangian Hofer geometry, significant advancements have been made, see for instance~\cite{Hu,KS,Us1,Us2,Us3,Su,Za3,Di}. Despite this, a basic question which seems to remain open is the following problem. We refer the reader to Problem 32 on Page~562 in McDuff and Salamon's book~\cite{MS} for the related Lagrangian Hofer diameter conjecture, which is very little known.

\begin{description}
  \item[Question B] How to generalize the aforementioned Khanevsky-Seyfaddini's result to other Lagrangian submanifolds (possibly with boundary)?
\end{description}

In what follows, we will see that Theorem~\ref{Mthm0} has applications to \textbf{Questions A} and \textbf{B}.

For a closed manifold $M$, there are isomorphisms between the wrapped Floer cohomology of the fiber $F_q\subset D^*M$ at $q\in M$ (resp. the disk conormal bundle $L_K\subset D^*M$ of a compact submanifold $K\subset M$) and 
 the Morse homology of the based loop space $\Omega_qM$
 (resp.  the space $\mathcal{P}_K(M)$ of paths in $M$ with endpoints in $K$)  over $\Z/2$, which preserve the splitting of the wrapped Floer and Morse complexes determined by the homotopy classes of the generators, see~\cite{APM,AM}.  Note that the wrapped Floer cohomology defined in this paper is generated by chords in the homotopy class of constant paths, see~Section~\ref{subsec:wFch}.  So, in particular, we have
 \[HW^*(F_q)\cong H_{-*}(\Omega_q^\circ M),\quad HW^*(L_K)\cong H_{-*}(\mathcal{P}_K^\circ(M)),\]
where $\Omega_q^\circ M$ (resp. $\mathcal{P}_K^\circ(M)$) is the component of $\Omega_q M$ (resp. $\mathcal{P}_K(M)$) consisting of elements in the homotopy class of constant paths. Consequently,  $HW^*(F_q)$ and $HW^*(L_K)$ never vanish. This, together with Theorem~\ref{Mthm0}, implies the following

\begin{thm}
The metric spaces $(\cL(F_q),\ga)$ and $(\cL(L_K),\ga)$ are unbounded.
\end{thm}

In particular, the above result gives a negative answer to \textbf{Question A}. Moreover, from the inequality
 between the Lagrangian spectral and Hofer metrics
 \[\ga(L_1,L_2)\leq \delta_H(L_1,L_2),\quad \forall L_1,L_2\in\cL(L)\]
and Theorem~\ref{Mthm0}, we get
\begin{thm}
Let $L$ be an admissible Lagrangian submanifold in the Weinstein domain $(M,d\theta)$. If $HW^*(L)\neq0$, then the metric space $(\cL(L),\delta_H)$ is unbounded.
\end{thm}

The above result partially solves  \textbf{Question B}. Specifically, we observe that

\begin{cor}\label{CHbd}
The metric spaces $(\cL(F_q),\delta_H)$ and $(\cL(L_K),\delta_H)$ are unbounded.
\end{cor}

Here we note that Usher obtained a similar result to Corollary~\ref{CHbd}{\footnote{Strictly speaking, on the unboundedness of Lagrangian Hofer metric on $\cL(F_q)$ or $\cL(L_K)$, Usher's result implied by~\cite[Theorem~1.2]{Us3} is a very special case of our result since there some geometric and topological conditions are required to be satisfied.}}. In particular, Theorem~1.2 in~\cite{Us3} implies that Corollary~\ref{CHbd} holds for the closed manifold $M$ which is either a compact semisimple Lie group with a bi-invariant metric or a sphere $S^n$ of dimension at least $3$ with the standard metric.%, see~\cite[Proposition~3.10]{Us3}.

On the other hand, a classical result that the symplectic cohomology of the unit closed ball $D=\overline{B}\subset\R^{2n}$ vanishes, i.e.  $SH^*(D)=0$, implies that the wrapped Floer cohomology of the standard Lagrangian $\overline{L}_0$ also vanishes, see~Theorem~10.6 in~\cite{Rit}. Hence, a simple application of  Theorem~\ref{Mthm0} implies the following
\begin{cor}
The spectral diameter of the metric space $(\cL(\overline{L}_0),\ga)$ is finite.
\end{cor}
Comparing this assertion with the previous Khanevsky-Seyfaddini's result, we find that the spectral metric $\gamma$ and the Lagrangian Hofer metric $\delta_H$ are not equivalent on the orbit space of an admissible Lagrangian in general.

\subsection{An overview of the proof of the main result}
Our main result (Theorem~\ref{Mthm0}) is a consequence of the following two theorems, which we will prove separately.

\begin{thm}\label{Mthm}
Let $(M^{2n},d\theta)$ be a Weinstein domain and $L^n\subset M$ an admissible Lagrangian submanifold. If the wrapped Floer cohomology of $L$ vanishes, then the diameter of the metric space $(\mathcal{L}(L),\ga)$ has the upper bound $2c_{HW}(L)$.

\end{thm}
In the above theorem $c_{HW}(L)$ is the\emph{ wrapped Floer capacity} of the admissible Lagrangian $L$ (cf.  Section~\ref{sec:wFc}), which is defined by
\[
c_{HW}(L)=\inf\big\{a>0|\iota_{-a}^L\circ\psi^f_{pss}(\mathds{1}_L)=0\big\},
\]
where $\mathds{1}_L\in H^0(L)$ is the fundamental class, $\iota_{a}^L$ is the natural map (cf.~Section~\ref{sec:cont}), and $\psi^f_{pss}$ is the PSS-map (cf. Section~\ref{subsec:pss}). This capacity, as defined by Borman and McLean in~\cite[Section~1.2]{BM}, is an analogue of the symplectic cohomology capacity which is known as \emph{Floer-Hofer-Wysocki capacity} (see~\cite{FHW}).

The strategy employed to prove Theorem~\ref{Mthm} bears resemblance to the approach used by Weber in~\cite{We} to compute the Biran-Polterovich-Salamon capacities as defined in~\cite{BPS}. Weber's method was subsequently generalized to the Finsler setting in~\cite{GX}, which appears to be the initial endeavor to establish a connection between relative symplectic capacities (BPS) and Finsler/convex geometry through Floer theory.

\begin{thm}\label{Mthm'}
Let $(M^{2n},d\theta)$ be a Weinstein domain, and let $L^n\subset M$ be an admissible Lagrangian submanifold. If the wrapped Floer cohomology of $L$ does not vanish, then the metric space $(\mathcal{L}(L),\ga)$ is unbounded.
\end{thm}% such that $\theta|_L=0$

A crucial technical result in the proof of Theorem~\ref{Mthm'} is the following
\begin{prop}\label{prop:tech}
Let $(M^{2n},d\theta)$ be a Weinstein domain and $L^n\subset M$ an admissible Lagrangian submanifold. For any Hamiltonian $H\in C^\infty([0,1]\times M)$ with compact support in $[0,1]\times int(M)$, it holds that $\ell(H,\mathds{1}_L)\leq 0$.
\end{prop}

The proof of Proposition~\ref{prop:tech} is based on Ganor-Tanny's ``barricade" technique~\cite{GS}, see Section~\ref{sec:barricades}. In particular, this proposition can be viewed as a relative version of Lemma~4.1 in~\cite{GS}. In a recent work~\cite{Ma}, Mailhot has used this technique to show the unboundness for the Hamiltonian spectral diameter of Liouville domains $M$ satisfying $SH^*(M)\neq0$.

We expect that the results and techniques presented in this paper could be extended to address \textbf{Question B} for a wider class of Lagrangian submanifolds beyond the exact case in a Liouville domain. Additionally, we hope that these results can be applied to investigate other related problems in Hofer geometry. Recently, Feng and Zhang~\cite{FZ} and Dawid~\cite{Da} extend our results in the setting of large-scale geometry on cotangent bundles.

%$C(g,L)$ which only depends on a Riemannian metric $g$ on $L$ and $L$ itself.
%the spectral distance $\ga(L_1,L_2)$ between two Lagrangians $L_1,L_2$ that are Hamiltonian isotopic the zero section inside $DT^*M$ of all closed manifolds

\medskip

\section*{Acknowledgements}I am grateful to Yaniv Ganor and Shira Tanny for useful discussions.  I would like to thank warmly Jun Zhang for his valuable comments on the first draft of this work and for the thrilling discussions during my stay at Institute of Geometry and Physics, University of Science and Technology of China (USTC). I am thankful to Adrian P. Dawid for the insightful discussions and for sharing his master’s thesis, supervised by Paul Biran, which employs a different method to demonstrate the unboundedness of the Lagrangian Hofer norm on Lagrangians isotopic to a fiber in the disk cotangent bundle of any sphere $S^n(n\geq 1)$. I appreciate the helpful discussions with Mike Usher, particularly for clarifying the grading issue in Lagrangian Floer theory. I would like to express my deep gratitude to Professor Leonid Polterovich for his generous help at the beginning of my research in symplectic topology. It was through his wonderful book~\cite{Po} that I was introduced to the field of Hofer Geometry. I also thank the anonymous referee so much for a very careful reading and valuable suggestions. W.G. is partially supported by NSFC 12271285 and 11701313.

\section{Wrapped Floer cohomology}
In this section, following~\cite{AS,Rit} we briefly recall the construction of wrapped Floer cohomology for admissible Lagrangian submanifolds.

\subsection{Liouville domain and admissible Lagrangians}
Let $(M^{2n},\omega=d\theta)$ be a Liouville domain with a Liouville vector field $V_\theta$ positively transverse to $\partial M$. Here, $V_\theta$ is defined by $\iota_{V_\theta}\omega=\theta$. Then $\alpha=\theta|_{\partial M}$ is a contact form on $\partial M$. Denote by $\varphi_{V_\theta}^t$ the flow of $V_\theta$. We extend $M$ to a complete manifold by setting
\[\widehat{M}=M\bigcup_{t\geq 0}\varphi_{V_\theta}^t(\partial M).\]

Using the Liouville flow $\varphi_{V_\theta}^{\log(r)}$, we can establish an identification between $\widehat{M}\setminus M$ and $(1,\infty)\times \partial M$, and then extend the one-form $\theta$ to $\widehat{M}$ by defining $\theta=\rho\alpha$ for $\rho\in [1,\infty)$. As a result, $(\widehat{M},d\theta)$ becomes a complete exact symplectic manifold. Throughout this paper, for $r\in(0,\infty)$ we set
\[M_r=\widehat{M}\setminus((r,\infty)\times \partial M).\]

\begin{df}\label{df:lag}

Let $L^n\subset (M,d\theta)$ be a connected, orientable, exact Lagrangian submanifold with Legendrian boundary $\partial L=L\cap\partial M$ such that $V_\theta$ is tangent to $TL$ along the boundary. We call such $L$ an \emph{admissible Lagrangian} if additional conditions are met: (1) $\theta|_L=dk_L$ for a function $k_L\in C^\infty(L,\R)$ which vanishes in a neighborhood of the boundary $\partial L$; (2) the relative Chern class $c_1(M,L)\in H^2(M,L;\Z)$ satisfies $2c_1(M,L)=0$.
\end{df}

We extend an admissible Lagrangian $L$ to a non-compact exact one by setting
\[
\widehat{L}=L\bigcup_{t\geq 0}\varphi_{V_\theta}^t(\partial L)
\]
with $\theta|_{\widehat{L}}=dk_L$ for the compactly supported function $k_L$ (by abuse of notion). Clearly, in the coordinates $(1,\infty)\times \partial M\subset \widehat{M}$, the non-compact Lagrangian $\widehat{L}\setminus L$ has the form $(1,\infty)\times \partial L$. %\textbf{In what follows for a fixed admissible Lagrangian $L$ we always assume that $\theta|_{\widehat{L}\setminus L}=0$ and $k_L=0$ on a neighborhood of $\widehat{L}\setminus L$ unless otherwise specified}. %The reason why one can do this is that we can extend $k_L$ to a compactly supported smooth function $f:M\to\R$, then the new one-form $\widetilde{\theta}=\theta-df$ has the same symplectic form for the Liouville form $\theta$ and  $\widetilde{\theta}|_{\widehat{L}}=0$.

For any Lagrangian submanifold $L \subset M$, the obstruction to the existence of a grading is given by the Malov class $m_L \in H^1(L;\Z)$, which maps to the relative first Chern class $2c_1(M,L)$ under the boundary map $H^1(L;\Z) \to H^2(M,L;\Z)$, see~\cite[Section~12]{Se0}. The following typical examples satisfy $2c_1(\widehat{M},\widehat{L})=0$.

\begin{ex}
All Lagrangian planes through 0 are admissible Lagrangians in a starshaped domain $(D^{2n},d\theta)$ in $\R^{2n}$ with $\theta=\frac{1}{2}\sum_{i=1}^n(y_idx_i-x_idy_i)$.
\end{ex}

\begin{ex}
For a fiberwise starshaped domain $W$ in the cotangent bundle  $T^*M$ of a closed manifold $M$, each cotangent fiber $W\cap T_q^*M$ is an admissible Lagrangian of the Liouville domain $(W,d\theta)$, where $\theta$ is the restriction of the canonical one-form $pdq$ of $T^*M$ to $W$.
\end{ex}

\begin{ex}\label{conormal}
Let $D^*M$ be the disk cotangent bundle of a closed manifold $M$. The disk conormal bundle
\[
L_K=\big\{(q,p)\in D^*M|_K\big|\langle p,v\rangle=0\;\forall v\in T_qK\big\}
\]
of a closed submanifold of $K\subset M$ is an admissible Lagrangian with respect to the one-form $pdq|_{D^*M}$.
\end{ex}

 An almost complex structure $J$ is said to be of \emph{contact type} on $[\rho_0,\infty)\times\partial M$ for some $\rho_0>0$ if $d\rho\circ J=-\theta$ for $\rho\geq \rho_0$. A family of time-dependent smooth almost complex structures $(J_t)_{t\in[0,1]}$ is said to be \emph{$d\theta$-compatible} if $\langle \cdot,\cdot\rangle=d\theta(\cdot,J_t\cdot)$ defines a family of Riemannian metrics on $\widehat{M}$, and $J^2=-Id$ in ${\rm End}(T\widehat{M})$. Denote by $\mathcal{J}_\theta$ the set of smooth families $(J_t)_{t\in[0,1]}$ of compatible almost complex structures that are of contact type and time-independent on $[1,\infty)\times \partial M$.

To define Lagrangian Floer cohomology on $(\widehat{M},d\theta)$, we need a $C^0$-bound for solutions $u:\R\times [0,1]\to \widehat{M}$ to the $s$-dependent Floer equation
\begin{equation}\label{e:feq}
\partial_s u+J_{t}^s\big(\partial_t u-X_{H_{t}^s}(u)\big)=0
\end{equation}
with Lagrangian boundary conditions
\begin{equation}\label{e:boundary}
u(s,0), u(s,1)\in \widehat{L}\quad \hbox{for any\;} s\in\R.
\end{equation}

The following technical lemma is a consequence of a maximum principle for the above Floer equation. The proof is
standard, see for instance~\cite[Lemma~D.2]{Rit}.

\begin{lem}\label{lem:bd}
Assume that for $\rho\geq \rho_0$, $J_{t}^s$ is of contact type and independent of parameters $s,t$, and $H^s_t=h_s(\rho)$ with $\partial _sh_s'(\rho)\leq 0$. Then the function $\rho\circ u$ cannot have local maxima unless it is constant. Consequently, any solution of (\ref{e:feq}) and (\ref{e:boundary}) with asymptotics $x_-,x_+$ must lie in the region $\rho\leq \max \{\rho(x_{\pm}),\rho_0\}$.
\end{lem}

\subsection{Admissible Hamiltonians and Hamiltonian chords}
Recall that the \emph{Reeb vector field} $R$ of a contact form $\alpha$ on $\partial M$ is determined by
\[
d\al(R,\cdot)=0\quad\hbox{and}\quad \al(R)=1.
\]

Fix an admissible Lagrangian $L\subset M$. A \emph{Reeb chord} of period $T$ is a map $\ga:[0,T]\to\partial M$ satisfying
\[
\dot{\ga}=R(\ga(t))\quad\hbox{and}\quad \ga(0),\ga(T)\in\partial L.
\]

The set of all positive periods of Reeb chords is denoted by $\mathcal{R}(\partial L,\theta)$, which is known to be a closed nowhere dense set in $(0,\infty)$. We define $\tau\in(0,\infty)$ to be the \emph{minimal positive period} among all Reeb chords, i.e.  $\tau=\min\mathcal{R}(\partial L,\theta)$. If there are no Reeb chords with ends on $\partial L$, we set $\tau=\infty$.

Set
\[
\cP(\hL,\hL)=\big\{\gamma\in C^1([0,1],\widehat{M})\big|\gamma(0),\gamma(1)\in\hL\big\}.
\]
Note that $\pi_0(\cP(\hL,\hL))$ admits a special component $\eta_0$ which contains all constant paths at points of $\hL$.

A \emph{Hamiltonian chord }of a smooth Hamiltonian function $H\in C^\infty([0,1]\times \widehat{M})$ is an element of $\cP(\hL,\hL)$ with properties
\[
\dot{x}=X_H(x(t)) \quad\hbox{and}\quad x(0),x(1)\in \widehat{L},
\]
where $X_H$ is the Hamiltonian vector field given by $d\theta(X_H,\cdot)=-dH_t$ with $H_t:=H(t,\cdot)$.

Clearly,  Hamiltonian chords of $H$ correspond to the intersection points in $\varphi^1_H(\widehat{L})\cap \widehat{L}$, where $\varphi_H^1$ is the time one map of the flow of $X_H$. Moreover, for any Hamiltonian chord of $x$ of $H$ in the region $(0,\infty)\times \partial M$ where $H=h(\rho)$, the function $\rho(x(t))$ is constant, and therefore $H(x(t))=h(\rho)$, and $x$ corresponds to the Reeb chord $\gamma$ given by $x(t)=(\rho,\gamma(t\cdot T))$ with period $T=|h'(\rho)|$.

We call $H\in C^\infty([0,1]\times \widehat{M})$ an \emph{admissible Hamiltonian}
 if there are constants $\rho_0>0$ and $a\in\R$ such that $H$ has the form
\begin{equation}\label{radial}
H(t,\rho,x)=\mu_H\rho+a\quad \hbox{for}\;\rho\geq \rho_0\;\hbox{on}\;(0,\infty)\times \partial M.
\end{equation}
Here, $\mu_H\notin \mathcal{R}(\partial L,\theta)$ is some non-negative number, called the \emph{slope} of $H$. Throughout the paper, we denote by $\mathcal{H}$ the set of admissible Hamiltonians, and $\cH_{<\tau}\subset\cH$ the set of admissible Hamiltonians $H$ which are linear in the region $\rho\geq 1$ with $0\leq \mu_H<\tau$.

%if there exists a function $h\in C^\infty\big(\R)$ with $0\leq h'(\rho)<\tau$ such that\[H|_{[0,1]\times[1,\infty)\times \partial M}=h\;\hbox{for}\;\rho\geq 1.\]

\subsection{The actions and indices of Hamiltonian chords}

Given an admissible Lagrangian $L\subset M$, an admissible Hamiltonian $H\in\cH$ and a free homotopy class $\eta\in\pi_0(\cP(\hL,\hL))$, we denote by $\mathcal{O}_\eta(L,H)$ the set of Hamiltonian chords belonging to the class $\eta$. We set
\[\ho=\bigcup_{\eta\in\pi_0(\cP(\hL,\hL))}\eho.\]

Given $\eta\in\pi_0(\cP(\hL,\hL))$, we choose a reference path $\gamma_\eta:[0,1]\to \widehat{M}$ representing the class $\eta$, and fix a symplectic trivialization $\mathfrak{t}_\eta:\gamma_\eta^*TM\to [0,1]\times\R^{2n}$ which maps the tangent space $T_{\gamma_\eta(i)}\widehat{L}$ to $ \{i\}\times\R^n\times\{0\}, i=0,1$. We consider the pair $(\gamma,v)$, where $v$ is a capping of $\gamma$, that is to say, a piecewise smooth map $v:[0,1]\times [0,1]\to \widehat{M}$ satisfies the following conditions: $v(s,i)\in\widehat{L}$ for $i=0,1$, $v(0,t)=\gamma_\eta(t)$, and $\gamma(1,t)=\gamma(t)$.

For the distinguished point class $\eta_0\in \pi_0(\cP(\hL,\hL))$,
we fix a point $pt\in \widehat{L}$, and choose a capping of each chord $x\in \mathcal{O}_{\eta_0}(L,H)$ as a half disk $v:\mathbb{D}_+\to M$ such that $v(e^{i\pi t})=x(t)$ for any $t\in[0,1]$, $v(t)\in \widehat{L}$ for any $t\in [-1,1]$, and $v(0)=pt$, where $\mathbb{D}_+$ is a half disk, i.e.  $\mathbb{D}_+=\{z\in\C:\;|z|\leq 1,\;\hbox{im}(z)\geq 0\}$.

For any $\eta\in\pi_0(\cP(\hL,\hL))$,  the Floer functional $\cA$ on the path component $\cP_\eta(\hL,\hL)$ is defined by
\begin{equation}\label{func}
\cA(x)=-\int v^*d\theta+\int^1_0 H(x(t))dt,
\end{equation}
where $v$ is any capping of $x\in\cP_\eta(\hL,\hL)$.  Due to the exactness of $\theta|_{\widehat{L}}$, $\cA(x)$ is independent of the choice of the capping disk $v$. In particular, by Stokes' theorem,
 for each chord $x\in\cP(\hL,\hL)$ belonging to $\eta_0$ (using the specified capping as above), one can rewrite
\[
\cA(x)=-\int x^*\theta+\int^1_0 H(x(t))dt+k_L\big(x(1)\big)-k_L\big(x(0)\big).
\]

%We denote $\widetilde{\cP}(\hL,\hL)$ the set of equivalence classes of such pairs $(x,v)$ with $x\in\cP(\hL,\hL)$. Then for each class $[x,v]$ we define its action as
%\begin{equation}\label{func}\cA([x,v])=-\int v^*d\theta+\int^1_0 H(x(t))dt\end{equation} which is independent of the choices of the capping $v$.

%Due to this fact, in the following we briefly denote $\cA(x)$ as the action of a capped chord $\overline{x}$ with $x\in \mathcal{O}_{\eta_0}(L,H)$.

It is easy to verify that the critical points of the action functional $\cA$ are exactly Hamiltonian chords with ends on $\widehat{L}$. We denote by $\spec(L,H)$ the set of actions of chords $x\in\mathcal{O}(L,H)$, which is referred to as the \emph{action spectrum} for $(L,H)$. It is well known that the set $\spec(H,L)$ is a closed nowhere dense subset of $\R$, see~\cite{LZ}.

For a Hamiltonian $H$ satisfying (\ref{radial}), if a chord $x\in\mathcal{O}_{\eta_0}(L,H)$ lies in $(\rho_0,\infty)\times \partial M$, then \begin{equation}\label{action}\cA(x)=k_L(x(1))-k_L(x(0))-\rho h'(\rho)+h(\rho).\end{equation} In particular, if both $x(0)$ and $x(1)$ lie in the neighborhood of $\widehat{L}\setminus L$  where $k_L$ vanishes, then the action of $x$ is equal to the $y$-intercept of the tangent line of the function $y=h(\rho)$ at $\rho$.

A chord $x\in \ho$ is said to be \emph{non-degenerate }if the vector spaces $T_{x(1)}\widehat{L}$ and $d\varphi_H^1(T_{x(0)}\widehat{L})$ are transverse. We call an admissible Hamiltonian $H\in\cH$ \emph{non-degenerate with respect to $\widehat{L}$} if all chords $x\in\ho$ are non-degenerate. The set of such Hamiltonians is denoted by $\cH^{reg}\subset\cH$. Throughout this paper, we denote $\cH_{<\tau}^{reg}:=\cH_{<\tau}\cap\cH^{reg}$.

To each non-degenerate chord $x\in \mathcal{O}_{\eta}(L,H)$, we now associate a Maslov-type index $\mu(x)$ as follows. Pick any capping $v:[0,1]\times [0,1]\to \widehat{M}$ of $x$ such that $v(0,t)=\gamma_\eta(t)$, $v(1,t)=x(t)$, $v(s,0), v(s,1)\in \widehat{L}$ for all $s,t\in[0,1]$. We denote the trivialization of the pullback bundle $v^*T\widehat{M}$ by $\mathfrak{t}_v$. This trivialization satisfies that the restriction of $\mathfrak{t}_v$ over $\gamma_\eta$ is the fixed trivialization $\trv$ as above, and $\mathfrak{t}_v$ maps $T_{v(s,i)}\widehat{L}$ to $\R^n\times\{0\}$.

The path $\Lambda(t)=\mathfrak{t}_{v(1,t)}(\varphi_H^t)_*T_{x(0)}\widehat{L}=\mathfrak{t}_{v(1,t)}(\varphi_H^t)_*\mathfrak{t}^{-1}_{v(1,0)}(\R^n\times\{0\})$ in the Lagrangian Grassmannian satisfies that $\Lambda(0)=\R^n\times\{0\}$ and $\Lambda(1)$ is transverse to $\R^n\times\{0\}$. Following~\cite[Section~5]{RS}, for a fixed Lagrangian $V$ in the Lagrangian Grassmannian $\mathcal{L}_n$ for $(\R^{2n},\sum_idx_i\wedge dy_i)$ and a path $\Lambda:[a,b]\to\mathcal{L}_n$, we denote by $\mu_{RS}(\Lambda;V)$ the Robbin-Salamon Maslov index. We define the Malov index of the chord $x$ as
\[
\mu(x)=-\mu_{RS}(\Lambda;\R^n\times\{0\})+\frac{n}{2}\in\Z
\]

Our assumption $2c_1(M,L)=0$ ensures that this integer is independent the choice of the capping $v$.
The grading $\mu$ is normalized so that if $H$ is a lift of a $C^2$-small Morse function $f$ on $\widehat{L}$ to a Weinstein neighborhood of $\widehat{L}$, then for any constant path $\gamma_q$ at a critical point $q$ of $f$, $\mu(\gamma_q)$ coincides with the Morse index of $q$.

%For a detailed construction of this index we refer to???. Here we remark that $\mu$ is a relative Maslov index which depends on the base point in the set of capped Hamiltonian chords. For our purposes we require that if $H_f$ is a lift of a $C^2$-small Morse function $f$ on $\widehat{L}$ to a Weinstein neighborhood of $\widehat{L}$, then for a constant chord $q\in\crit(f)$ with constant capping disk $\overline{q}$ we have $\mu(\overline{q})=m_f(q)$.

\subsection{Lagrangian Floer cohomology}\label{sec:LFC}
Given $J\in\mathcal{J}_\theta$ and $H\in\cH^{reg}$, consider the solutions $u:\R\times[0,1]\to M$ to  the Floer equation
\begin{equation}\label{f1}
\partial_s u+J_{t}\big(\partial_t u-X_{H_{t}}(u)\big)=0
\end{equation}
with boundary conditions $u(\R,i)\in \widehat{L},i=0,1$. The \emph{energy} of a solution $u$ to~(\ref{f1}) is defined as
\[
E_J(u)=\frac{1}{2}\int^1_0\int^1_0\big|\partial_su\big|_J^2+\big|\partial_tu-X_H(u)\big|_J^2dsdt,
\]
where $|\cdot|_J$ is the norm with respect to the metric $\omega(\cdot,J\cdot)$. We sometimes use the notation $E(u)$ whenever the almost complex structure is clear from the context.

Denote by $\widehat{\mathcal{M}}_{H,J}(x,y)$ the space of Floer solutions $u$ of finite energy such that
\[
\lim\limits_{s\to-\infty}u(s,t)=x,\qquad \lim\limits_{s\to+\infty}u(s,t)=y
\]
There is a natural $\mathbb{R}$-translation action in $s$-direction on $\widehat{\mathcal{M}}_{H,J}(x,y)$. The quotient space under this action is denoted by $\mathcal{M}_{H,J}(x,y)$. Solutions of (\ref{f1}) can be thought as negative gradient flow lines for $\cA$ in an $L^2$-metric on $\ho$. For each $u\in \widehat{\mathcal{M}}_{H,J}(x,y)$, we linearize (\ref{f1}) to obtain a Fredholm operator $D_{H,J,u}$ in suitable Sobolev spaces. Since $H\in\cH$ and $J\in\mathcal{J}_\theta$, it follows from Lemma~\ref{lem:bd} that all $u\in\mathcal{M}_{H,J}(x,y)$ have a  uniform $C^0$-bound in $\R$-component of the coordinates in $[1,\infty)\times \partial M$.
Besides, in our situation that the symplectic form $d\theta$ is exact and $\widehat{L}\subset \widehat{M}$ is an exact submanifold, bubbling of pseudo-holomorphic disks or spheres does not occur in the limits of families of solutions of~(\ref{f1}). Hence, there is a dense subspace $\mathcal{J}^{reg}\subset \mathcal{J}_\theta$ of almost complex structures such that for each $J\in\mathcal{J}^{reg}$, $D_{H,J,u}$ are onto for all $u\in \widehat{\mathcal{M}}_{H,J}(x,y)$ (see~\cite{Fl1,FHS}). In this situation, we call $(H,J)$ \emph{Floer-regular} or \emph{regular}. As a consequence, the space
$\widehat{\mathcal{M}}_{H,J}(x,y)$, as well as $\mathcal{M}_{H,J}(x,y)$, is a smooth manifold whose component representing  each homotopy class $[u]$  has dimension (see~\cite{Fl1,FOOO,Au})
\begin{equation}\label{dim}
\dim_{[u]}\widehat{\mathcal{M}}_{H,J}(x,y)=\mu(x)-\mu(y).
\end{equation}

Let $CF^*(\widehat{L},H)$ be the vector space over $\mathbb{Z}/2$ generated by chords $x\in\ho$.
When $\mu(x)=\mu(y)+1$, one can define the Floer differential $$d_{H,J}:CF^*(\widehat{L},H)\to CF^*(\widehat{L},H)$$ by counting isolated points in $\mathcal{M}_{H,J}(x,y)$ mod $2$, i.e.
$$d_{H,J}(y)=\sum_{\mu(x)=\mu(y)+1}\sharp_2\mathcal{M}_{H,J}(x,y)\cdot x.$$
This map has square zero, i.e.  $d_{H,J}^2=0$, and hence $(CF^*(\widehat{L},H),d_{H,J})$ is a cochain complex over the coefficient $\mathbb{Z}/2$, called the \emph{Lagrangian Floer complex} of $(L,H)$. The \emph{Lagrangian Floer cohomology for $(\widehat{L},H)$} is defined to be the quotient space $\mathrm{Ker}(d_{H,J})/\mathrm{Im}(d_{H,J})$, which is denoted by $HF^*(\widehat{L};H)$. It can be shown that the cohomology $HF^*(\widehat{L};H)$ is independent of $J\in\cJ^{reg}$ up to canonical isomorphisms, see~\cite{SZ}. Due to this, we suppress it from the notation. From the construction of Lagrangian Floer cohomology,  we have the decomposition
\[
HF^*(\widehat{L};H)=\bigoplus_{\eta\in\pi_0(\cP(\hL,\hL))} HF^*_\eta(\widehat{L};H).
\]

\subsection{The filtered wrapped Floer cohomology}\label{subsec:wFch}
For $H\in\cH^{reg}$, we define the\emph{ wrapped Floer cohomology} of $(L,H)$ to be the Lagrangian Floer cohomology $HF^*_{\eta_0}(\widehat{L};H)$ (recall that $\eta_0$ is the point class in $\pi_0(\cP(\hL,\hL))$). For simplifying notations, hereafter we denote $$\hc=\mathcal{O}_{\eta_0}(L,H),\quad CW^*(L,H)=CF^*_{\eta_0}(\widehat{L},H),\quad HW^*(L,H)=HF^*_{\eta_0}(\widehat{L};H).$$

For $a\in\R\cup\{\pm\infty\}$, we define
\[
CW^*_{(a,+\infty)}(L,H)=\big\{x\in\hc\big|\cA(x)>a\big\}.
\]
Since the action $\cA$ does not increase along its negative gradient flows, the differential $d_{H,J}$ increases the action $\cA$. Hence, the vector space $CW^*_{(a,+\infty)}(L,H)$ is a subcomplex of the wrapped Floer complex $CW^*(L,H)$.
For $b>a$, the homology of the quotient
\[
CW^*_{(a,b]}(L,H)=CW^*_{(a,+\infty)}(L,H)/CW^*_{(b,+\infty)}(L,H)
\]
is called the \emph{filtered wrapped Floer cohomology of $(L,H)$} with the \emph{action window} $(a,b]$. We denote it by $HW^*_{(a,b]}(L,H)$.

If $a<b<c$, the inclusion map $CW^*_{(b,c]}(L,H)\to CW^*_{(a,c]}(L,H)$ induces a map
\begin{equation}\label{win1}
\iota^{a,c}_{b,c}:HW^*_{(b,c]}(L,H)\longrightarrow HW^*_{(a,c]}(L,H),
\end{equation}
and the quotient map $CW^*_{(a,c]}(L,H)\to CW^*_{(a,b]}(L,H)$ induces a map
\begin{equation}\label{win2}
\pi^{a,c}_{a,b}:HW^*_{(a,c]}(L,H)\longrightarrow HW^*_{(a,b]}(L,H).
\end{equation}
These maps are called \emph{action window maps}.
In particular, we have the natural maps induced by the inclusion and quotient maps
\[
\iota_a: HW^*_{(a,\infty)}(L,H)\longrightarrow HW^*(L,H),
\quad
\pi_a:HW^*(L,H)\longrightarrow HW^*_{(-\infty,a]}(L,H).
\]

For any admissible Hamiltonian $H\in\cH$, when $a,b\in\R\cup\{\pm\infty\}\setminus\spec(L,H)$, one can define the wrapped Floer cohomology for the pair $(L,H)$ by setting
\[
HW^*_{(a,b]}(L,H):=HW^*_{(a,b]}(L,\widetilde{H}),
\]
where $\widetilde{H}\in\cH^{reg}$ is a $C^2$-small perturbation of $H$ and has the same slope as $H$ at infinity. This definition does not depend on the choices of perturbations $\widetilde{H}$, see for instance~\cite{BPS,We,GX}.

\subsection{Energy estimates and continuation maps}\label{sec:cont}

Let $s\mapsto (H_{t}^s,J_{t}^s)$ be a path of families of admissible Hamiltonians and almost complex structures, which is constant at the ends and satisfies
\begin{equation}\notag
\big(H^{\pm\infty},J^{\pm\infty}\big)=(H^{\pm},J^\pm)\in \cH^{reg}\times\cJ^{reg}.
\end{equation}

We consider the solutions $u:\R\times[0,1]\to M$ to the $s$-dependent Floer equation
\begin{equation}\label{e:sfeq}
\partial_s u+J_{t}^s\big(\partial_t u-X_{H_{t}^s}(u)\big)=0
\end{equation}
with the boundary conditions $u(\R,0), u(\R,1)\subset L$. For non-degenerate chords $x_{\pm}\in\mathcal{C}(L,{H^\pm})$, we denote by
$\mathcal{M}_{H^s,J^s}(x_-,x_+)$ the set of finite energy solutions to (\ref{e:sfeq}). If the conditions concerning $(H,J)$ in Lemma~\ref{lem:bd} hold, as described in Section~\ref{sec:LFC}, for a generic path $J^s\subset\cJ_\theta$ the spaces $\mathcal{M}_{H^s,J^s}(x_-,x_+)$ are finite dimensional manifolds of local dimension given by~(\ref{dim}). There is another way to achieve regularity, that is to say,  we fix a path of almost complex structures $J^s\in\mathcal{J}_\theta,s\in \R$, then perturb the homotopy $H^s$ in the function space of homotopies connecting $H^-$ to $H^+$
such that the linearized operators $D_{H,J,u}$ for the resulting homotopy $H$ are onto. In both cases above, we call the pair $(H,J)$ \emph{Floer-regular} or \emph{regular}.

The energy of any solution to (\ref{e:sfeq}) satisfies the identity
\begin{equation}\label{E}
E(u)=\mathcal{A}_{L,H^-}\big(x_-\big)-\mathcal{A}_{L,H^+}\big(x_+\big)+\int_{\R\times [0,1]}\big(\partial_sH^s_t\big)\big(u(s,t)\big)dsdt.
\end{equation}

%with action in the window $(a,b]$

For $x_+\in\mathcal{C}(L,H^+)$, we define the map
\[
\Phi_{H^+H^-}:\big(CW^*(L,H^+),d_{H^+,J^+}\big)\longrightarrow \big(CW^*(L,H^-),d_{H^-,J^-}\big),
\]
\[
\Phi_{H^+H^-}(x_+)=\sum_{x_-\in\mathcal{C}(L,H^-)\atop \mu(x_-)=\mu(x_+)}\sharp_{\Z_2}\mathcal{M}^0_{H^s,J^s}(x_-,x_+)x_-,
\]
where $\mathcal{M}^0_{H^s,J^s}(x_-,x_+)$ denotes the zero dimensional component of $\mathcal{M}_{H^s,J^s}(x_-,x_+)$. This a chain map, and we call $\Phi_{H^+H^-}$ a \emph{continuation map} from $HW^*(L,H^+)$ to $HW^*(L,H^-)$.

From the energy identity (\ref{E}), one can see that if
\begin{equation}\label{e:MONO}
\int_{\R\times [0,1]}\big(\partial_sH^s_t\big)\big(u(s,t)\big)dsdt\leq 0,
\end{equation}
then $\Phi_{H^+H^-}$ preserves the action filtration, and hence induces a homomorphism
\[
\Phi_{H^+H^-}:HW^*_{(a,b]}(L,H^+)\longrightarrow HW^*_{(a,b]}(L,H^-)
\]
for $a,b\in[-\infty,\infty)$.

In case that $\mu_{H^+}\leq \mu_{H^-}$, the linear homotopy $H^s$ from $H^-$ to $H^+$
\begin{equation}\label{eq:hmtp}
H^s=\beta(s)H^++\big(1-\beta(s)\big)H^-
\end{equation}
is monotone at infinity, i.e.  $\partial_sH^s_t\leq 0$ in the region $\rho>R$ with sufficiently large $R>0$ where both $H^+$ and $H^-$ are linear with respect to the radial coordinate. Here $s\mapsto\beta(s)$ is a smooth cutoff function on $\R$ such that $\beta(s)=0$ for $s\leq 0$, $\beta(s)=1$ for $s\geq 1$, and $\beta'(s)\geq 0$. So we have the continuation map
\[\Phi_{H^+H^-}:HW^*(L,H^+)\longrightarrow HW^*(L,H^-)\]
defined as above, and we call it a \emph{monotone homomorphism}. It can be shown that monotone homomorphisms are independent of the choices of monotone homotopies $(H^s,J^s)$ used to define them, see~\cite{SZ}. Moreover, these maps are natural, meaning that $\Phi_{HH}=Id$ and
\[
\Phi_{GH}\circ\Phi_{FG}=\Phi_{FH}
\]
for admissible Hamiltonians $F\leq G\leq H$.

If $H\leq K$, then for any $a\in\R\cup\{-\infty\}$ 
the natural maps commute with the monotone homomorphisms
\begin{equation}\label{diag:IQ}
\begin{split}
\xymatrix{
&
HW^*_{(a,\infty)}(L,H)
\ar[r]^{\;\iota_a}
\ar[d]^{\Phi_{HK}}
&
HW^*(L,H)
\ar[d]^{\Phi_{HK}}
\ar[r]^{\pi_a\quad}
&
HW^*_{(-\infty,a]}(L,H)
\ar[d]^{\Phi_{HK}}
\\
&
HW^*_{(a,\infty)}(L,K)
\ar[r]^{\;\iota_a}
&
HW^*(L,K)
\ar[r]^{\pi_a\quad}
&
HW^*_{(-\infty,a]}(L,K)
}
\end{split}
\end{equation}

Using Lemma~\ref{lem:bd}, one can prove the following.

\begin{lem}~\label{lem:slope}
Assume that $H,K\in\cH$ are two admissible Hamiltonians with slopes $\mu_H$ and $\mu_K$ at infinity, respectively. If
$\mu_H=\mu_K$, then the continuation map $\Phi_{HK}$ is an isomorphism. When $(\mu_H,\mu_K)\cap \mathcal{R}(\partial L,\theta)=\emptyset$, we still have the isomorphism $\Phi_{HK}:HW^*(L,H)\stackrel{\cong}{\to} HW^*(L,K)$.
\end{lem}

The proof of Lemma~\ref{lem:slope} is parallel to~\cite[Theorem~2.5]{We}, see also~\cite[Lemma~4.1]{GX}.

The following lemma is similar to the one in the Hamiltonian Floer homology case, see for instance~\cite{Vi,BPS,Go1}.
\begin{lem}~\label{lem:mm}%monotone isomorphism
Let $(H^s)_{s\in[0,1]}$ be a homotopy connecting $K$ to $H$ satisfying~(\ref{e:MONO}). Assume that $-\infty \leq a<b<\infty$, and that $b_s$ is a continuous family of numbers such that $a<b_s\notin \spec(L,H^s)$ for every $s\in[0,1]$ with $b_0=b$ and $b_1=c$. Then we have the isomorphism
$$\Phi_{HK}: HW^*_{(a,b]}(L,H)\stackrel{\cong}{\longrightarrow} HW^*_{(a,c]}(L,K).$$
\end{lem}

Let $X\subset M$ be a compact subset of $X$. We set\[\cH^X=\big\{H\in\cH|H|_{[0,1]\times X}<0\;\hbox{and}\;H\;\hbox{satisfies}\;(\ref{radial})\;\hbox{for some}\;\rho_0\geq 1\big\}.\]

Now we introduce the partial-order relation $\preceq$ on the set $\cH$ by
\[
H\preceq K\Longleftrightarrow H\leq K.
\]

Using this relation we get a partially ordered system $(HW,\chi)$ of $\Z_2$-vector spaces over $\cH^X$. Here, $HW$ assigns to every $H\in\cH^X$ the $\Z_2$-vector space $HW^*_{(a,b]}(L,H)$, and $\chi$ assigns to any two Hamiltonians $H,K\in\cH^X$ with $H\preceq K$ the monotone continuation map $\Phi_{HK}$. Since $(\cH^X,\preceq)$ is a directed system, we have
\begin{df}\label{df:wh}
Let $L\subseteq (M,d\theta)$ be an admissible Lagrangian.
For $a<b$, the \emph{wrapped Floer cohomology of $L$ relative to $X$} is defined by
\[
HW^*_{(a,b]}(L,X)=\lim_{\substack{\longrightarrow \\ H\in\cH^X}}HW^*_{(a,b]}(L,H).
\]
\end{df}

%Taking the direct limit with respect to monotone Hamiltonians in (\ref{win1}) and (\ref{win2}) we obtain
%\[\iota^{a,c}_{b,c}:HW^*_{(b,c]}(L,X)\longrightarrow HW^*_{(a,c]}(L,X),\]
%\[\pi^{a,c}_{a,b}:HW^*_{(a,c]}(L,X)\longrightarrow HW^*_{(a,b]}(L,X).\]

When $X=M$, the cohomology defined above is the usual wrapped Floer cohomology, which is denoted by
$HW^*_{(a,b]}(L)$. If $H\in\cH$ is negative on $M$, we have the natural homomorphism
\[
\sigma_H:HW^*_{(a,b]}(L,H)\longrightarrow HW^*_{(a,b]}(L)
\]
which sends each element to its equivalence class.
Clearly, for any $H,K\in\cH^M$ with $H\preceq K$, we have the commutative diagram:

\begin{equation}\label{diam:hk}
\begin{split}
\xymatrix{
HW^*_{(a,b]}(L,H)
\ar[dr]_{\sigma_H}
 \ar[rr]^{\Phi_{HK}} &&% arrow cannot occupy space
HW^*_{(a,b]}(L,K)
\ar[dl]^{\sigma_K}
\\
&HW^*_{(a,b]}(L)
}
\end{split}
\end{equation}
Hence, taking the direct limits in (\ref{diag:IQ}) yields
\begin{equation}\label{diag:IQw}
\begin{split}
\xymatrix{
&
HW^*_{(a,\infty)}(L,H)
\ar[r]^{\;\iota_a}
\ar[d]^{\sigma_H}
&
HW^*(L,H)
\ar[d]^{\sigma_H}
\ar[r]^{\pi_a\quad}
&
HW^*_{(-\infty,a]}(L,H)
\ar[d]^{\sigma_H}
\\
&
HW^*_{(a,\infty)}(L)
\ar[r]^{\;\iota_a^L}
&
HW^*(L)
\ar[r]^{\pi_a^L\quad}
&
HW^*_{(-\infty,a]}(L)
}
\end{split}
\end{equation}

Since there exists a cofinal family $(H_k)_{k\in\N}$ of Hamiltonians in $\cH^M$ such that the Hamiltonian chords of every $H_k$ lying in $M$ are constant ones (whose actions are negative obviously), for any $a>0$ we have
\[
HW^*(L):=HW^*_{(-\infty,\infty)}(L)=HW^*_{(-\infty,a]}(L).
\]

We emphasize here that the wrapped Floer cohomology $HW^*(L)$ defined in this paper is equivalent to the definition of the wrapped Floer cohomology in~\cite{AS,Rit}, where admissible Hamiltonians are required to be positive on $M$. This is because adding a constant to any Hamiltonian function does not change the Hamiltonian flow nor Floer equations, but only shifts the action of a Hamiltonian chord by the constant.

\subsection{The ring structure of wrapped Floer cohomology}
Following~\cite[Section~6.12]{Rit}, we outline the construction of the product structure of wrapped Floer cohomology $HW^*(L)$.

Given regular pairs $(H^i,J^i),i=0,1,2$ in the definition of $HW^*(L,H^i)$, we shall define a product
\begin{equation}\label{e:pd}
*_F:HW^*(L,H^0)\otimes HW^*(L,H^1)\longrightarrow HW^*(L,2H^2)
\end{equation}
which makes $HW^*(L)$ into a ring by taking direct limits.

Consider a $2$-dimensional disk $\mathcal{D}$ with three boundary points $z^i\in\partial \mathcal{D},i=0,1,2$ removed. Let $j$ be a complex structure on $\mathcal{D}$.  Near every boundary puncture we equip a \emph{strip-like end} which can be biholomorphically mapped onto the semi-infinite strips
\[
Z_\pm=\R_\pm\times[0,1]
\]
with the standard complex structure, i.e.  $j\partial_s=\partial_t$. More precisely, for $i=0,1$ we consider a positive strip-like end near $z_i$ which is a holomorphic embedding
\[
\kappa_i:\R_+\times[0,1]\longrightarrow \mathcal{D}
\]
satisfying
\[
\kappa_i^{-1}(\partial \mathcal{D})=\R_+\times\{0,1\}\quad\hbox{and}\;\lim_{s\to+\infty}\kappa_i(s,\cdot)=z_i.
\]
Near $z_2$ we equip a negative strip-like end in a similar way. Moreover, we require that these strip-like ends (the images of $\kappa_i$) are pairwise disjoint.

To define the moduli space of the product structure, we need to choose a compatible perturbation data. Let $\alpha\in\Omega^1(\mathcal{D})$ be a $1$-form whose restriction to the boundary $\partial \mathcal{D}$ is zero, and which satisfies $\kappa_i^*\alpha=dt, i=0,1, \kappa_2^*\alpha=2dt$ and $d\alpha\leq 0$ (such form $\alpha$ exists, see~\cite{Rit}). Let $J^\mathcal{D}$ be a family of almost complex structures parameterized by $\mathcal{D}$, which satisfies $J^\mathcal{D}_{\kappa_i(s,t)}=J^i_t$ for $i=0,1,2$, and is of contact type at the cylindrical ends of $\widehat{M}$. Let $H^\mathcal{D}$ be a family of Hamiltonians $(H^\mathcal{D}_z)_{z\in \mathcal{D}}$ parameterized by $\mathcal{D}$ such that for every $z\in \mathcal{D}$, $H_z$ is an admissible Hamiltonian with $H_z^\mathcal{D}=\mu\rho+a$ outside $M$ for two constants $\mu>0$ and $a\leq0$, and $H^\mathcal{D}_{\kappa_i(s,t)}=H^i_t$ for $i=0,1,2$. Besides,  for each $x\in \widehat{M}$ we impose the \emph{non-positive condition} on the $2$-form $\beta(x):=d(H^\mathcal{D}_z(x)\alpha)$:
\begin{equation}\label{nonpos}
\beta(x)(V,jV)\leq 0
\end{equation}
for every tangent vector $V$ along $\mathcal{D}$. Such $2$-form $\beta$ exists by a suitable choice of $H^\mathcal{D}$ for special Hamiltonians $H^i,i=0,1,2$, see~\cite[Section~2.6]{EO}.
%Under this condition, fixing $x\in \widehat{M}$ and integrating $\beta(x)$ on the complement of a small neighborhood of the punctures yields \[\int^1_0 H^0_t(x)dt+\int^1_0 H^1_t(x)dt\leq \int^1_0 H^2_t(x)dt.\]In particular, we have $\mu_{H^0}+\mu_{H^1}\leq \mu_{H^2}$. %Conversely, if $H^0+H^1\leq H^2$ on $[0,1]\times \widehat{M}$,  one can find a $1$-form $\alpha$ and a family of Hamiltonians $H^\mathcal{D}$ so that the non-positive condition holds.

Consider the inhomogeneous $\overline{\partial}$-equation
\begin{equation}\label{e:inhom}
\begin{cases}
u:\mathcal{D}\longrightarrow \widehat{M},\quad u(\partial \mathcal{D})\subset \widehat{L},\\
\lim_{s\to\pm \infty}u(\kappa_i(s,\cdot))=z_i,\quad i=0,1,2,\\
\big(du_z-X_{H^\mathcal{D}}(u(z))\otimes \alpha_z\big)\circ j-J^\mathcal{D}_{z}(u(z))\circ \big(du_z-X_{H^\mathcal{D}}(u(z))\otimes \alpha_z\big)=0.
\end{cases}
\end{equation}
where $X_{H^\mathcal{D}}$ is the Hamiltonian vector field generated by $H^\mathcal{D}_z$, and each $z_i$ is a chord of the Hamiltonian $H_i$ with respect to $\widehat{L}$.

The energy of each solution $u$ to (\ref{e:inhom}) is defined by
\[
E(u)=\frac{1}{2}\int_\mathcal{D}|du-X\otimes \alpha|^2{\rm Vol}_\mathcal{D},
\]
where $|\cdot|$ is a norm induced by the symplectic form $d\lambda$, the complex structures $J^\mathcal{D}$ and $j$, and ${\rm Vol}_\mathcal{D}$ is the volume form on $\mathcal{D}$.

The following no escape lemma follows from \cite[Lemma 7.2]{AS}, see also~\cite[Lemma~D.6]{Rit} or~\cite[Lemma~3.1]{EO}, we include
the proof for completeness.
\begin{lem}\label{lem:bound}
Under the above assumptions, $u(\mathcal{D})$ lie in $M$ for all solutions $u$ to (\ref{e:inhom}).
\end{lem}

\begin{proof}
It suffices to prove that if the solution $u$ lands outside $M$ then $u(\Sigma)\subset\partial M$, where $\Sigma:=u^{-1}(\widehat{M}\setminus int(M))$.
We notice that $\Sigma$ is a compact surface with corners which divide the boundary $\partial \Sigma$ into two pieces: the piece landing in the boundary $\partial M$ and the one landing in $\widehat{L}$. Write $\partial \Sigma=\partial_b \Sigma\cup \partial_l \Sigma$ according to these two pieces. %Clearly, by our assumption $\partial_b \Sigma\neq\emptyset$.
We denote by $j$ the restriction of the complex structure from $\mathcal{D}$ to $\Sigma$.

By Stokes' theorem and the non-positive condition on $\beta$, we have the energy estimate:
\begin{eqnarray}\label{e:ee}
E(u|_\Sigma)&=&\frac{1}{2}\int_\Sigma \big|du-X_{H^\mathcal{D}}\otimes \alpha\big|^2\hbox{Vol}_\Sigma\notag\\
&=&\int_\Sigma u^*d\theta-u^*(dH^\mathcal{D})\wedge \alpha\notag\\
&=&\int_\Sigma d\big(u^*\theta-(u^*H^\mathcal{D})\alpha\big)+\beta(u)\notag\\
&\leq&\int_\Sigma d\big(u^*\theta-(u^*H^\mathcal{D})\alpha\big)\notag\\
&=&\int_{\partial\Sigma} u^*\theta-(u^*H^\mathcal{D})\alpha\notag\\
\end{eqnarray}

Since $\alpha|_{\partial \mathcal{D}}=0$, for any connected component $\varpi$ of $\partial_l\Sigma$ we have $(u^*H^\mathcal{D})\alpha|_\varpi=0$. And since $u^*\theta|_{\widehat{L}}=u^*dk_L$, by Stokes' theorem we get $\int_\varpi u^*\theta=0$ for circles $\varpi$, while for intervals $\varpi$, $\int_\varpi u^*\theta=k_L(u(p))-k_L(u(q))$ for corners $p,q\in\partial_b\Sigma\cap \partial_l\Sigma$. So this integral also vanishes thanks to the assumption that $k_L|_{\widehat{L}\setminus L}=0$.
It therefore follows from (\ref{e:ee}) and the contact condition $d\rho\circ J=-\theta$ outside $M$ that
\begin{eqnarray}\label{e:ee'}%\max_{0\leq i\leq 2}
E(u|_\Sigma)&\leq&\int_{\partial_b\Sigma} u^*\theta-(u^*H^\mathcal{D})\alpha\notag\\
&\leq&\int_{\partial_b\Sigma} \theta\circ\big(du-X_{H^\mathcal{D}}(u)\otimes \alpha\big)-a\alpha
\notag\\
&=&\int_{\partial_b\Sigma} \theta\circ\big(du-X_{H^\mathcal{D}}(u)\otimes \alpha\big)-a\int_{\Sigma}
d\alpha\notag\\
&\leq&\int_{\partial_b\Sigma} (\theta\circ J)\circ\big(du-X_{H^\mathcal{D}}(u)\otimes \alpha\big)\circ(-j)\notag\\
&=&\int_{\partial_b\Sigma} d\rho\circ\big(du-X_{H^\mathcal{D}}(u)\otimes \alpha\big)\circ(-j)\notag\\
&=&\int_{\partial_b\Sigma} d\rho\circ du\circ(-j),
\end{eqnarray}
where we have used $\theta(X_{H^\mathcal{D}})+a=H^\mathcal{D}$ outside $M$ in the second inequality, and $d\rho(X_H)=0$ in the last equality.
Let $\xi$ be a tangent vector to $\partial_b\Sigma$ which gives rise to the boundary orientation. Then $j\xi$ points into $\Sigma$, and thus $du(j\xi)$ does not point outwards along $\partial M$, so $d\rho\circ du(j\xi)\geq 0$. Then it follows from (\ref{e:ee'}) that $E(u|_\Sigma)=0$. This implies that each connected component of $u|_\Sigma$ is contained in a single orbit of $X_{H^\mathcal{D}}$. If $\partial_b \Sigma\neq\emptyset$, then since $X_{H^\mathcal{D}}$ is tangent to $\partial M$, this orbit must be contained in $\partial M$.
\end{proof}

Now we define the moduli space $\mathcal{M}(H^\mathcal{D},J^\mathcal{D},\alpha;z_0,z_1,z_2)$ to be the space of the solutions $u$ to (\ref{e:inhom}). We call the perturbation datum $(H^\mathcal{D},J^\mathcal{D})$ \emph{regular} if for all choices of chords $z_i$ and points $u\in\mathcal{M}(H^\mathcal{D},J^\mathcal{D},\alpha;z_0,z_1,z_2)$, the linearized operators $D_u$ are surjective. It can be shown that there is a residual set of compatible perturbation data $(H^\mathcal{D},J^\mathcal{D})$ such that all $u\in\mathcal{M}(H^\mathcal{D},J^\mathcal{D},\alpha;z_0,z_1,z_2)$  are regular. In this case, $\mathcal{M}(H^\mathcal{D},J^\mathcal{D},\alpha;z_0,z_1,z_2)$ is a smooth manifold, and its dimension equals to $\mu(z_2)-\mu(z_1)-\mu(z_0)$. When this dimension is zero, this manifold is
compact and hence consists of a finite number of points. Therefore we have a chain map
\[
*:CW^*(L,H^0)\otimes CW^*(L,H^1)\longrightarrow CW^*(L,2H^2),
\]
\[
z_0*z_1=\sum_{\mu(z_2)=\mu(z_0)+\mu(z_1)}\sharp_{\Z_2}\mathcal{M}(H^\mathcal{D},J^\mathcal{D},\alpha;z_0,z_1,z_2)z_2.
\]
This can be proved by the usual transversality and gluing arguments, combined with a $C^0$-estimate for Floer trajectories (see Lemma~\ref{lem:bound}). Moreover, by the standard cobordism arguments, one can show that the induced map on homology is independent of the choices of conformal structures on $\mathcal{D}$ and
perturbation data $(H^\mathcal{D},J^\mathcal{D})$, see~\cite[Section~6.12 and Section~16]{Rit}. So we have a well-defined bilinear
map on homology
\[
*_F:HW^*(L,H^0)\otimes HW^*(L,H^1)\longrightarrow HW^*(L,2H^2).
\]

In a local holomorphic coordinate $s+it$ for $(\mathcal{D},j)$, we have
\[
\frac{1}{2}|du-X\otimes \alpha|{\rm Vol}_\mathcal{D}=u^*d\theta-u^*dH^\mathcal{D}_z\wedge \alpha=d\big(u^*(\theta-H^\mathcal{D}\alpha)\big)+\beta(u)
\]
from which, applying the boundary condition $u(\partial \mathcal{D})\subset \widehat{L}$ and Stokes' theorem, we get
\begin{eqnarray}
0\leq E(u)&=&\mathcal{A}_{L,H^2}(z_2)-\mathcal{A}_{L,H^0}(z_0)-\mathcal{A}_{L,H^1}(z_1)+\int_\mathcal{D}\beta(u(z))\notag\\
&\leq&\mathcal{A}_{L,2H^2}(z_2)-\mathcal{A}_{L,H^0}(z_0)-\mathcal{A}_{L,H^1}(z_1)\notag
\end{eqnarray}
where in the second equality we have used the non-positive condition~(\ref{nonpos}).
So the product $*_F$ induces a map on the filtered wrapped Floer cohomology
\begin{equation}\label{eq:filt0}
*_F:HW^*_{(a_0,b_0]}(L,H^0)\otimes HW^*_{(a_1,b_1]}(L,H^1)\longrightarrow HW^*_{(a_2,b_2]}(L,2H^2).
\end{equation}
for any $a_i<b_i,i=0,1$ with $a_2=a_0+a_1$ and $b_2=\max\{a_0+b_1,a_1+b_0\}$.

\begin{figure}[H]
	\centering
\includegraphics[scale=0.6]{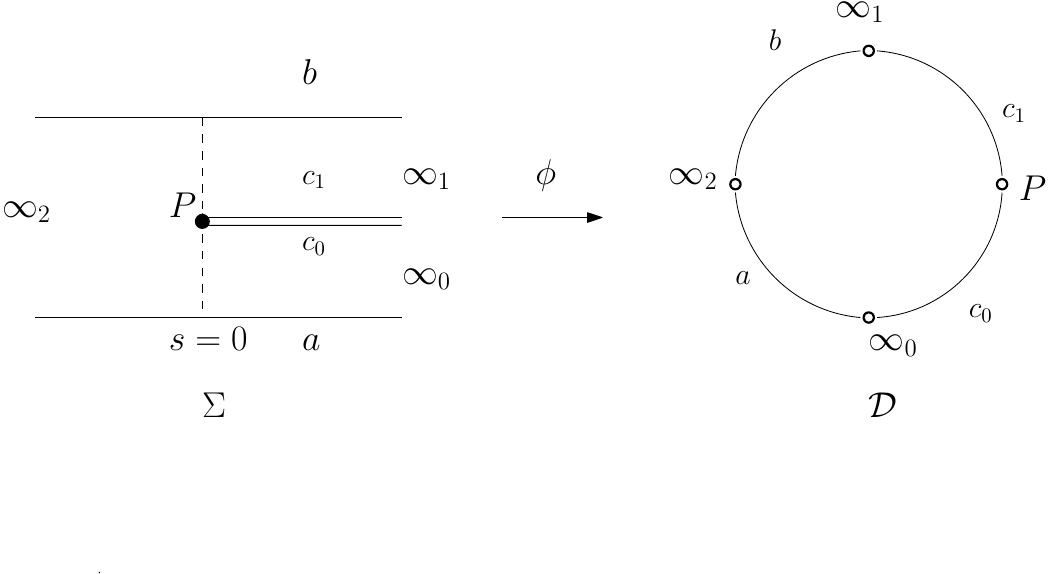}
\caption{The biholomorphic map $\phi$ between two surfaces.}\label{fig:hol}
\end{figure}

For our purpose, we consider the open subset
\[
\Sigma=\R\times(0,2)\setminus[0,\infty)\times \{1\}
\]
and equip it with the conformal structure such that one can map $\Sigma$ holomorphically onto the interior of $\mathcal{D}$. Besides, we require that, under this biholomorphism, the negative (positive) ends of $\Sigma$ correspond to the negative (positive) punctures
of $\mathcal{D}$, see~Figure~\ref{fig:hol}.

 Given $H,K\in\cH^{reg}$ which are linear outside $M$ and satisfy the following condition \[\frac{\partial^rH_t}{\partial t^r}\bigg|_{t=1}=\frac{\partial^r K_t}{\partial t^r}\bigg|_{t=0},\quad \hbox{for all}\;r\in\N,\] we define the concatenation $H*K$ as
\begin{equation}\label{conc}
(H*K)_t=
\begin{cases}
H_{t},& t\in[0,1],\\
K_{t-1} & t\in[1,2].
\end{cases}\quad
\end{equation}

Now we choose the following parameterized Hamiltonian
\begin{equation}\notag
H^\Sigma=
\begin{cases}
H_{t},& (s,t)\in [0,\infty)\times(0,1),\\
K_{t-1} & (s,t)\in [0,\infty)\times(1,2),\\
(H*K)_t & (s,t)\in (-\infty,0]\times(0,2),
\end{cases}
\end{equation}
and the family of complex structures $(J^\Sigma_{(s,t)})_{(s,t)\in\R\times[0,2]}$ which are of contact type outside $M$ and
\begin{equation}\notag J^\Sigma_{(s,t)}=\begin{cases}J_{t}^0,& (s,t)\in [1,\infty)\times(0,1),\\J_{t-1}^1 & (s,t)\in [1,\infty)\times(1,2),\\J_t^2 & (s,t)\in (-\infty,-1]\times(0,2),\end{cases}\end{equation}
where $J^i\in\mathcal{J}_\theta,i=0,1,2$ are complex structures such that $HW^*(L,H)$, $HW^*(L,K)$ and $HW^*(L,H*K)$ are well-defined, respectively. The aforementioned biholomorphism $\phi:\Sigma\to\mathcal{D}^\circ$ gives rise to the corresponding perturbation datum on the interior $\mathcal{D}^\circ$, and hence one on the whole disk by continuity, which we denote by $(H^\mathcal{D},J^\mathcal{D})$. In this case, if we take $\alpha=\phi_*dt$ then $d\alpha=0$ and the non-positive condition~(\ref{nonpos}) holds, obviously.
Moreover, for a generic choice of $J^\mathcal{D}$, the datum $(H^\mathcal{D},J^\mathcal{D})$ is regular. Hence, by~(\ref{eq:filt0}) we obtain the product
\begin{equation}\label{eq:filt}
*_F:HW_{(a,\infty)}^*(L,H)\otimes HW_{(b,\infty)}^*(L,K)\longrightarrow HW_{(a+b,\infty)}^*(L,H*K).
\end{equation}
for any $a,b\in\R\cup\{-\infty\}$.

\underline{}

\section{Spectral invariants from wrapped Floer theory}

\subsection{Morse cohomology}\label{subsubsec:morse}
\begin{df}\label{adapt}
Fix a Riemannian metric $g$ on $\widehat{L}$.
A function $f\in C^\infty(\widehat{L})$ is said to be \emph{adapted} to $L$ if the following conditions hold:
\begin{itemize}
  \item[(1)] No critical points of $f$ occur in $\widehat{L}\setminus L\cong \partial L\times(1,\infty)$;
  \item[(2)] The gradient vector field $\nabla f$ of $f$ with respect to $g$ points outward along $\partial L$;
  \item[(3)] $f|_L$ is a Morse function, and $(f,g)$ is a Morse-Smale pair.
\end{itemize}
\end{df}

Fix an adapted function $f$, and denote $\crit_k(f)$ the set of critical points $x$ of $f$ with Morse index $m_f(x)=k$. We define a differential $\delta$ on the free $\Z_2$-modules $CM^k(L,f,g)=\langle\crit_k(f)\rangle_{\Z_2}$  by counting isolated negative gradient flow lines of $\nabla f$, i.e.
\[\delta:CM^k(L,f,g)\longrightarrow CM^{k+1}(L,f,g),\]
\[\delta q=\sum_{m_f(p)=k+1}\sharp_2\mathcal{M}_{p,q}(f,g)\cdot p,\]
where the moduli space $\mathcal{M}_{p,q}(f,g)$ is given by
\[\mathcal{M}_{p,q}(f,g):=\big\{\gamma\in C^\infty(\R,L)\big|\dot{\gamma}=-\nabla f(\gamma(t));\;\gamma(-\infty)=p,\;\gamma(\infty)=q\big\}\big/\R.\]
The cohomology of the complex $(CM^k(L,f,g),\delta)$ is called the \emph{Morse cohomology}  for the pair $(f,g)$ which we denote by $HM^*(L,f,g)$. It can be shown that $HM^*(L,f,g)$ is isomorphic to the singular cohomology of $H^*(L)$ over $\Z/2$,
and it does not depend on the Morse-Smale pair $(f,g)$, see for instance~\cite{Sch2}.

\subsection{Lagrangian PSS morphism}\label{subsec:pss}

The Piunikhin-Salamon-Schwarz (PSS) homomorphism was firstly introduced to compare Hamiltonian Floer homology with singular homology, see~\cite{PSS}. After that, it was adapted to the Lagrangian setting for the case of cotangent bundles by Kati\'{c} and Milinkovi\'{c}~\cite{KM}, and later in more generality by Biran and Cornea~\cite{BC0}, Albers~\cite{Al}, Leclercq and Zapolsky~\cite{LZ}.
The Lagrangian PSS morphism is useful for us in the present paper since it respects the product structures between Morse cohomology and wrapped Floer cohomology, and plays an important role in the properties of spectral invariants. %see~Proposition~\ref{prop:spectral}.

%\subsubsection{Lagrangian PSS map}\label{subsubsec:pss}

%Let $H\in\cH_{<\tau}^{reg}$.  $$\cH_{<\tau}^{reg}=\big\{H\in\cH_{<\tau}\big|H|_M\;\hbox{is a $C^2$-small Morse function}\big\}.$$

\begin{figure}[H]
	\centering
\includegraphics[scale=0.7]{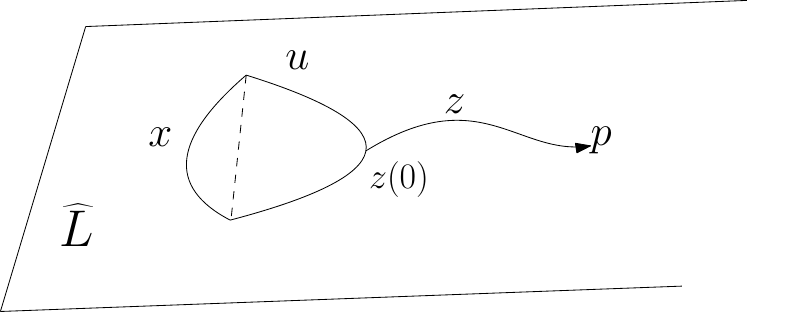}
	\caption{Counting isolated spiked disks defining the PSS-map.}\label{fig:pss}
\end{figure}

Let $f$ be a function adapted to $L$ and $H\in\cH$ a non-degenerate admissible Hamiltonian that is linear outside $M$.
Following~\cite[Section~2.2.3]{Le}, we now briefly recall the construction of the Piunikhin-Salamon-Schwarz homomorphism from $HM^*(L,f,g)$ to $HW^*(L,H)$ by counting isolated spiked Floer strips (see Figure~\ref{fig:pss}).

Let $\chi:\R\to [0,1]$ be a smooth cutoff function which satisfies: $\chi(s)=1$ if $s\leq 0$, $\chi(s)=0$ if $s\geq 1$, and $\chi'(s)\leq 0$ for all $s$. %where $R$ is a fixed positive number. %For simplicity  we set \[H_R(s,t,x)=\chi(s)H(t,x).\]

Let $p\in L$ be a critical point of $f$ and let $x\in\mathcal{C}(L,H)$.  Let $(J^s)_{s\in\R}\subset \cJ_{\theta}$ be a family of complex structures so that $J^s_t$ is constant outside $s\in[0,1]$ and is independent of $s,t$ for $s\geq 1$.

We consider the moduli space $\M(p,x;H,J,\chi,f,g)$ of the pairs of maps
\[z:[0,\infty)\longrightarrow L,\quad u:\R\times [0,1]\longrightarrow \widehat{M}\]
which satisfy
\begin{equation}\label{e:Mpss}
\begin{cases}
\frac{dz}{dt}=-\nabla f(z(t)),\\
\partial_s u+J_{t}\big(\partial_t u-\chi(s)X_{H}(u)\big)=0\;\hbox{with}\;E(u)<\infty,\\
z(\infty)=p,\; u(-\infty,t)=x(t),\\
u(s,0), u(s,1)\in L,\\
z(0)= u(+\infty).
\end{cases}
\end{equation}

Since both the critical points of $f$ and Hamiltonian chords of $H$ lie in the compact set $M$, and outside of it $H$ is linear and $\partial_s(\chi(s)H)\leq 0$, and $-\nabla f$ points inward along $\partial L$, by Lemma~\ref{lem:bd} every pair $(z,u)$ in the above space is entirely contained in $M$.

Up to generic choices of $(f,g,J)$, the space $\M(p,x;H,J,\chi,f,g)$ is a smooth manifold with dimension
$\mu(x)-m_f(p)$ (cf.~\cite{Al}), where $m_f(p)$ denotes the Morse index of $f$ at $p${\footnote{Note that the Floer strip u is holomorphic for $s\geq 1$ because of the choice of the cut-off function $\chi$. The finite energy condition for $u$ implies that $u$ can be extended continuously to the point $u(+\infty)$ in $L$, meaning that topologically $u$ is a half disk.}.

We define the morphism $\psi^H_f:CM^k(L,f,g)\to CW^k(L,H)$ on generators by
\[
\psi^H_f(p)=\sum_{\mu(x)=m_f(p)}\sharp_{\Z_2}\M(p,x;J,H,\chi,f,g)\cdot x.
\]
and extend this map by linearity over $\Z_2$. This extended map is a chain map and induces the PSS-type homomorphism
$$\psi_{f}^H:HM^k(L,f,g)\longrightarrow HW^k(L,H).$$

Recall that, for any $H\in\cH$, by definition $HW^*(L,H)=HW^*(L,\widetilde{H})$ for a sufficiently small regular perturbation $\widetilde{H}\in\cH^{reg}$ of $H$. Due to $H^*(L)\cong HM^*(L,f,g)$, we simplify the notation and denote the corresponding homomorphism from $H^*(L)$ to $HW^*(L,H)$ as $\psi_{pss}^H$.

\begin{lem}[{\cite[Lemma~2.3]{Le}}]\label{lem:iso}
If $H,K\in\cH$ with slopes $\tau_H,\tau_K$ such that $\tau_H\leq \tau_K$, then the following diagram commutes:
\begin{equation}\notag%
\xymatrix{
HW^*(L,H)
 \ar[rr]^{\Phi_{HK}} &&% arrow cannot occupy space
HW^*(L,K)
\\
&H^*(L)
\ar[ul]^{\psi_{pss}^H}
\ar[ur]_{\psi_{pss}^K}}
\end{equation}

\end{lem}

Moreover, by the standard compactness and gluing arguments~\cite[Section~5]{Rit}, one can show that the following diagram communicates:
\begin{equation}\label{TQFT}
\begin{split}
\xymatrix{
&
H^*(L)\otimes H^*(L)
\ar[r]^{\qquad\cup}
\ar[d]^{\psi^H_{pss}\otimes\psi^K_{pss}}
&
H^*(L)
\ar[d]^{\psi^{G}_{pss}} \\
&
HW^*(L,H)\otimes HW^*(L,K)
\ar[r]^{\;\qquad*_F}
&
HW^*(L,G)}
\end{split}
\end{equation}
where $H,K,G\in\cH$ are admissible Hamiltonians such that the product as in~(\ref{e:pd}) is well-defined.

%$H(1,\cdot)=K(0,\cdot)$ with all time derivatives so that the Hamiltonian $H*K$ as defined in~(\ref{conc}) is an admissible Hamiltonian.

Moreover, if $H\in\cH_{<\tau}^{reg}$ is a $C^2$-small Morse function in $M$, then all intersections $\varphi_{H}(\widehat{L})\cap \widehat{L}$ lie in $M\setminus\partial M$, and $\varphi_{H}(L)$ lies in a Weinstein neighborhood that is symplectically diffeomorphic to $T^*L$. Under this identification, $\varphi_{H}(L)$ is the graph of one form $d\zeta$ with $\zeta\in C^\infty(L)$, and
the generators of the complex $CW^*(L,H)$ correspond to the critical points of $\zeta$ on $L$, and the Floer strips $u$ for $H$ correspond to the negative gradient flow lines of $\zeta$ by $u\to u(s,0)$. In this case $HW^*(L,H)$ computes $H^*(L)$, and $\psi^H_f$ is indeed an isomorphism, see~\cite[Section~15.6]{Rit}. It follows from Lemma~\ref{lem:slope} and Lemma~\ref{lem:iso} that for any $H\in\cH_{<\tau}^{reg}$, we have the PSS-type isomorphism
\[\psi^H_f:HM^*(L,f,g)\stackrel{\cong}{\longrightarrow} HW^*(L,H).\]

\subsection{The wrapped Floer cohomology for compactly supported Hamiltonians}

Recall that the set $\cH_c(M)$ consists of Hamiltonians $H\in C^\infty([0,1]\times M)$ which satisfy $\cup_{t\in[0,1]}\supp(H_t)\subset M\setminus\partial M$. Without loss of generality we may assume that $H_t=H(t,\cdot)$ vanishes near $t=0,1$. Indeed, one can replace $H$ by $H'(t,\cdot)=\chi'(t)H(\chi(t),\cdot)$, where $\chi:[0,1]\to[0,1]$ is a monotone map with $\chi(t)=0$ near $t=0$ and $\chi(t)=1$ near $t=1$, then the Hamiltonian flow of $H'$ is a reparametrization of that of $H$, and in particular $\varphi_{H'}^1=\varphi_{H}^1$.

For each $H\in\cH_c(M)$, we take an admissible Hamiltonian $\widehat{H}\in \cH_{<\tau}^{reg}$ such that $\widehat{H}|_M$ is a $C^2$-small perturbation of $H$ and then define the wrapped Floer cohomology of $H$ as
\begin{equation}\label{cwfh}
HW^*(L,H):=HW^*(L,\widehat{H}).
\end{equation}

By Lemma~\ref{lem:slope}, the wrapped Floer cohomology $HW^*(L,\widehat{H})$ only depends on the slope of $\widehat{H}$ at infinity. So the above definition of wrapped Floer cohomology for $(L,H)$ is well-defined, i.e.  different choices of Hamiltonians $\widehat{H}$ yield isomorphic cohomologies $HW^*(L,\widehat{H})$.

\begin{rmk}
In the above definition of the wrapped Floer cohomology for a Hamiltonian $H\in\cH_c(M)$, we specify the range of an auxiliary Hamiltonian $\widehat{H}$ where $\widehat{H}$ is linear. Without this restriction there is no problem to define $HW^*(L,H)$ for any admissible Hamiltonian $H\in\cH$ with slope less than $\tau$ at infinity, but this point will be crucial for us to show the continuity property of the Lagrangian spectral invariant (to be defined later) for $H\in\cH_{<\tau}$ with respect to the Hofer distance, see Lemma~\ref{lem:continu}.
\end{rmk}

\subsection{The wrapped Floer capacity}\label{sec:wFc}
Let $\mathcal{R}^-(\partial L,\theta)$ be the symmetric set of $\mathcal{R}(\partial L,\theta)$ with respect to zero, i.e.
$$\mathcal{R}^-(\partial L,\theta):=\{-a|a\in\mathcal{R}(\partial L,\theta)\}.$$
Note that $HW^*_{(a,b]}(L)$ changes only when the action window crosses $\mathcal{R}^-(\partial L,\theta)$. This means that
\[
HW^*_{(a,b]}(L)\cong HW^*_{(c,d]}(L)\quad\hbox{if}\;(a,b]\cap \mathcal{R}^-(\partial L,\theta)=(c,d]\cap \mathcal{R}^-(\partial L,\theta).
\]
For $\delta\in(0,\tau)$, we have $(-\delta,\infty)\cap\mathcal{R}^-(\partial L,\theta)=\emptyset$, and
\[
\psi^f_{pss}:H^*(L)\stackrel{\cong}{\longrightarrow} HW^*(L,f)\cong HW^*_{(-\delta,\infty)}(L).
\]
where $f\in \cH_{<\tau}^{reg}$.

\begin{df}\label{df:wFc}
The \emph{wrapped Floer capacity} for $L$ is defined by
\[
c_{HW}(L)=\inf\big\{a>0|\iota_{-a}^L\circ\psi^f_{pss}(\mathds{1}_L)=0\big\}
\]
where $\mathds{1}_L\in H^0(L)$ is the fundamental class, and $\iota_{a}^L$ is the inclusion map (see~(\ref{diag:IQw})). And by convention we set
$c_{HW}(L)=\infty$ if $\iota_{-a}^L\circ\psi^f_{pss}(\mathds{1}_L)\neq0$ for all $a>0$.
\end{df}

Since $HW^*(L)$ admits a ring structure and the PSS-map preserves the ring structures of $H^*(L)$ and $HW^*(L,f)$, and so for $\iota_{-a}^L\circ\psi^f_{pss}$, the capacity $c_{HW}(L)$ is finite if and only if $HW^*(L)=0$.

%Given $\delta<0$, we take $H^\delta\in\cH_{<\tau}$ such $H^\delta$ is $C^2$-small Morse function on $M$ and its restriction to $L$ satisfies $$\delta<H^\delta|_L<0.$$ Using the PSS-map we have the isomorphism \[\psi_{pss}^{H^\delta}:H^*(L)\xrightarrow{\cong}HW^*(L,H^\delta)\]

\subsection{Spectral invariants from wrapped Floer cohomology}

%Let $H\in C^\infty([0,1]\times M)$ be a smooth function supported in $[0,1]\times M\setminus\partial M$.
We pick a non-degenerate admissible Hamiltonian $H\in \cH_{<\tau}^{reg}$. For $a\in\R$, we consider the following short exact sequence
\[
0\longrightarrow CW^*_{(a,\infty)}(L,H)\stackrel{\iota_a}{\longrightarrow} CW^*(L,H)\stackrel{\pi_a}{\longrightarrow} CW^*_{(-\infty,a]}(L,H)\longrightarrow 0.
\]
This induces the long exact sequence
\[
\longrightarrow HW^*_{(a,\infty)}(L,H)\stackrel{\iota_a}{\longrightarrow} HW^*(L,H)\stackrel{\pi_a}{\longrightarrow} HW^*_{(-\infty,a]}(L,H)\longrightarrow
\]

\begin{df}\label{SIV}
Let $0\neq\alpha\in H^*(L)$. Its associated \emph{Lagrangian spectral invariant} for $H$ is defined by
\[
\ell(H,\alpha):=\sup\big\{a\in\R\big|\pi_a\circ\psi_{pss}^H(\alpha)=0\big\}.
\]

\end{df}
Clearly, the above definition is equivalent to the following definition
\[
\ell(H,\alpha):=\sup\big\{a\in\R\big|\psi_{pss}^H(\alpha)\in\im(\iota_a)\big\}.
\]

\begin{rmk}
If $H$ is a non-degenerate admissible Hamiltonian with slope $\mu_H\notin\mathcal{R}(\partial L,\theta)$, in a similar way one can define the corresponding Lagrangian spectral invariant for $H$ which we will not pursue in this paper.
\end{rmk}

The following lemma implies that the spectral invariant $\ell$ defined as above depends continuously on the restrictions of Hamiltonians in $\cH_{<\tau}^{reg}$ to $M$ with respect to Hofer distance. Its proof is an adaption of that in~\cite[Section~3.2]{LZ} to the wrapped setting, and  similar arguments have already appeared in~\cite{Sch,Oh2,Le}.

\begin{lem}\label{lem:continu}
For $H,K\in \cH_{<\tau}^{reg}$ and non-zero class $\alpha\in H^*(L)$,
\[\int^1_0\min_M\big(H_t-K_t\big)dt\leq \ell(H,\alpha)-\ell(K,\alpha)\leq \int^1_0\max_M\big(H_t-K_t\big)dt.\]
\end{lem}

\begin{proof}
For an arbitrary $\epsilon>0$, we will prove the following inequality
\begin{equation}\label{e:aim}
\ell(H,\alpha)-\ell(K,\alpha)\geq\int^1_0\min_M\big(H_t-K_t\big)dt-\epsilon.
\end{equation}
Once this inequality is proved, by exchanging the role of $H$ and $K$ we get the inequality in the opposite direction.
To prove (\ref{e:aim}) we proceed in three steps:

\noindent \textbf{Step 1.}
Assume that $H^0,H^1\in \cH_{<\tau}^{reg}$ are the same on $M_r$ for some $r\in(0,1)$, and satisfy
\begin{equation}\label{e:fun}
H^i(t,\rho,x)=h^i(\rho)\quad \hbox{on}\;[0,1]\times[r,\infty)\times \partial M,\; i=0,1
\end{equation}
for two smooth functions $h^0,h^1\in C^\infty((0,\infty),\R)$ with $0\leq(h^0)'(\rho)\leq(h^1)'(\rho)<\tau$. We claim that  $\ell(H^0,\alpha)=\ell(H^1,\alpha)$.

%We pick the constant homotopy of Floer data $(H^s,J^s)$ with $H^s$
%Without loss of generality, we may assume that\[H^i(t,\rho,x)=\lambda_i\rho+c_i\quad \hbox{on}\;[0,1]\times[1,\infty)\times \partial M,\; i=0,1\] with $0\leq\lambda_0\leq\lambda_1<\tau$.
Let $s\mapsto\beta(s)$ be a smooth cutoff function on $\R$ such that $\beta(s)=0$ for $s\leq 0$, $\beta(s)=1$ for $s\geq 1$, and $\beta'(s)\geq 0$.
Consider the homotopy
\begin{equation}\notag
H^s=\beta(s)H^0+\big(1-\beta(s)\big)H^1
\end{equation}
which satisfies \[\frac{\partial^2 H^s}{\partial s\partial\rho}(t,\rho,x)=\beta'(s)\big((h^0)'(\rho)-(h^1)'(\rho)\big)\leq 0 \quad\hbox{on}\;[0,1]\times[r,\infty)\times \partial M.\]
It follows from Lemma~\ref{lem:bd} that solutions to the Floer equation (\ref{e:sfeq}) connecting $x_-\in\mathcal{C}(L,H^1)$ to $x_+\in\mathcal{C}(L,H^0)$ can not escape from $M_r$. Since the continuation map $\Phi_{H^0H^1}$ is independent of the choices of homotopies $(H^s,J^s)$ used to define it, and Hamiltonian chords of $H^0$ and $H^1$ are the same and lie in $M_r\subset M$, one can choose generic families of almost complex structures $J^s$ such that the continuation homomorphism $\Phi_{H^0H^1}$ induced by $(H^s,J^s)$ is the identity map on the chain level. Therefore, we have $\ell(H^0,\alpha)=\ell(H^1,\alpha)$.

\noindent \textbf{Step 2.} When $H^0,H^1\in \cH_{<\tau}^{reg}$ have the same slope, we will prove that (\ref{e:aim}) holds for these two Hamiltonians,
\begin{equation}\label{ineq}
\ell(H^1,\alpha)-\ell(H^0,\alpha)\geq\int^1_0\min_M\big(H^1_t-H^0_t\big)dt-\epsilon
\end{equation}
To this end we only need to prove that for each $\epsilon>0$ one can find a regular homotopy $(H^s,J^s)$ connecting $(H^0,J^0)$ to $(H^1,J^1)$ such that for all $a\in\R$, the corresponding continuation homomorphism
\[
\Phi_{H^0H^1}:\big(CW^*(L,H^0),d_{H^0,J^0}\big)\longrightarrow \big(CW^*(L,H^1),d_{H^1,J^1}\big)
\]
maps $CW^*_{(a,\infty)}(L,H^0)$ into $CW^*_{(a+b,\infty)}(L,H^1)$ for $b=\int^1_0\min_M\big(H^1_t-H^0_t\big)dt-\epsilon$. Then the following commutative diagram

\begin{equation}\notag%decomposition of homorphisms
\begin{split}
\xymatrix{H^*(L)
\ar[r]^{\psi^{H^0}_{pss}\quad}
\ar[dr]_{\psi^{H^1}_{pss}\;}
&
HW^*(L,H^0)
\ar[r]^{\pi_a\quad}
\ar[d]^{\Phi_{H^0H^1}}
&
HW^*_{(-\infty,a]}(L,H^0)
\ar[d]^{\Phi_{H^0H^1}} \\
&
HW^*(L,H^1)
\ar[r]^{\pi_{a+b}\quad}
&
HW^*_{(-\infty,a+b]}(L,H^1)}
\end{split}
\end{equation}
implies (\ref{ineq}).

As before, we first consider the special homotopy of Hamiltonians \[H^s=\beta(s)H^0+\big(1-\beta(s)\big)H^1.\]
Although in this case $H^s$ in general is not regular, one can pick a regular homotopy of Floer datum $(K,J)$ such that $K^s,H^s$ are the same on $[0,1]\times\widehat{M}\setminus M$ for all $s$, $K^s$ is independent of $s$ outside of the interval $[0,1]$, and
\[
\max_{(s,t,z)\in [0,1]\times[0,1]\times M}\bigg|\frac{\partial_s K^s_t}{\partial s}(z)-\frac{\partial_s H^s_t}{\partial s}(z)\bigg|\leq \epsilon,
\]
see~Proposition~\ref{prop:reg} for a proof of this fact. Here we may further require that $J^s\equiv J\in\cJ_\theta$ which is of contact type and time-independent outside of $M$.

Let $u$ be a solution to the Floer equation (\ref{e:sfeq}) for $(K,J)$ connecting $x\in\mathcal(L,H^1)$ to $y\in\mathcal(L,H^0)$. Since $H^0,H^1$ are linear functions in $\R$-coordinate outside of $M$ with the same slope, it follows from Lemma~\ref{lem:bd} that $u$ is contained in $M$. By the energy identity (\ref{E}), we have
\begin{eqnarray}
&&\mathcal{A}_{L,H^0}(y)-\mathcal{A}_{L,H^1}(x)=-E(u)+\int_{\R\times [0,1]}\big(\partial_sK^s_t\big)\big(u(s,t)\big)dsdt\notag\\
&\leq&\int_{\R\times [0,1]}\big(\partial_sH^s_t\big)\big(u(s,t)\big)dsdt+
\int_{\R\times [0,1]}\bigg(\frac{\partial_s K^s_t}{\partial s}-\frac{\partial_s H^s_t}{\partial s}\bigg)\big(u(s,t)\big)dsdt\notag\\
&\leq&\int^1_0\int^{\infty}_{-\infty}\beta'(s)\max_{z\in M}\big(H^0_t(z)-H^1_t(z)\big)dsdt+\epsilon\notag\\
&=&\int^1_0\max_{z\in M}\big(H^0_t(z)-H^1_t(z)\big)dt+\epsilon\notag
\end{eqnarray}
which implies
\[\mathcal{A}_{L,H^1}(x)\geq \mathcal{A}_{L,H^0}(y)+\int^1_0\min_{z\in M}\big(H^1_t(z)-H^0_t(z)\big)-\epsilon.\]

\begin{figure}[H]
	\centering
	\includegraphics[scale=0.5]{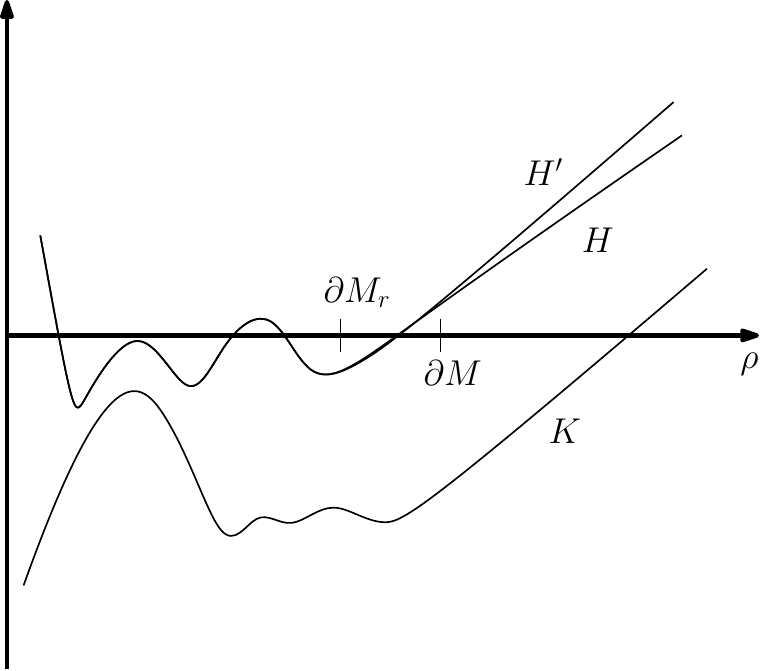}
	\caption{Modifying $H$ into $H'$ near $\partial M$ that has the same slope with $K$.}\label{fig:deformH}
\end{figure}

%due to the denseness $\cH_{<\tau}^{reg}$ in $\cH_{<\tau}$
\noindent \textbf{Step 3.} If $H,K\in \cH_{<\tau}^{reg}$,  we modify $H$ near $\partial M$ into
a regular Hamiltonian $H'\in \cH_{<\tau}^{reg}$ which satisfies
\begin{itemize}
  \item $H,H'$ are the same on $M_r$ for some $r\in(0,1)$ and satisfy (\ref{e:fun});
  \item $H', K$ have the same slope on $[1,\infty)\times \partial M$
    \item and the estimate
\begin{equation}\label{appr}
\max_{(t,z)\in[0,1]\times M}\big|H'_t(z)-H_t(z)\big|<\frac{\epsilon}{2}.
\end{equation}

\end{itemize}
\noindent We refer to Figure~\ref{fig:deformH} for a schematic graph of these Hamiltonians.

By \textbf{step 1},
\begin{equation}\label{1st}
\ell(H,\alpha)=\ell(H',\alpha).
\end{equation}
By \textbf{step 2},
\begin{equation}\label{2st}
\int^1_0\min_M\big(H'_t-K_t\big)dt-\frac{\epsilon}{2}\leq \ell(H',\alpha)-\ell(K,\alpha).%\leq \int^1_0\max_M\big(H'_t-K_t\big)dt+\frac{\epsilon}{2}.
\end{equation}
Combining (\ref{appr}), (\ref{1st}) and (\ref{2st}) yields (\ref{e:aim}). This completes the proof.
\end{proof}

Thanks to Lemma~\ref{lem:continu}, one can now $C^\infty$-continuously extend $\ell$ to a map $\ell:\cH_{<\tau}\to\R$ because $\cH_{<\tau}^{reg}$ is $C^\infty$-dense in $\cH_{<\tau}$. Using this extended map, one can define a map $$\ell:\mathcal{H}_c(M)\times H^*(L)\setminus\{0\}\to\R.$$ Explicitly, for $H\in\mathcal{H}_c(M)$ we take a sequence of Hamiltonians $(H^k)_{k\in\N}\subset\cH_{<\tau}^{reg}$ such that $H^k|_M\to H$ in $C^\infty$-topology as $k\to\infty$, then define
\begin{equation}\label{df:siv1}\ell(H,\al)=\lim_{k\to\infty}\ell(H^k,\al).\end{equation}

Notice that there is another way to define the spectral invariant for $H\in\cH_{<\tau}$. Indeed, we have already defined the wrapped Floer cohomology $HW^*(L,H)$ for all $H\in\mathcal{H}_c(M)$, see~(\ref{cwfh}), and for $a,b\notin\spec(L,H)$ the filtered wrapped Floer cohomology $HW^*_{(a,b]}(L,H)$ can be defined in a similar vain. Moreover, as in the regular case, we have the long exact sequence
\[
\longrightarrow HW^*_{(a,\infty)}(L,H)\stackrel{\iota_a}{\longrightarrow} HW^*(L,H)\stackrel{\pi_a}{\longrightarrow} HW^*_{(-\infty,a]}(L,H)\longrightarrow
\]
Hence, we may define
\begin{equation}\label{df:siv2}
\ell(H,\alpha)=\sup\big\{a\in\R\setminus\spec(L,H)\big|\pi_a\circ\psi_{pss}^H(\alpha)=0\big\}.
\end{equation}
It is not hard to verify that these two definitions of the spectral invariants for $H\in\mathcal{H}_c(M)$ coincide.

Clearly, for all $H,K\in \mathcal{H}_c(M)$ and non-zero class $\alpha\in H^*(L)$ we have
\begin{equation}\label{e:cont}
\int^1_0\min_M\big(H_t-K_t\big)dt\leq \ell(H,\alpha)-\ell(K,\alpha)\leq \int^1_0\max_M\big(H_t-K_t\big)dt.
\end{equation}

Since the wrapped Floer homology $HW^*(L,H)$ for $H\in\cH$ is precisely the Lagrangian Floer cohomology $HF^*_{\eta_0}(\varphi^1_H(\widehat{L}),\widehat{L})$, and the generators of $HW^*(L,H)$ for $H\in\cH_{<\tau}$ are not genuine wrapped chords, as in~\cite{Le}, we refer to the map $$\ell:\mathcal{H}_c(M)\times H^*(L)\setminus\{0\}\to\R$$ as the \emph{Lagrangian spectral invariant} for the pair $(L,H)$.

Recall that for $H,K\in\mathcal{H}_c(M)$, the composition map $\varphi_H^t\circ\varphi_K^t$ and the inverse map $(\varphi_H^t)^{-1}$ are generated  by
\[
H_t\sharp K_t=H_t+K_t\circ (\varphi_H^t)^{-1},\quad \overline{H}_t=-H_t\circ\varphi_H^t,
\]
respectively. Recall also that a smooth Hamiltonian $H_t$ on a symplectic manifold $(M,\omega)$ is called \emph{normalized} if $\int_MH_t\omega^n=0$ for all $t$.

In the following, we list some basic properties of $\ell$ that can be proved by adapting the methods used in~\cite{LZ} with minor changes.

\begin{prop}\label{prop:spectral}
	The map $\ell: \cH_c(M)\times H^*(L)\setminus\{0\} \to\R$ has the following properties:

\begin{itemize}
\item[(1)] \textbf{Continuity}: $\ell(H,\al)$ is Lipschitz in $H$ in the $C^0$-topology.
\item[(2)] \textbf{Spectrality}: $\ell(H,\al)\in\spec(L,H)$.
\item[(3)] \textbf{Normalization}: If $c$ is a function of time then
\[
\ell(H+c,\al)=\ell(H,\al)+\int^1_0c(t)dt.
\]
We define $\ell(\al,0)=0$ for all $\al\in H^*(L)\setminus\{0\}$.
\item[(4)] \textbf{Monotonicity}: $\ell(H,\al)\geq\ell(K,\al)$ for any $0\neq\al\in H_*(L)$ provided that $H(t,z)\geq K(t,z)$ for all $(t,z)\in[0,1]\times M$.
\item[(5)] \textbf{Homotopy invariance}: $\ell(H,\al)=\ell(K,\al)$, when $\varphi_H=\varphi_K$ in the universal covering of the group of Hamiltonian diffeomorphisms with compact support in $M\setminus\partial M$, where $H$ and $K$ are normalized.
\item[(6)] \textbf{Anti-triangle inequality}: $\ell(H\sharp K,\al\cup \beta)\geq \ell(H,\al)+\ell(K,\al)$.
\item[(7)] \textbf{Non-positivity}: $\ell(H,\mathds{1}_L)+\ell(\overline{H},\mathds{1}_L)\leq0$.
 \item[(8)] \textbf{Lagrangian control}: For all $H\in \cH_c(M)$ we have
	$$\int^1_0\min\limits_{L}H_tdt\leq \ell(\al,H)\leq\int^1_0\max\limits_{L}H_tdt.$$
 \item[(9)] \textbf{Symplectic invariance}: $\ell(H,\al)=\ell'(H\circ\phi^{-1},\phi^*(\al))$ for any symplectomorphism $\phi$ with compact support in $M\setminus\partial M$
     %(the group consisting of such maps is denoted by $\sym_c(M,\omega)$)
     satisfying $L'=\phi(L)$, where
	$\ell':H^*(L')\setminus\{0\}\times \cH_c(M)\to\R$
	is the corresponding spectral invariant.
\end{itemize}

\end{prop}

\begin{proof}
Property (1) follows from (\ref{e:cont}) immediately.

To prove property (2), arguing by contradiction we assume that $\ell_\alpha=\ell(H,\al)\notin\spec(L,H)$. Since $\spec(L,H)$ is a nowhere dense and closed subset of $\R$, there exists $\delta>0$ such that \[(\ell_\alpha-\delta,\ell_\alpha+\delta)\cap\spec(L,H)=\emptyset.\]
Hence, the quotient map
\[HW^*_{(-\infty,\ell_\al+\delta]}(L,H)\longrightarrow HW^*_{(-\infty,\ell_\al-\delta]}(L,H)\]
is an isomorphism. This implies $\pi_{\ell_\al+\delta}\circ\psi_{pss}^H(\alpha)=0$ because by definition of $\ell(H,\al)$ we have $\pi_{\ell_\al-\delta}\circ\psi_{pss}^H(\alpha)=0$. Therefore, by definition we have $\ell(H,\al)\geq\ell_\al+\delta$ which is absurd.

Property (3) follows from properties (1) and (2). Indeed, for the homotopy $H^s=H+sc$ with $s\in[0,1]$, the Hamiltonian chords are the same for all $s$, hence
\[\spec(L,sH)=\bigg\{a+s\int^1_0c(t)dt\big|a\in\spec(L,H)\bigg\}.\]
Then by property (2) we have $\ell(H+sc,\al)-s\int^1_0c(t)dt\in\spec(L,H)$. Since the action spectrum $\spec(L,H)$ is a closed nowhere dense set in $\R$, and due to property (1), $\ell(H+sc,\al)$ is continuous with respect to $s$, $\ell(H+sc,\al)-s\int^1_0c(t)dt$ must be constant. In particular, for $s=0,1$ we have
\[\ell(H+c,\al)-\int^1_0c(t)=\ell(H,\al)\]
which implies property (3).

Property (4) follows from the following commutative diagram

\begin{equation}\notag%
\begin{split}
\xymatrix{H^*(L)
\ar[r]^{\psi^{K}_{pss}\quad}
\ar[dr]_{\psi^{H}_{pss}}
&
HW^*(L,K)
\ar[r]^{\pi_a\quad}
\ar[d]^{\Phi_{KH}}
&
HW^*_{(-\infty,a]}(L,K)
\ar[d]^{\Phi_{KH}} \\
&
HW^*(L,H)
\ar[r]^{\pi_{a}\quad}
&
HW^*_{(-\infty,a]}(L,H)}
\end{split}
\end{equation}

Property (5) can be deduced from the following lemma:
\begin{lem}[{\cite[Lemma~31]{LZ}}]\label{lem:homotopy}
If there exists a homotopy $(H^s)_{s\in[0,1]}$ between two Hamiltonians $H$ and $K$ in $\cH_c(M)$ so that $\varphi_{H^s}^1=\varphi_H^1$ for all $s\in[0,1]$ and $\int^1_0dt\int_M(K_t(p)-H_t(p))\omega_p^n=0$, then $\spec(L,H)=\spec(L,K)$.
\end{lem}

We now show property (5) assuming Lemma~\ref{lem:homotopy}. Let $(H^s)_{s\in[0,1]}$ be a homotopy as in Lemma~\ref{lem:homotopy}. By property (2) and Lemma~\ref{lem:homotopy} we have
\[\ell(H^s,\al)\in\spec(L,H^s)=\spec(L,H)\quad\forall s\in[0,1].\]
It follows from property~(1) and the fact that $\spec(L,H)$ is a closed nowhere dense set in $\R$ that $\ell(H^s,\al)$ is independent of $s$. So we have $\ell(L,H)=\ell(L,K)$.

\noindent {\bf Proof of Lemma~\ref{lem:homotopy}.} Since $\varphi_{H^s}^1=\varphi_H^1$ for all $s\in[0,1]$, if $p\in(\varphi_H^1)^{-1}(L)\cap L$, then $x^s(t)=\varphi_{H^s}^t(p)$ is a Hamiltonian chord of $H^s$ for each $s$, i.e.    $x^s\in\crit(\mathcal{A}_{L,H^s})$. Hence,
\begin{eqnarray}\label{e:dif}
\mathcal{A}_{L,H^1}(x^1)-\mathcal{A}_{L,H^0}(x^0)&=&\int^1_0\frac{d}{ds}\mathcal{A}_{L,H^s}(x^s)ds\notag\\
&=&\int^1_0d\mathcal{A}_{L,H^s}(x^s)[\partial_sx^s]ds+\int^1_0\int^1_0\partial_sH^s_t(x^s(t))dsdt\notag\\
&=&\int^1_0\int^1_0\partial_sH^s_t(x^s(t))dsdt=: I.
\end{eqnarray}

Since $\varphi_H^1=\varphi_K^1$, the map
\[
\begin{split}
\crit(\mathcal{A}_{L,H})&\longrightarrow \crit(\mathcal{A}_{L,K})\\
\varphi_{H}^t(p)&\longmapsto\varphi_{K}^t(p)\quad\forall\;p\in(\varphi_H^1)^{-1}(L)\cap L
\end{split}
\]
is bijective. To finish the proof, it suffices to show the last expression in (\ref{e:dif}) vanishes.

Let $\lambda:[0,1]\to [0,1]$ be a smooth function such that $\lambda(t)=0$ if $0\leq t\leq 1/8$, and $\lambda(t)=1$ if $t\geq 1/4$.
Consider the homotopy $G^s$ given by the concatenation of $H^s$ and $H$ with respect to the time variable
\[
G^s(t,x)=
\begin{cases}\lambda'(t)H^s\big(\lambda(t),x\big),\quad  \ \ & t\in[0,\frac{1}{2}],\\
-\lambda'\big(1-t\big)H\big(\lambda(1-t),x\big), \quad \ \ & t\in[\frac{1}{2},1].
\end{cases}
\]

For each $s\in[0,1]$, $G^s$ generates the flow
\[\varphi_{G^s}^t=\begin{cases}\varphi_{H^s}^{\lambda(t)},\qquad  \ \ & t\in[0,\frac{1}{2}],\\
\varphi_H^{\lambda(1-t)}\circ\varphi_H^{-1}\circ \varphi_{H^s}^1=\varphi_H^{\lambda(1-t)}, \qquad \ \ & t\in[\frac{1}{2},1]
\end{cases}\]
which is a loop in $\ham_c(M,d\theta)$. Thus, fixing $s\in[0,1]$, for every $p\in M$, $y^s_p(t)=\varphi_{G^s}^t(p)$ is a critical point of the functional
\[
\mathcal{A}_{G^s}(x)=\int^1_0G^s(x(t))dt-\int x^*\theta.
\]
Clearly, $\mathcal{A}_{G^0}(y^0_p)=0$ since $H^0=H$. So we have
\begin{eqnarray}\label{e:p}
\mathcal{A}_{G^1}(y^1_p)&=&\mathcal{A}_{G^1}(y^1_p)-\mathcal{A}_{G^0}(y^0_p)\notag\\
&=&\int^1_0\frac{d}{ds}\mathcal{A}_{G^s}(y^s_p(t))ds\notag\\
&=&\int^1_0d\mathcal{A}_{G^s}(y^s_p)[\partial_s y^s_p(t)]dt+\int^1_0\int^1_0\partial_sG^s_t(y^s_p(t))dsdt\notag\\
&=&\int^1_0\int^1_0\partial_sH^s_t\big(\varphi_{H^s}^t(p)\big)dsdt
\end{eqnarray}
where in the last equality we have used $y^s_p\in\crit(\mathcal{A}_{G^s})$. Since the functional $\mathcal{A}_{G^1}$ must be constant along every connected critical submanifold and since the map $p\mapsto y_p^1$ is a smooth embedding from $M$ to $\crit(\mathcal{A}_{G^1})$, the last term in (\ref{e:p}) is independent of $p\in M$.

Now we show that the last expression in (\ref{e:dif}) vanishes. Indeed, putting $\omega=d\theta$ we have
\begin{eqnarray}\label{vanish}
I\cdot\int_M\omega^n&=&\int^1_0\int^1_0\partial_sH^s_t\big(\varphi_{H^s}^t(p)\big)dsdt\cdot\int_M\omega^n\notag\\
&=&\int^1_0\int^1_0\bigg(\int_M\partial_s H^s_t\big(\varphi_{H^s}^t(p)\big)\omega^n_p\bigg)dsdt\notag\\
&=&\int^1_0\int^1_0\bigg(\int_M\partial_s H^s_t(p)\omega^n_p\bigg)dsdt\notag\\
&=&\int^1_0dt\int_M\big(K_t(p)-H_t(p)\big)\omega_p^n=0\notag\\
\end{eqnarray}
where in the third equality we used the fact that $\varphi_{H^s}^t$ preserve the symplectic form $\omega$ for all $s,t$, and in the last equality we used the condition that $H,K$ have the same mean value. This completes the proof of Lemma~\ref{lem:homotopy}.

To prove property (6), we notice that (\ref{TQFT}) implies the following commutative diagram

\begin{equation}\notag%
\begin{split}
\xymatrix{
&
H^*(L)\otimes H^*(L)
\ar[r]^{\qquad\cup}
\ar[d]^{\psi^H_{pss}\otimes\psi^K_{pss}}
&
H^*(L)
\ar[d]^{\psi^{H* K}_{pss}} \\
&
HW^*(L,H)\otimes HW^*(L,K)
\ar[r]^{\;\qquad*_F}
&
HW^*(L,H*K)}
\end{split}
\end{equation}
It follows from~(\ref{eq:filt})
that $\ell(H* K,\al\cup \beta)\geq \ell(H,\al)+\ell(K,\al)${\footnote{Though $H*K$ depends on the time variable in $[0,2]$, the spectral invariants can be defined similarly by using Hamiltonians parametrized by a time interval which is not necessarily $[0, 1]$ and are invariant with respect to the time reparametrization.}}, where $H_t,K_t$ are supposed to be zero near $t=0,1$ by reparametrazing in time. Since $\varphi^t_{H\sharp K}$ and $\varphi_{H*K}^{2t}$ are homotopic relative to the ends at $t=0,1$ and $\int^1_0\int_M(H\sharp K)_t(d\theta)^ndt=2\int^1_0\int_M(H*K)_{2t}(d\theta)^ndt$,  it follows from Lemma~\ref{lem:homotopy} that the desired inequality holds.

Property (7) follows from property (3) and (6). Indeed,
\[
\ell(H,\mathds{1}_L)+\ell(\overline{H},\mathds{1}_L)\leq\ell(H\sharp \overline{H},\mathds{1}_L)=\ell(0,\mathds{1}_L)=0.
\]

%by property (1) it suffices to prove this property for non-degenerate Hamiltonians $H\in\cH_{<\tau}^{reg}(M)$.

To prove property (8), we first consider the case that the restriction of $H$ to $L$ is a function of time, i.e.  $H_t|_L=c(t)$. Let $H^s=sH, s\in[0,1]$. Then $H^s|_L=sc(t)$, hence the Hamiltonian chords with ends in $L$ are constant chords lying in $L$. So we have
\[\spec(L,H^s)=\bigg\{s\int^1_0c(t)dt\bigg\},\quad \forall s\in[0,1].\]
By property (2), $\ell(H,\al)=\int^1_0c(t)dt$. For $H\in\cH_c(M)$, we set
\[c(t)=\max_{x\in L}H(t,x)\]
and pick a function $K\in\cH_c(M)$ such that $K_t|_L=c(t)$ and $H\leq K$. Then $\ell(K,\al)=\int^1_0c(t)dt$.
It follows from property (4) that
\[
\ell(H,\al)\leq \ell(K,\al)=\int^1_0\max_{x\in L}H(t,x)dt.
\]
The proof of the other direction of the inequality in property (8) is similar.

Before giving a sketch of the proof of property (9), we notice that whenever $\phi$ is a symplectomorphism with compact support in $M\setminus\partial M$,  $L'=\phi(L)$ is an admissible Lagrangian submanifold of $(M,d\theta)$ provided that it is so for $L$. Set $H'=H\circ\phi^{-1}$. By property (1), we only need to prove property (9) for the case that $H$ is a restriction of a regular Hamiltonian $\widehat{H}\in\cH^{reg}_{<\tau}$ to $M$.
It is easy to see that $\phi$ induces a canonical isomorphism between the wrapped Floer cohomologies
\[
\phi^*:HW^*(L',H')\longrightarrow HW^*(L,H).
\]
Since $\phi$ preserves the actions of the critical points of the respective action functionals, the standard cobordism argument yields
 the following commutative diagram

\begin{equation}\notag%decomposition of homorphisms
\begin{split}
\xymatrix{H^*(L')
\ar[r]^{\psi_{pss}^{H'}\quad}
\ar[d]^{\phi^*}
&
HW^*(L',H')
\ar[r]^{\pi_a\quad}
\ar[d]^{\phi^*}
&
HW^*_{(-\infty,a]}(L',H')
\ar[d]^{\phi^*} \\
H^*(L)
\ar[r]^{\psi_{pss}^{H}\quad}
&
HW^*(L,H)
\ar[r]^{\pi_{a}\quad}
&
HW^*_{(-\infty,a]}(L,H)}
\end{split}
\end{equation}
which implies property (9).
\end{proof}

\subsection{The Lagrangian spectral metric}
Recall that for any symplectic manifold $(M,\omega)$, the Hofer's bi-invariant metric on the Hamiltonian diffeomorphism group $\ham_c(M,\omega)$ is defined by
\[d_H(\varphi,id)=\inf\big\{\|H\|\big|\varphi=\varphi_H^1\big\},\quad d_H(\varphi,\psi)=d_H(\psi^{-1}\circ\varphi,id)\]
where
\[
\|H\|=\int^1_0\big(\sup_{x\in M}H(t,x)-\inf_{x\in M}H(t,x)\big)dt.
\]

Let $L$ be a Lagrangian submanifold of $(M,\omega)$. Recall that
$\cL(L)=\{\varphi(L)|\varphi\in\ham_c(M,\omega)\}$, and that a pseudo-metric on $\cL(L)$ is given by
$$\delta_H(L_1,L_2)=\inf\big\{\|H\|\big|\varphi^1_H(L_1)=L_2,\;H\in C_c^\infty([0,1]\times M)\big\}.$$

For our purpose, for a Lagrangian submanifold $L\subseteq M$ we also use the following pseudo-norm
\[
\|H\|_L=\int^1_0\big(\sup_{x\in L}H(t,x)-\inf_{x\in L}H(t,x)\big)dt.
\]
Note that the pseudo-metric $\delta_H$ on $\cL(L)$ can be also defined by using $\|H\|_L$, and this equivalence is obtained in~\cite{Us4}.

Following~\cite{Vi,Oh3,KS}, for an admissible Lagrangian submanifold $L^n\subset (M^{2n},\omega)$, we define the \emph{Lagrangian spectral pseudo-norm} by
$$\ga(L,H)=-\ell(H,\mathds{1}_L)-\ell(\overline{H},\mathds{1}_L),\quad \forall H\in\cH_c(M).$$

By property (7) in Proposition~\ref{prop:spectral},  $\ga(L,H)$ is a non-negative function on $\cH_c(M)$.
This gives rise to a function $\ga:\cL(L)\times \cL(L)\to [0,\infty)$, which is called the \emph{Lagrangian spectral pseudo-metric}, by setting
$$\ga(L_1,L_2)=\inf\ga(L,H),\quad \forall L_1, L_2\in\cL(L),$$
with $L_1=\varphi_F^1(L),L_2=\varphi_G^1(L)$, where the infimum runs over ${H\in\cH_c(M)}$ such that $\varphi_H^1(L)=\varphi_F^{-1}\circ\varphi_G^1(L)$.

\begin{lem}\label{lem:smallham}
Let $f:L\to\R$ be a smooth $C^2$-small function with compact support in $L\setminus \partial L$, and let $H\in\cH_c(M)$ be an autonomous Hamiltonian which coincides with the lift of $f$ to a Weinstein neighborhood of $L$ in $M$. More precisely, identifying $T^*L$ with some Weinstein neighborhood of $L$ in $M$, we let $H=f\circ\pi$ on a co-disk bundle $D_RT^*L\subset T^*L$ of radius $R>0$ containing $L^f:=\{(q,\partial_q f(a))\in T^*L|q\in L\}$, and $H=0$ outside $T^*_{R+1}L$ in $M$, where $\pi:T^*L\to L$ is the natural projection map.
Then $\gamma(L,H)=\|H\|_L$.
\end{lem}
\begin{proof}

It suffices to prove the above assertion for a $C^2$-small Morse function $f:\widehat{L}\to\R$ adapted to $L$ and the corresponding lift $H\in\cH_\tau^{reg}$. Clearly, the PSS-map
\[\psi_{pss}^H:H^*(L)\cong HM^*(L,f,g)\longrightarrow HW^*(L,H)\]
is a chain-level isomorphism sending critical points $x\in\crit(f)$ to the corresponding constant chords $x\in\hc$.
Here $g$ is a metric on $L$ such that $(f,g)$ is a Morse-Smale pair. Under the above map, $\mathds{1}_L\in H^*(L)$ corresponds to
\[\bigg[\sum^k_{i=1}x_i\bigg]\in HW^*(L,H)\]
where $x_i$ are critical points of $f=H|_L$ with Morse index zero. So we have
\[
\ell(H,\mathds{1}_L)=\min\big\{H(x_i)|i=1,\ldots,k\big\}=\min_{x\in L} H(x).
\]
Similarly, we get
\[
\ell(\overline{H},\mathds{1}_L)=\min_{x\in L} \overline{H}(x)=-\max_{x\in L} H(x).
\]
Therefore, $\gamma(L,H)=\|H\|_L$ and hence the continuity property of $\ell$ concludes the desired result.
\end{proof}

\begin{thm}\label{thm:nontrivality}
The pseudo-metric $\ga$ on $\cL(L)$ satisfies
\begin{itemize}
\item[(a)] $\ga(L_1,L_2)=0$ if and only if $L_1=L_2$;
\item[(b)] $\ga(L_1,L_2)=\ga(L_2,L_1)$;
\item[(c)] $\ga(\varphi(L_1),\varphi(L_2))=\ga(L_1,L_2)$ for all $\varphi\in\ham_c(M,d\theta)$;
\item[(d)] $\ga(L_1,L_2)\leq\ga(L_1,L_3)+\ga(L_2,L_3)$;
\item[(e)] $\ga(L_1,L_2)\leq \delta_H(L_1,L_2)$.
%\item[(e)] $\ga(L_1,L_2)=\ga'(\phi(L_1),\phi(L_2))$ for all $\phi\in\sym_c(M,d\theta)$, where $\ga':\cL(\phi(L))\times \cL(\phi(L))\to [0,\infty)$ is the corresponding pseudo-metric for $\phi(L)$.
\end{itemize}

\end{thm}

For proving property~(a) in Theorem~\ref{thm:nontrivality} (which is Theorem~\ref{nontrivality} of the introduction), we need a Lagrangian version of Schwarz's result~\cite[Proposition~5.1]{Sch}.
\begin{lem}\label{lem:displace}
Let $K\in\cH_c(M)$ and let $U\subseteq M\setminus\partial M$ be a subset such that $\cup_{t\in[0,1]}{\rm supp}(K_t)\subseteq U$. Suppose that the Hamiltonian time-one map of $H\in\cH_c(M)$ displaces $L$ from $U$, i.e. $\varphi_H(L)\cap U=\emptyset$.
Then we have $\gamma(L,K)\leq 2\gamma(L,H)$.

\end{lem}

We shall prove Theorem~\ref{thm:nontrivality} assuming Lemma~\ref{lem:displace}.

\begin{proof}
By definition, property~(a) is equivalent to the following statement:
\begin{equation}\notag
\ga(L,\varphi_H(L))>0\Longleftrightarrow L\neq\varphi_H(L)\quad \hbox{for}\; H\in\cH_c(M).
\end{equation}
Assume that $\varphi_H(L)\neq L$, then one can find an open subset $U$ of $M\setminus\partial M$ such that $U\cap L\neq\emptyset$ and $\varphi_H(L)\cap U=\emptyset$. Pick a $C^2$-small function $f\in C^\infty(L)$ which is supported in $U\cap L$, and let
$K\in\cH_c(M)$ be the corresponding lift of $f$ so that $\|K\|_L>0$. It follows from Lemma~\ref{lem:smallham} and Lemma~\ref{lem:displace} that $2\gamma(L,\varphi_H(L))\geq \gamma(L,K)=\|K\|_L>0$. This completes the proof of property~(a).

The symmetry property~(b) and the $\ham$-invariant property (c) follow from the definition of $\ga$. The triangle inequality~(d) follows from property (6) in Proposition~\ref{prop:spectral}. The inequality~(e) is implied by property (8) in Proposition~\ref{prop:spectral} (or by (\ref{e:cont})). %The bi-invariant property (e) follows from property (9) in Proposition~\ref{prop:spectral} and the fact that for $\varphi_H\in\ham_c(M,\theta)$ and $\phi\in\sym_c(M,\theta)$ it holds that\[\phi\circ\varphi_H^t\circ\phi^{-1}=\varphi_{H^\phi}^t\quad\forall t\]where $H^\phi(t,z)=H(t,\phi^{-1}(z))$.
\end{proof}

\begin{proof}[{The proof of Lemma~\ref{lem:displace}}]
By our assumption concerning $K$, we have $\varphi_K^t(U)=U$ and $\varphi_K^t=id$ outside of $U$. Since $\varphi_H$ displaces $L$ from $U$,
we have
\[
\varphi_K\circ \varphi_H(L)\cap L=\varphi_H(L)\cap L\subseteq M\setminus U.
\]
As a consequence, there is a bijection between the chords of $H$ and that of $K\sharp H$ by sending $x\in\hc$ to $y\in\mathcal{C}(L,K\sharp H)$ with $y(t)=\varphi_K^t(x(t))$. Let $\overline{x}$ be a capping of $x$.
We put $u(s,t)=\varphi_K^{st}(x(t))$ where $(s,t)\in [0,1]\times[0,1]$. Then $\overline{x}\sharp u$ is a capping of $y$ denoted by $\overline{y}$. Using (\ref{func}), a direct calculation shows that
$\mathcal{A}_{L,K\sharp H}(y)=\mathcal{A}_{L,H}(x)$, see~\cite[Lemma~18]{Go2} for a completely same computation.
Hence, $\spec(L,K\sharp H)=\spec(L,H)$. Now we pick a parameter $\epsilon\geq0$ and consider $\epsilon K$ in place of $K$.
By a similar argument as above we obtain $\spec(L,\epsilon K\sharp H)=\spec(L,H)$. It follows from the spectrality property of $\ell$ that
\[
\ell(\epsilon K\sharp H,\mathds{1}_L)\in\spec(L,\epsilon K\sharp H).
\]
Since $\spec(L,H)$ is a nowhere dense and closed subset of $\R$, the continuity property of $\ell$ implies that
$\ell(\epsilon K\sharp H,\mathds{1}_L)$ does not depend on $\epsilon$. In particular, we get
$\ell(K\sharp H,\mathds{1}_L)=\ell(H,\mathds{1}_L)$.

Note that $\overline{K}$ is also supported in $U$, it holds that
\begin{equation}\label{inveq}
\ell(\overline{H}\sharp \overline{K},\mathds{1}_L)=\ell(\overline{H},\mathds{1}_L).
\end{equation}
Indeed, we have\[\varphi_{\overline{H}}\circ\varphi_{\overline{K}}(L)\cap L=\varphi_{\overline{H}}(L)\cap L,\] and 
if $\varphi_{\overline{H}}(q)$ belongs to the intersection from above, then $q\notin U$. Moreover, there is a bijective map between $\mathcal{C}(L,\overline{H})$ and $\mathcal{C}(L,\overline{H}\sharp \overline{K})$ by sending $x\in\mathcal{C}(L,\overline{H})$ to $y\in\mathcal{C}(L,\overline{H}\sharp \overline{K})$ with $y(t)=\varphi_{\overline{H}}^t\circ\varphi_{\overline{K}}^t(x(0))$.
Then an analogous argument as above yields (\ref{inveq}).
So we have
\begin{eqnarray}
\gamma(L,K\sharp H)&=&\ell(K\sharp H,\mathds{1}_L)+\ell(\overline{K\sharp H},\mathds{1}_L)\notag\\
&=&\ell(K\sharp H,\mathds{1}_L)+\ell(\overline{H}\sharp \overline{K},\mathds{1}_L)\notag\\
&=&\ell(H,\mathds{1}_L)+\ell(\overline{H},\mathds{1}_L)
=\gamma(L,H)\notag
\end{eqnarray}
It follows from the anti-triangle property of $\ell$ and the definition of $\gamma$ that
\[\gamma(L,K)=\gamma(L,K\sharp H\sharp \overline{H})\leq \gamma(L,K\sharp H)+\gamma(L,\overline{H})=2\gamma(L,H).\]
\end{proof}

The following statement implies that the Lagrangian spectral invariant for $H\in \cH_{<\tau}$ depends only on the Lagrangian submanifolds $\widehat{L}$ and $\varphi^{-1}_H(\widehat{L})$ and not on $H$. In the weakly exact case for a closed Lagrangian submanifold $L$ (i.e. $\omega|_{\pi_2(M,L)}=\mu|_{\pi_2(M,L)=0}$), this fact was shown in~\cite{Le}. The proof given here is essentially due to Usher in~\cite[Proposition~6.2]{Us3}.

\begin{prop}\label{prop:shift}
Assume that $G,H\in\cH_{<\tau}^{reg}$ satisfy that $\varphi^{-1}_G(\widehat{L})=\varphi^{-1}_H(\widehat{L})$ and $G-H\in\cH_c(M)$. Then for suitably generic paths of almost complex structures $J,J'\in\mathcal{J}^{reg}$, there is a shift isomorphism $$\Psi_L:(CW^*(L,G),d_{G,J})\to (CW^*(L,H),d_{H,J'})$$ which preserves the filtration and grading, i.e. for any $k\in\Z$ and $a\in\R$,
\[
\Psi_L \big(CW^{k}_{(a,\infty)}(L,G)\big)= CW^{k}_{(a,\infty)}(L,H).
\]
Moreover, the following diagram
\begin{equation}\notag%
\xymatrix{
CW^*(L,G)
 \ar[rr]^{\Psi_L} &&% arrow cannot occupy space
CW^*(L,H)
\\
&CM^*(L,f,g)
\ar[ul]^{\psi_{f}^G}
\ar[ur]_{\psi_{f}^H}}
\end{equation}
commutes in cohomology.

\end{prop}

\begin{proof}
Let $\psi_t=\varphi^t_H\circ(\varphi^t_G)^{-1}$. Since $G-H\in\cH_c(M)$, the isotopy $\{\psi_t\}$ is supported in $M\setminus\partial M$.
The condition $\varphi^{-1}_G(\widehat{L})=\varphi^{-1}_H(\widehat{L})$ implies that $\psi_1(\widehat{L})=\widehat{L}$. Let $K\in \cH_c(M)$ be the Hamiltonian which generates $\{\psi_t\}$. Then we have $H=K\sharp G$, or equivalently
\[
G(t,z)=(H-K)(t,\psi_t(z))\quad \forall (t,z)\in [0,1]\times \widehat{M}.
\]
Let $x\in \cP(\hL,\hL)$ be any path which represents the point class $\eta_0$ in $\pi_0(\cP(\hL,\hL))$. Then $\Psi x$ defined by $(\Psi x)(t)=\psi_t(x(t))$ represents some class $\eta\in\pi_0(\cP(\hL,\hL))$. Moreover, $x\in\mathcal{O}_{\eta_0}(L,G)$ if and only if $\Psi x\in \mathcal{O}_{\eta}(L,H)$. Therefore, one can define a map $\Psi_L$ from $CF^*_{\eta_0}(\widehat{L},H)$ to $CF^*_\eta(\widehat{L},H)$ by sending a generator $x$ to $\Psi x$.

We take $J\in\mathcal{J}_\theta$ and consider the perturbed Cauchy-Riemann operator
\[
\overline{\partial}_{G,J}(u):=\partial_su+J_t(\partial_tu-X_G),
\]
where $u:\R\times [0,1]\to \widehat{M}$. Since the isotopy $\{\psi_t\}$ is the identity map on $\widehat{M}\setminus M$, the path $J'$ of almost complex structures  defined by $J_t'=(\psi_t)_*J_t(\psi_t^{-1})_*$ also belongs to $\mathcal{J}_\theta$. For any $(s,t)\in\R\times [0,1]$, a direct computation shows that
\[
(\psi_t)_*\big((\overline{\partial}_{G,J}(u))_{(s,t)}\big)=(\overline{\partial}_{H,J'}(\Psi u)\big)_{(s,t)},
\]
where $(\Psi u)(s,t)=\psi_t(u(s,t))$. Note that $u(\R\times\{i\})\subset \widehat{L}$ if and only if $(\Psi u)(\R\times\{i\})\subset \widehat{L}$ because $\psi_0$ is the identity map and $\psi_1(\widehat{L})=\widehat{L}$. Consequently, there exists one to one correspondence between the moduli spaces of flowlines defining the differential
$d_{G,J}$ on $CF^*_{\eta_0}(\widehat{L},G)$ and the ones defining $d_{H,J'}$ on $CF^*_\eta(\widehat{L},H)$. Moreover, this correspondence retains Fredholm regularity of those Floer strips. So if $J$ is chosen to be regular so that $d_{G,J}$ is well-defined, then $\Psi_L$ is an isomorphism of chain complexes. Then we get $HF^*_\eta(\widehat{L},H)\neq 0$ due to the fact that $HF^*_\eta(\widehat{L},H)\cong HF^*_{\eta_0}(\widehat{L},H)\cong H^*(L)\neq 0$ by the PSS-type isomorphism for $H\in\cH_{<\tau}$.

We claim that $\eta=\eta_0$. Indeed, by the continuation maps in Floer theory, we have that $HF^*_\eta(\widehat{L},H)\cong HF^*_\eta(\widehat{L},H')$ where $H'\in\cH_{<\tau}^{reg}$ is some Hamiltonian with the property that $H'|_M$ be $C^2$-small. Notice that all Hamiltonian paths for $H'$ with endpoints in $\widehat{L}$ will be contained within a Darboux-Weinstein neighborhood of $L$ and hence will be homotopic rel endpoints to a path entirely contained in $L$. This implies that $HF^*_\eta(\widehat{L},H')=0$ whenever $\eta\neq\eta_0$ since generators in the class $\eta$
do not exist, in contradiction to the non-vanishing of $HF^*_\eta(\widehat{L},H)$ shown as above. In other words, we get the isomorphism $\Psi_L:(CW^*(L,G),d_{G,J})\to (CW^*(L,H),d_{H,J'})$.

Now we turn to the grading issue. We will show that $\Psi_L$ preserves the gradings of two complexes. Recall that in each component $\eta$ of $\cP(\hL,\hL)$ we have fixed a reference path $\gamma_\eta$ and a symplectic trivialization $\mathfrak{t}_\eta$. We denote $\gamma_0:=\gamma_{\eta_0}$ (which is equal to $pt$ in $L$) and take a homotopy $w:[0,1]\times[0,1]\to \widehat{M}$ such that $w([0,1]\times\{i\})\subset \widehat{L}$ for $i=0,1$, and $w(0,\cdot)=\gamma_0$ and $w(1,\cdot)=\Psi(\gamma_0)$. This homotopy exists due to the previous paragraph. Fix a symplectic trivialization $\mathfrak{t}_w:w^*T\widehat{M}\to [0,1]^2\times\R^{2n}$ which restricts to the fixed trivialization $\mathfrak{t}_0$ over $\gamma_0$ and which maps $T^*_{w(s,i)}\widehat{L}$ to $\R^n\times\{0\}$ for all $s\in[0,1]$ and $i=0,1$. For $x\in\mathcal{O}_{\eta_0}(L,G)$, using $(\psi_t)_*$ we push forward the trivialization $\mathfrak{t}_v$ over the homotopy $v:[0,1]\times[0,1]\to \widehat{M}$ from $\gamma_0$ to $x$ used to define the Maslov index $\mu(x)$ to obtain a trivialization $\Psi_*\mathfrak{t}_v$ over the homotopy $\Psi v$ restricting to the trivialization $\Psi_*\mathfrak{t}_0$ over $\Psi\gamma_0$. Taking into account the difference of the trivializations $\Psi_*\mathfrak{t}_v|_{\Psi\gamma_0}$ and $\mathfrak{t}_w|_{\Psi\gamma_0}$ over $\Psi\gamma_0$, one can obtain a trivialization over the concatenation of
the maps $w$ and $\Psi v$ to compute $\mu(\Psi x)$ as follows. Let $\mathcal{L}_w=(\mathfrak{t}_w)^{-1}(\{1\}\times[0,1]\times\R^n\times\{0\})$ be the Lagrangian subbundle of $(\Psi\gamma_0)^*T\widehat{M}$. Then $\Psi_*\mathfrak{t}_0(\mathcal{L}_w)$ is a loop of Lagrangian subspaces of $\R^{2n}$ since $\mathcal{L}_w$ coincides with $T_{\Psi\gamma_0(i)}\widehat{L}$ over $i=0,1$. By the catenation and homotopy properties of the Robbin-Salamon-Maslov index, one can show that $\mu(\Psi x)-\mu(x)$ coincides with the Maslov index of this loop, and hence is independent of the choice of $x\in\mathcal{O}_{\eta_0}(L,G)$. Therefore, there is an integer $\nu\in\Z$ such that for any $k\in\Z$ we have the isomorphism
$$HW^{k}(L,G)\longrightarrow HW^{k+\nu}(L,H).$$
induced by $\Psi_L$. Notice that the gradings of $HW^*(L,G)$ and $HW^*(L,H)$ are concentrated in $0\leq *\leq n$ and both cohomologies (being isomorphic to $H^*(L)$) do not vanish, we obtain $\nu=0$.

Next we consider the action shift.  For any $x\in\mathcal{O}_{\eta_0}(L,G)$, concatenating $w$ and $\Psi v$ yields a homotopy from $\gamma_0$ to $\Psi x$. A routine computation shows that
\[\mathcal{A}_{L,H}(\Psi x)=\mathcal{A}_{L,G}(x)+\mathcal{A}_{L,K}(\Psi\gamma_0),\]
see for instance~\cite[Proposition~6.2]{Us3} or~\cite[Proposition~31]{KS}. So  $\Psi_L$ shifts the filtrations of these two complexes by the constant $\kappa=\mathcal{A}_{L,K}(\Psi\gamma_0)$. Since $\psi_1(\widehat{L})=\widehat{L}$, each $q\in \widehat{L}$ gives a Hamiltonian chord $\Psi q\in\crit(\mathcal{A}_{L,K})$. Hence $\mathcal{A}_{L,K}\circ\Psi$ is constant along the connected Lagrangian $\widehat{L}$. Note that if $q\in \widehat{L}\setminus \supp (K)$ then $\mathcal{A}_{L,K}(\Psi q)=0$. Therefore $\kappa=0$.

As for the commutative diagram, the isomorphism induced by $\Psi_L$ in cohomology indeed coincides with the continuation map $\Phi_{GH}$, see for instance~\cite[Remark~5.5]{HLL}.
\end{proof}

\begin{cor}\label{cor:indep}
If $G,H\in\cH_c(M)$ satisfy that $\varphi^{-1}_G(L)=\varphi^{-1}_H(L)$, then for any non-zero class $\alpha\in H^*(L)$, we have $\ell(H,\alpha)=\ell(G,\alpha)$.

\end{cor}% and $\int_MG_t\omega^n=\int_MH_t\omega^n$ for all $t\in[0,1]$

\begin{proof}
Pick $F\in\cH_{<\tau}$ with the properties that $F|_M$ is $C^2$-small and $\varphi_G\circ\varphi_F(\widehat{L})$ intersects $\widehat{L}$ transversally. Then we have $G\sharp F, H\sharp F\in\cH_{<\tau}^{reg}$ and
$\varphi^{-1}_{G\sharp F}(\widehat{L})=\varphi^{-1}_{H\sharp F}(\widehat{L})$ since $\varphi^{-1}_G(L)=\varphi^{-1}_H(L)$.
It follows from Proposition~\ref{prop:shift} that $\ell(H\sharp F,\alpha)=\ell(G\sharp F,\alpha)$ for any non-zero class $\alpha\in H^*(L)$, where
$K\in\cH_c(M)$ is the Hamiltonian generating the isotopy $\varphi^t_{H\sharp F}\circ(\varphi^t_{G\sharp F})^{-1}=\varphi^t_{H}\circ(\varphi^t_{G})^{-1}$.
Note that $G\sharp F|_M\to G$ and $H\sharp F|_M\to H$ in $C^2$-topology provided that $F|_M\to 0$ in $C^2$-topology. The continuation property of the spectral invariant $\ell$ implies that $\ell(H,\alpha)=\ell(G,\alpha)$ by taking $F|_M\to 0$ in $C^2$-topology.
\end{proof}

\section{Ganor-Tanny barricades on Liouville domains}\label{sec:barricades}
In~\cite{GS} Ganor and Tanny constructed special perturbations of a Hamiltonian homotopy supported in a contact incompressible boundary domain (CIB) which, together with an almost complex structure, has \emph{barricades} in CIB such that Floer trajectories cannot enter or exit CIB.
Following closely their construction, we introduce barricades in the above sense on Liouville domains with admissible Lagrangian submanifolds.

In what follows, we say that $H\in C^\infty(\R\times[0,1]\times \widehat{M})$ is \emph{ stationary} for $|s|>R>0$ if $\partial_sH$ is supported in $[-R,R]\times[0,1]\times M$ and $H^s_t=\mu\rho+b$ outside $M$ where $H^s_t:=H(s,t,\cdot)$, and $\mu\geq0,b\in\R$ are two constants for all $s\in\R$ and $t\in[0,1]$. We denote the interior of a set $X$ by $X^\circ$.

\begin{df}%$W\cap L$ is an admissible Lagrangian for $W$
Let $(M,d\theta)$ be a Liouville domain and $L$ an admissible Lagrangian submanifold. Let $W_0=M_{r'}$ and $W_1=M_r$  with the properties that $L$ intersects $\partial W_i$ transversely, and for $\theta|_L=dk_L$, $k_L$ vanishes near $\partial W_i\cap L$ for $i=0,1$. Let $(H^s)_{s\in\R}\subset\cH_{<\tau}$ be a homotopy, stationary for sufficiently large $|s|$, from $H^-$ to $H^+$ and $J=(J_t)_{t\in[0,1]}$ a family of almost complex structures such that $J\in\mathcal{J}_\theta$. The pair $(H,J)$ is said to have a \emph{barricade} for $L$ on $\Omega=W_1\setminus W_0^\circ$ if the Hamiltonian chords of $H^\pm$ do not intersect the boundaries $\partial W_0$ and $\partial W_1$, and any solution $u:\R\times [0,1]\to \widehat{M}$ to (\ref{e:feq}) satisfying (\ref{e:boundary}) with asymptotic chords $x_\pm\in \mathcal{C}(L,H^\pm)$ at $\pm \infty$ satisfies
\begin{enumerate}
  \item[I.] if $x_-\subset W_0$ then $\im(u)\subset W_0$.
  \item[II.] if $x_+\subset W_1$ then $\im(u)\subset W_1$.
\end{enumerate}

\end{df}

Allowed and Forbidden Floer strips are illustrated in Figure~\ref{fig:bar}. We say that a pair $(H,J)$ of a Hamiltonian $H\in\cH_{<\tau}$ and a family of almost complex structures $(J_t)_{t\in[0,1]}\in\mathcal{J}_\theta$ has a barricade for $L$ on $\Omega$ if the pair of the constant homotopy $H^s\equiv H$ and $J=(J_t)_{t\in[0,1]}$ has a barricade in the above sense.

\begin{figure}[H]
	\centering
	\includegraphics[scale=0.5]{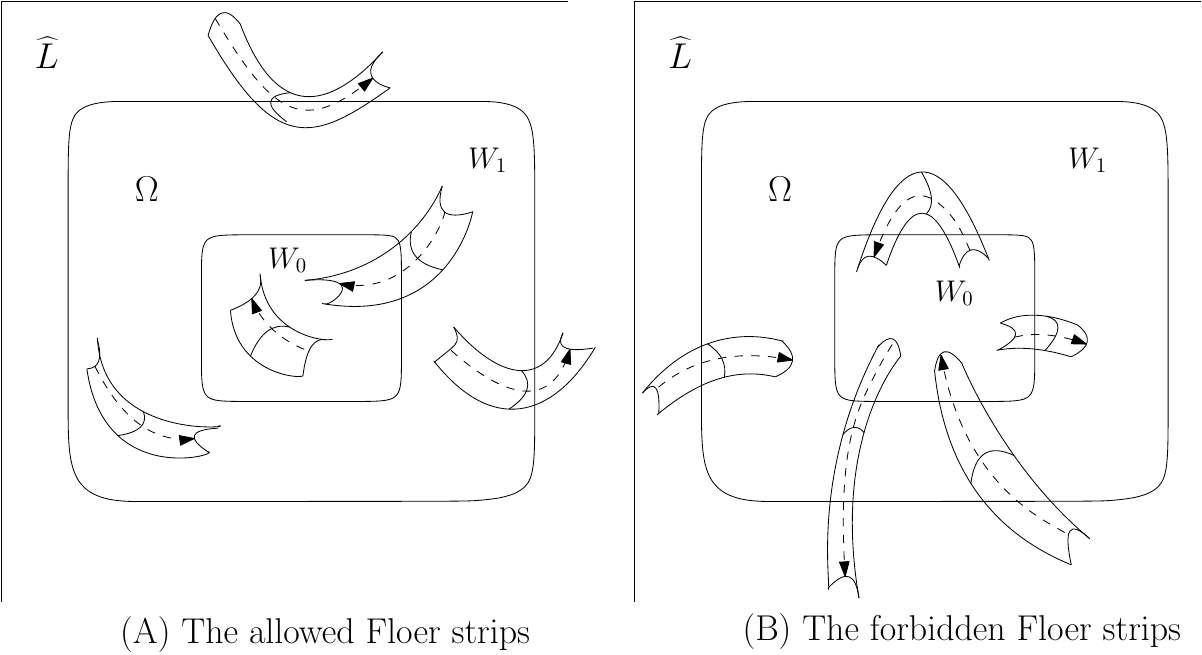}
	\caption{Floer strips in a barricade for $L$ on $\Omega$.}\label{fig:bar}
\end{figure}

\begin{thm}\label{barricade}
Let $(H^s)_{s\in\R}$ be a homotopy in $\cH_{<\tau}$, stationary for sufficiently large $|s|$, from $H^-$ to $H^+$ such that
 \begin{itemize}
   \item $H$ vanishes on $\R\times[0,1]\times M_{r}\setminus M_{r'}^\circ$ for two numbers $r,r'\in(0,1)$ with $r>r'$,
   %\item $H^s_t=h(\rho)$ on $\widehat{M}\setminus M_{r}$ for a smooth function $h:(0,\infty)\to\R$,
   %\item $h$ is $C^2$-small Morse function on $[r,1]\times\partial M$.
   \item $H^s_t=f^s(\rho)$ on $\widehat{M}\setminus M_r$ (with $\partial_s\partial_\rho f^s(\rho)= 0$ outside $M$),
   \item  for every $s\in\R$, $f^s$ is a $C^2$-small function and $\partial_sf^s\leq 0$ on $[r,1]\times\partial M$.
 \end{itemize}
 Then there exists a $C^\infty$-small perturbation $h$ of $H$ and a family of almost complex structures $(J_t)_{t\in[0,1]}$ such that the pairs $(h,J)$ and $(h^\pm,J)$ are Floer-regular and have a barricade on $\Omega$.
\end{thm}

The proof of Theorem~\ref{barricade} is an adaptation of the construction of barricades on CIB given by Ganor and Tanny in~\cite[Section~3]{GS}. The main tools in the proof involve $C^0$-estimates for Floer strips, transversality, Gromov-Floer compactness, and other standard techniques in Floer theory. The crucial insights from Ganor and Tanny are two aspects:
\begin{enumerate}
  \item One can control the size of the support of perturbations of a homotopy $H$ with Floer-regular ends $H^\pm$ with respect to $J\in\mathcal{J}_\theta$ to achieve transversality, as shown in Proposition~\ref{prop:reg}.
  \item Barricades are robust under $C^\infty$-small perturbations that have support in $I\times[0,1]\times M$, where $I\subset\mathbb{R}$ is a compact interval, as stated in Proposition~\ref{prop:survive}.
\end{enumerate}

We postpone the proof of Theorem~\ref{barricade} to Section~\ref{sec:bar} since a detailed analysis of the construction of barricades, which may be of independent interest, is necessary to complete the proof.

\begin{rmk}\label{rk:varmin}
From the construction of barricades, one can obtain the additional properties for the perturbation $h$ and the almost complex structure $J$ in the following:
\begin{itemize}
  \item The family of almost complex structures $J\in\cJ_\theta$ can be chosen to be time-independent.
If the restriction of the end $H^-\in\cH$ of $H$ to some Weinstein neighborhood of $L$ in $M$ is a lift of some time-independent $C^2$-small Morse function on $L$, and has a cylindrical bump  on $M_{r}\setminus M_{r'}^\circ$, then the perturbation $h$ can be chosen to satisfy $h^-=H^-$. %(see~Definition~\ref{df:bump})
  \item If the homotopy $H$ is constant on  some domain, then the perturbation $h$ can be chosen such that on this domain, $h^\pm$ have the same Hamiltonian $1$-chords.
  \item If $H^s$ is a constant homotopy for $s\notin[a,b]$, then the perturbation $h$ of $H$ can be chosen such that $supp(\partial_s h)\subset [a,b]\times[0,1]\times M$.
\end{itemize}

\end{rmk}

\section{The proof of the main theorem}

The basic idea of proving Theorem~\ref{Mthm} is to find a sequence of cofinal Hamiltonians which are radial and convex to compute the filtered wrapped Floer cohomology for an admissible $L$. This technique was used repeatedly in~\cite{FHW,Vi2,BPS,We,GX,BK}, etc.

\subsection{The proof of Theorem~\ref{Mthm}}

\begin{proof}

Note that if $HW^*(L)$ vanishes then the $c_{HW}(L)$ is finite.
In what follows we will show that for any $H\in\cH_c(M)$ it holds that $\ell(H,\mathds{1}_L)\geq -c_{HW}(L)$, and hence $\ga(L,H)\leq 2c_{HW}(L)$. So the diameter of metric space $(\mathcal{L}(L),\ga)$ has the upper bound $2c_{HW}(L)$.

We proceed in three steps:

\noindent \textbf{Step 1}. We construct a sequence of cofinal Hamiltonians to compute
the filtered wrapped Floer cohomology $HW_{(-a,\infty)}^*(L)$ with $a\in(0,\infty)\setminus\mathcal{R}(\partial L,\theta)$.

Given $\mu\in (0,\infty)$, we choose parameters $\epsilon,\delta\in(0,1)$, and denote by $\mathcal{F}_\mu$ the set of admissible functions $F\in\cH$ which satisfy
\begin{itemize}
  \item $F=-\epsilon$ on $M_{1-\delta}$;
  \item $F(\rho,x)=h(\rho)\;\hbox{for some smooth convex function $h$ on}\;(0,\infty)\times \partial M$;
  \item $F(\rho,x)=\mu(\rho-1)-\delta\;\hbox{on}\;\widehat{M}\setminus M\cong(1,\infty)\times \partial M$.
\end{itemize}

Here we ask that $\epsilon,\delta$ are sufficiently small. By~(\ref{action}) and our assumption that $k_L$ vanishes near $\partial L$, for every $F\in\mathcal{F}_\mu$, the absolute value of actions of non-constant Hamiltonian chords for $(L,F)$ approximate arbitrarily the periods of the corresponding Reeb chords.

For $a\in(0,\infty)\setminus\mathcal{R}(\partial L,\theta)$, we will show that for every $F\in \mathcal{F}_a$ the natural homomorphism
\[
\sigma_{F}:HW^*_{(-a,\infty)}(L,F)\longrightarrow HW^*_{(-a,\infty)}(L)
\]
is an isomorphism.

We first consider a sequence $\{\mu_k\}$ of increasing numbers  such that $\mu_1=a$ and
$\mu_k\to\infty$ as $k\to\infty$, and pick a sequence of functions $\{F_k\}$ with $F_1:=F$ such that $F_k\in\mathcal{F}_{\mu_k}$ and $F_k\leq F_{k+1}$ for all $k\in\N$. Moreover, we require that $\{F_k\}$ is an upward exhausting sequence of functions in the set
\[
\mathcal{H}^a=\{H\in\cH|-a\notin\spec(L,H)\}
\]
as illustrated in Figure~\ref{fig:cof}.
Indeed, since $\mathcal{R}(\partial L,\theta)$ is a nowhere dense set in $(0,\infty)$, there exist two positive numbers $\underline{\mu}<a<\overline{\mu}$ with $(\underline{\mu},\overline{\mu})\cap \mathcal{R}(\partial L,\theta)=\emptyset$. Associated to the function sequence $\{F_k\}$, we have two sequences of numbers
\[
A_k:=\frac{a-\epsilon_k}{1-\delta_k},\quad B_k:=\frac{a}{1+\delta_k}
\]
where $\epsilon_k\to0,\delta_k\to0$ as $k\to0$. Clearly, $\{A_k\},\{B_k\}$ are two sequences which converge to $a$. Without loss of generality, we may assume that $A_k,B_k\in(\underline{\mu},\overline{\mu})$ and $\theta|_{\widehat{L}\cap (\widehat{M}\setminus M_{1-\delta_k})}=0$ for all $k\in\N$. By our choice of $F_k$, the non-constant chords of each $F_k$ lie in the region between $\rho=1-\delta_k$ and $\rho=1$. Let $x$ be one of such chords. If the action of the chord $x$ is $-a$, then by~(\ref{action}) the line tangent to the graph of $F_k$ at $\rho=\rho(x_k)$ must pass through the point $(0,-a)$ and has the slope between $A_k$ and $B_k$, in contradiction to $(\underline{\mu},\overline{\mu})\cap \mathcal{R}(\partial L,\theta)=\emptyset$. Therefore, $F_k\in\mathcal{H}^a$ for all $k$.

\begin{figure}[H]
	\centering
\includegraphics[scale=0.4]{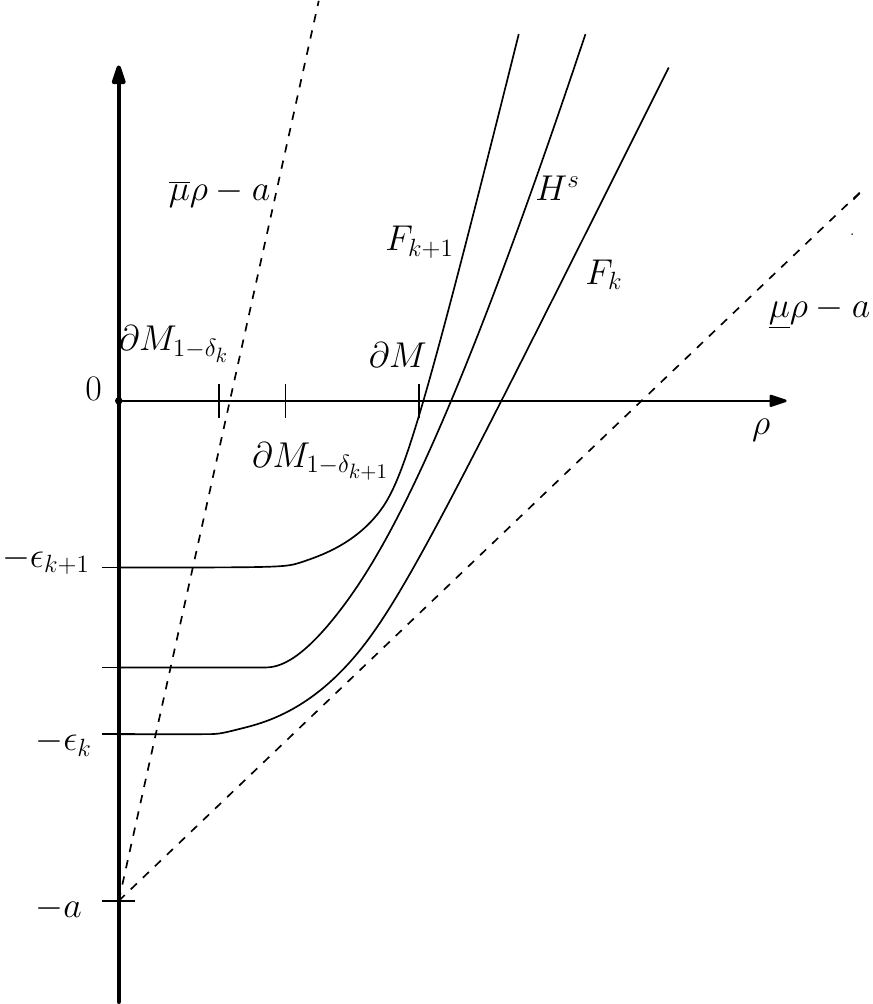}
	\caption{The monotone homotopy connecting $F_{k+1}$ to $F_k$.}\label{fig:cof}
\end{figure}

Next we consider a sequence of monotone homotopies $(H^s)_{s\in[0,1]}$ connecting $F_{k+1}$ to $F_{k}$ such that $-a\notin\spec(L,H^s)$ for every $s\in[0,1]$, see Figure~\ref{fig:cof}.

By Lemma~\ref{lem:mm}, the monotone homomorphisms
\[
\Phi_{F_kF_{k+1}}:HW^*_{(-a,\infty)}(L,F_k)\longrightarrow HW^*_{(-a,\infty)}(L,F_{k+1})
\]
are isomorphisms, and hence we have the isomorphism
\[
\sigma_{F_1}:HW^*_{(-a,\infty)}(L,F_1)\longrightarrow HW^*_{(-a,\infty)}(L)
\]
which is given by the natural homomorphism for the direct limit.

\noindent \textbf{Step 2}.
Given $c\in(c_{HW}(L),\infty)\setminus\mathcal{R}(\partial L,\theta)$, for any $F\in\mathcal{F}_{c}$ we will show that the PSS-map $\psi_{pss}^F$ is trivial, i.e. $\psi_{pss}^F(\mathds{1}_L)=0$.

Fix $\eta\in(0,\min\{\tau,c\})$. We pick a function $f\in\mathcal{F}_{\eta}$ such that $f\leq F$,  then  we have the following commutative diagram

\begin{equation}\notag%
\begin{split}
\xymatrix{H^*(L)
\ar[r]^{\psi^{f}_{pss}}
\ar[dr]_{\psi^{F}_{pss}}
&
HW^*(L,f)
\ar[r]_{\cong\quad}^{\sigma_f\quad}
\ar[d]^{\Phi_{fF}}
&
HW^*_{(-\eta,\infty)}(L)
\ar[d]^{\iota} \\
&
HW^*(L,F)
\ar[r]_{\cong\quad}^{\sigma_F\quad}
&
HW^*_{(-c,\infty)}(L)}
\end{split}
\end{equation}
where the map $\iota$ is induced by the action window map. By \textbf{Step 1} the natural homomorphisms $\sigma_f$ and $\sigma_F$ are isomorphisms.
Since $c>c_{HW}(L)$, by definition we have $\iota=0$.
Hence $\psi_{pss}^F(\mathds{1}_L)=0$.

\noindent \textbf{Step 3}. Let $H\in\cH_c(M)$.
For every positive number $c>c_{HW}(L)$ with $-c\notin \mathcal{R}^-(\partial L,\theta)\cup\spec(L,H)$, we will show that $\ell(H,\mathds{1}_L)\geq -c$.
This gives rise to $\ell(H,\mathds{1}_L)\geq-c_{HW}(L)$.

From now on we fix $F\in\mathcal{F}_{c}$ as in \textbf{Step 2}.
Assume that $H$ is compactly supported in $M_r$ for some $r\in(0,1)$.
We pick $H_i\in\cH, i=1,2$ such that
\begin{itemize}
  \item $H_1=H_2=H$ on $M_r$;
  \item $H_i=h_i(\rho)$ \;\hbox{for two smooth convex functions $h_i$ on}\;$(r,\infty)\times \partial M$ and $0< h_1'\leq h_2'$;
  \item $h_1,h_2$ are linear on $\widehat{M}\setminus M$ with $h_1'=\eta$ and $h_2'=c$ where $\eta$ is fixed as in \textbf{Step 2}.
\end{itemize}

By definition we have $\ell(H_1,\mathds{1}_L)=\ell(H,\mathds{1}_L)$. So it suffices to show that $\pi_{-c}\circ\psi_{pss}^{H_1}(\mathds{1}_L)=0$. To this end we consider the following
commutative diagram

\begin{equation}\notag%
\begin{split}
\xymatrix{&H^*(L)
\ar[r]^{\psi^{H_1}_{pss}\quad}
\ar[d]^{\psi^{F}_{pss}}
\ar[dr]^{\psi^{H_2}_{pss}}
&
HW^*(L,H_1)
\ar[r]^{\pi_{-c}\quad}
\ar[d]%^{\Phi_{H_1H_2}}
&
HW^*_{(-\infty,-c]}(L,H_1)
\ar[d]^{\Phi_{H_1H_2}}_{\cong} \\
&
HW^*(L,F)
\ar[r]^{\Phi_{FH_2}}_{\cong}
&
HW^*(L,H_2)
\ar[r]^{\pi_{-c}\quad}
&
HW^*_{(-\infty,-c]}(L,H_2)}
\end{split}
\end{equation}
Here the continuation map $\Phi_{FH_2}$ is an isomorphism because $F$ and $H_2$ have the same slope at infinity, and the monotone homomorphism $\Phi_{H_1H_2}$ is an isomorphism since there exists a homotopy $H_s$ between $H_1$ and $H_2$ such that $-c\notin\spec(L,H_s)$ for all $s\in[0,1]$.
Therefore, by \textbf{Step 2} the desired inequality follows from the above diagram.
\end{proof}

\subsection{The proof of Theorem~\ref{Mthm'}}
The key to proving Theorem~\ref{Mthm'} is the ``barricade" argument for Floer trajectories, originally developed by Ganor and Tanny~\cite{GS}. This argument leads to the following result, which is  an adaptation of \cite[Lemma B]{Ma} to the wrapped setting.

\begin{prop}\label{prop:nonpos}
Let $(M^{2n},d\theta)$ be a Weinstein domain and $L^n\subset M$ an admissible connected Lagrangian submanifold. For any Hamiltonian $H\in C^\infty([0,1]\times M)$ with compact support in $[0,1]\times int(M)$ it holds that $\ell(H,\mathds{1}_L)\leq 0$.
\end{prop}% satisfying $\theta|_L=0$

We now complete the proof of Theorem~\ref{Mthm'} assuming Proposition~\ref{prop:nonpos}. The strategy of the proof is to construct a family of compactly supported Hamiltonians $H$ on $M$ that have sufficiently large oscillations to force $-\ell(H,\mathds{1}_L)$ to be large enough. %similar to the computation of the homological BPS capacity in~\cite[Section~4]{We} or~\cite[Section~9.4]{GX}.

\begin{proof}[The proof of Theorem~\ref{Mthm'}]

By our assumption that $\theta|_L=dk_L$ and $k_L=0$ on a neighborhood of $\partial L$, we can modify the Liouville one-form such that $\theta|_L=0$ as follows. Indeed, by the connectedness of $L$
we can first extend $k_L$ to a smooth function on $L$ and next to a compactly supported smooth function $f:M\to\R$, then add $-df$ to $\theta$. With respect to this new one-form $\theta$, $L$ is still an exact Lagrangian but with $\theta|_L=0$. Note that this progress does not change the symplectic form, and hence the actions of chords of any Hamiltonian with ends in $L$ remain the same. Besides, it only changes
$\theta$ in the interior of $M$, but not on the boundary $\partial M$, so the Reeb vector field does not change. Hence, with respect to this new form, the corresponding spectral invariant is not changed, and this is what we are ultimately interested in. In the following proof we assume that $\theta|_{\widehat{L}}=dk_L$ with $k_L\equiv 0$ on $\widehat{L}$.

Fix $\rho_0\in(0,1)$ and let $\delta\in (0,\rho_0)$. Pick a sequence of numbers $(a_k)$ such that $a_k\in (0,\infty)\setminus\mathcal{R}(\partial L,\theta)$ and $a_k\to \infty$ as $k\to\infty$. For each $k\in\N$, we take a piecewise linear curve $c_k$. More precisely, we begin with a horizontal line segment with the starting point $(0,-a_k(\rho_0-\delta))$. Upon reaching the point $(\delta,-a_k(\rho_0-\delta))$, we follow the line with slope $a_k$ until meeting the $\rho$-axis, then follow it to the right until we arrive at the point $(0,1)$, and finally we go along the line with slope $\eta\in(0,\tau)$. We smooth out the corners of each curve $c_k$, and hence obtain a sequence of functions $H_k\in\cH_{<\tau}$ whose graphs are these smooth curves  as illustrated in Figure~\ref{fig:fun}.

\begin{figure}[H]
	\centering
	\includegraphics[scale=0.5]{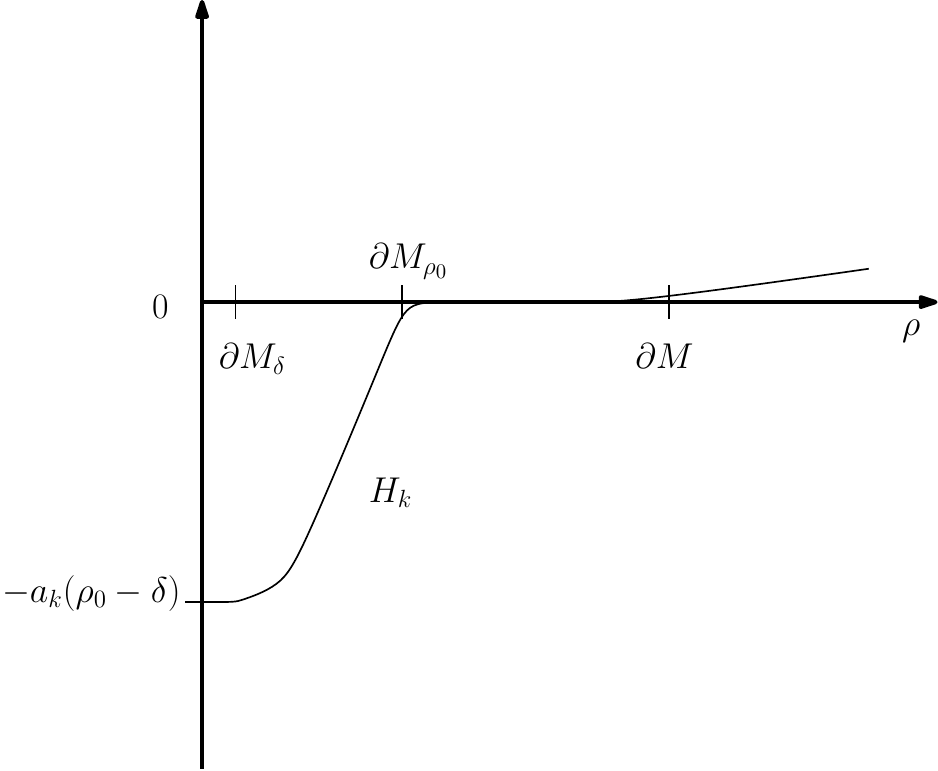}
	\caption{The function $H_k$.}\label{fig:fun}
\end{figure}

The chords of $H_k$ are divided into four classes according to the regions which they lie in. The only non-constant $1$-period chords of $H_k$ arise near $\rho=\delta,\rho_0$. The actions of non-constant chords for $h_k$ are $-\rho h_k'(\rho)+h_k(\rho)$. Let $\eta_k$ denote the distance between $a_k$ and $\mathcal{R}(\partial L,\theta)$. Then $\eta_k>0$ since the Reeb periods are a discrete points in $[0,\infty)$. The non-constant chords of $H_k$ near $\rho=\delta, \rho_0$ have actions in \[[-a_k\rho_0+\eta_k\delta, -\tau\delta-a_k(\rho_0-\delta)]\;\hbox{and}\;[-\rho_0(a_k-\eta_k),-\rho_0\tau]\] respectively. The constant chords of $H_k$ have actions
close to $-a_k(\rho_0-\delta)$ in the region $M_\delta$ and to $0$ outside $M_{\rho_0}$.

For sufficiently small $\delta>0$ there is a positive number $\epsilon_k\in(0,1)$ such that
$a_k\delta<\epsilon_k<\eta_k\rho_0$. Fixing such $\delta$ we can separate the action values of the $1$-period chords lying in $M_\delta$ and those outside of\; $M_\delta$. More precisely, the action values of chords lying in two regions near $\partial M_\delta$ and $\partial M_{\rho_0}$ belong to $(-\infty,-a_k\rho_0+\epsilon_k)$ and $(-a_k\rho_0+\epsilon_k,\infty)$ respectively.

Now we deform $H_k$ by the monotone homotopy $H_k^s=h^s_k(\rho)$ as indicated in Figure~\ref{fig:hmty} to a function $F_k$ which is convex outside of\; $M_{\delta}$. The graph of this new function is obtained by following $H_k$ until we arrive at $\rho=\rho_0$ and keep going linearly with slope $a_k$. Since the values of intersections of vertical axis with the tangent lines to the graphs of $h_k^s$ are larger than $-a_k\rho_0+\epsilon_k$ where non-constant chords occur, by Lemma~\ref{lem:mm} we  obtain the monotone isomorphism
\[
\Phi_{H_kF_k}:HW^*_{(-\infty,-a_k\rho_0+\epsilon_k]}(L,H_k)\stackrel{\cong}{\longrightarrow}HW^*_{(-\infty,-a_k\rho_0+\epsilon_k]}(L,F_k).
\]

Next we deform $F_k$ into a function $F_k'$ with the property $F_k'|_{M}<0$ by a monotone homotopy $(F^s)_{s\in[0,1]}$ as  indicated in Figure~\ref{fig:Fk}. Note that during the homotopy, $-a_k\rho_0+\epsilon_k$ is not the action of a Hamiltonian chord with period one for every pair $(L,F^s)$.
By Lemma~\ref{lem:mm} again, we obtain the isomorphism
\[
\Phi_{F_k'F_k}:HW^*_{(-\infty,-a_k\rho_0+\epsilon_k]}(L,F_k')\stackrel{\cong}{\longrightarrow}
HW^*_{(-\infty,-a_k\rho_0+\epsilon_k]}(L,F_k).
\]
Note that the inverse of the above map is indeed induced by the continuation map 
\[
\Phi_{F_kF_k'}:HW^*(L,F_k)\longrightarrow HW^*(L,F_k')
\]
which is well-defined since $\mu_{F_k}=\mu_{F_k'}$. 

\begin{figure}[htbp]
\centering
\subfloat[{Figure 7a}]{
\label{fig:hmty}
\includegraphics[scale=0.4]{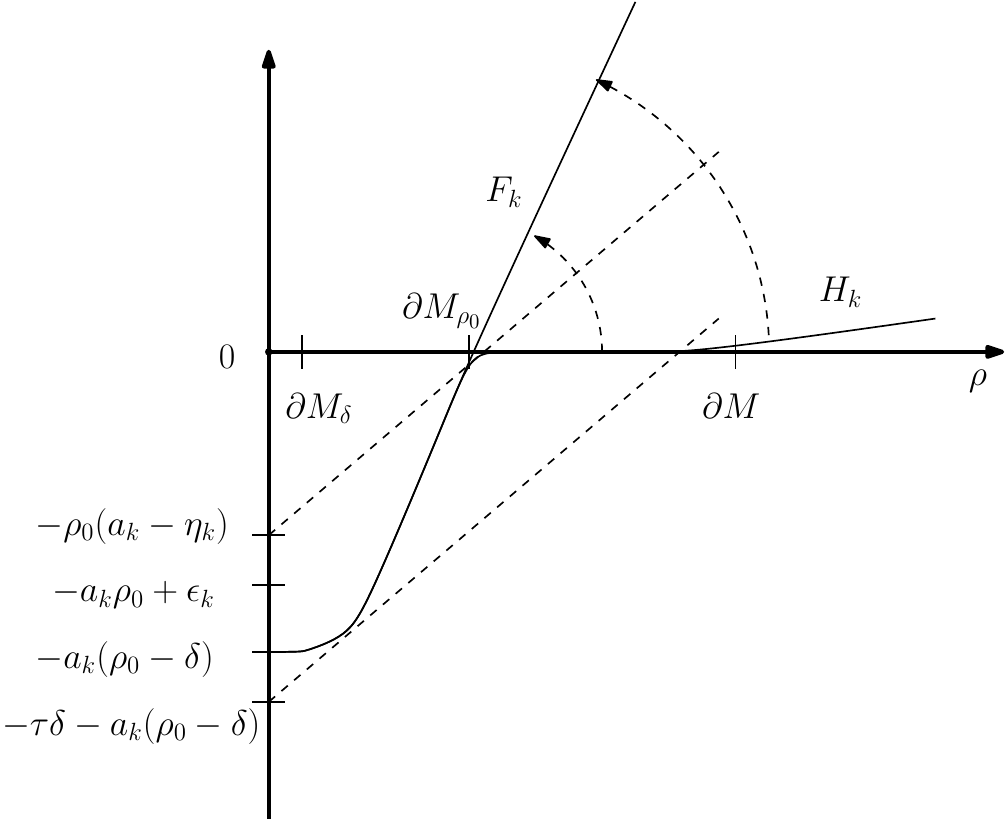}}
\subfloat[{Figure 7b}]{
\label{fig:Fk}
\includegraphics[scale=0.4]{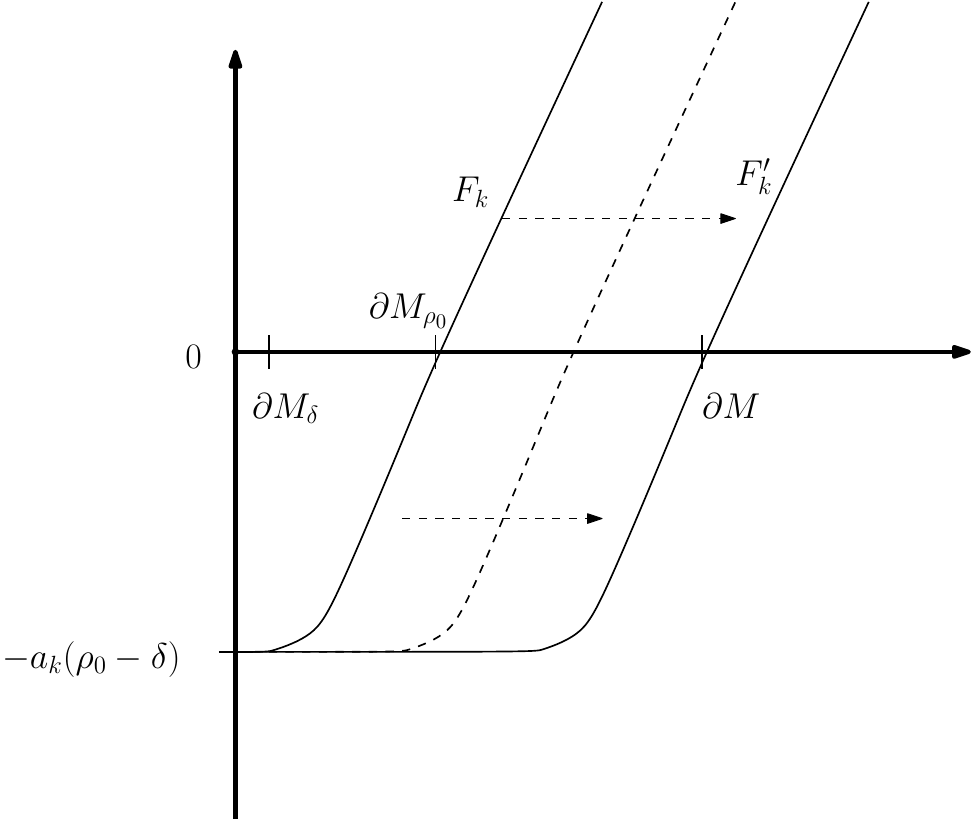}}
\caption{The homotopy connecting $F'_k$ to $H_k$.}
\end{figure}

So we get an isomorphism
\[
\Psi_{H_kF_k'}\triangleq\Phi_{F_k'F_k}^{-1}\circ \Phi_{H_kF_k}:HW^*_{(-\infty,-a_k\rho_0+\epsilon_k]}(L,H_k)\stackrel{\cong}{\longrightarrow}
HW^*_{(-\infty,-a_k\rho_0+\epsilon_k]}(L,F_k')
\]
which is induced by the continuation map $\Phi_{H_kF_k'}:HW^*(L,H_k)\to HW^*(L,F_k')$ (the linear homotopy from $H_k$ to $F_k'$ is monotone at infinity).

Fix $c>0$. For simplicity we set $A_k=-a_k\rho_0+\epsilon_k$ (tending to $-\infty$ as $k\to \infty$).
Since by our construction all chords of $F_k'$ have actions less than the negative number $A_k$,  we have
$HW^*_{(-\infty,A_k]}(L,F_k')= HW^*_{(-\infty,c]}(L,F_k')=HW^*_{(-\infty,\infty)}(L,F_k')$.

Then we have the following diagram:
\begin{equation}\label{dia:comm}
\begin{split}
\xymatrix{
&
HW^*(L,H_k)
\ar[r]^{\pi_{A_k}\qquad}
\ar[d]_{\Phi_{H_kF_k'}}
&
HW^*_{(-\infty,A_k]}(L,H_k)
\ar[d]^{\Psi_{H_kF_k'}}_{\cong} \\
&
HW^*(L,F_k')
\ar[r]^{\pi_{A_k}\qquad}
\ar[d]_{\sigma_{F_k'}}
\ar[dr]^{\pi_c}
&
HW^*_{(-\infty,A_k]}(L,F_k')
\\
&
HW^*(L)
\ar[dr]^{\pi_c^L}_{\cong}
&
HW^*_{(-\infty,c]}(L,F_k')
\ar[d]^{\sigma_{F_k'}}
\ar[u]_{\pi^{-\infty,c}_{-\infty,A_k}}^{\cong}
\\
& &
HW^*_{(-\infty,c]}(L)
}
\end{split}
\end{equation}

The top rectangular block in (\ref{dia:comm}) commutes, because (\ref{diag:IQ}) implies  that the following diagram commutes
\begin{equation}\notag%
\xymatrix{
&
HW^*(L,H_k)
\ar[r]^{\pi_{A_k}\qquad}
\ar[d]_{\Phi_{H_kF_k}}^{\cong}
&
HW^*_{(-\infty,A_k]}(L,H_k)
\ar[d]^{\Phi_{H_kF_k}}_{\cong} \\
&
HW^*(L,F_k)
\ar[r]^{\pi_{A_k}\qquad}
&
HW^*_{(-\infty,A_k]}(L,F_k)
\\
&
HW^*(L,F_k')
\ar[r]^{\pi_{A_k}}
\ar[u]^{\Phi_{F_k'F_k}}_{\cong}
&
HW^*_{(-\infty,A_k]}(L,F_k')
\ar[u]_{\Phi_{F_k'F_k}}^{\cong}
}
\end{equation} %\[\Phi_{H_kF_k}\circ\pi_{A_k}=\pi_{A_k}\circ\Phi_{H_kF_k'},\quad \pi_{A_k}\circ\Phi_{F_k'F_k}=\Phi_{F_k'F_k}\circ\pi_{A_k}.\]

Obviously, the middle triangle in (\ref{dia:comm}) commutes. The bottom parallelogram block in~(\ref{dia:comm}) commutes due to (\ref{diag:IQw}). By (\ref{diam:hk}), we get
\[
\sigma_{H_k}=\sigma_{F_k'}\circ\Phi_{H_kF_k'}:HW^*(L,H_k)\longrightarrow HW^*(L).
\]

Since $\sigma_{H_k}\circ\psi_{pss}^{H_k}$ respects the ring structures on $H^*(L)$ and $HW^*(L)$, we deduce from our assumption $HW^*(L)\neq 0$ that
$\sigma_{H_k}(\psi_{pss}^{H_k}(\mathds{1}_L))\neq0$. It follows from the above diagram that
\[
\pi_{A_k}\big(\psi_{pss}^{H_k}(\mathds{1}_L)\big)\neq0
\]
which implies $\ell(H_k,\mathds{1}_L)\leq A_k$.
Consequently, by our construction of the sequence $\{H_k\}$ and the continuation property of the spectral invariant $\ell$, one can find a sequence of Hamiltonians $\{H'_k\}\subset \cH_c(M)$ such that $\ell(H'_k,\mathds{1}_L)\to-\infty$ as $k\to\infty$. Put $L_k=\varphi_{H'_k}^{-1}(L)$ for all $k\in\N$. By Proposition~\ref{prop:nonpos} and Corollary~\ref{cor:indep}, we have
\begin{eqnarray}\notag
\gamma(L_k,L)&=&\inf\big\{-\ell(H,\mathds{1}_L)-\ell(\overline{H},\mathds{1}_L)\big|\varphi_H^1(L_k)=L,\;H\in\cH_c(M)\big\}\notag\\
&\geq&\inf\big\{-\ell(H,\mathds{1}_L)\big|\varphi_H^{-1}(L)=L_k,\;H\in\cH_c(M)\big\}\notag\\
&=&-\ell(H_k',\mathds{1}_L)\notag
\end{eqnarray}
which concludes the desired result by letting $k\to\infty$.
\end{proof}

\begin{figure}[H]
	\centering
	\includegraphics[scale=0.5]{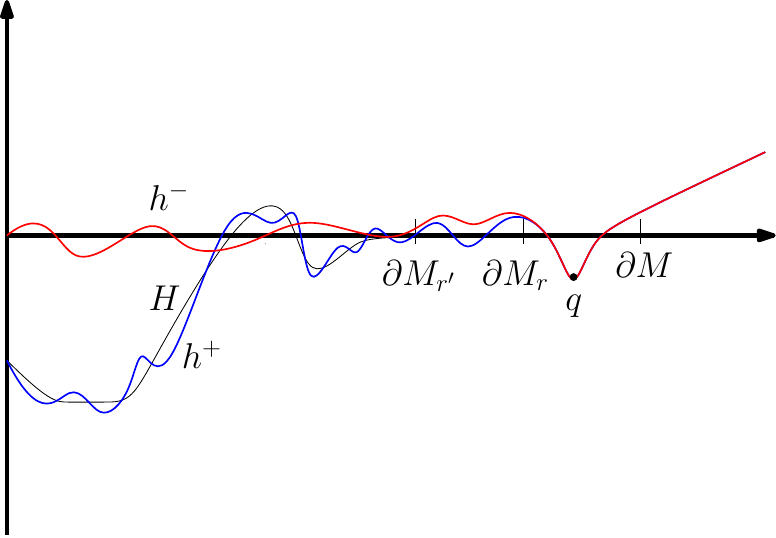}
	\caption{The various Hamiltonians in the proof of Proposition~\ref{prop:nonpos}.}\label{fig:varfmin}
\end{figure}

\begin{proof}[The proof of Proposition~\ref{prop:nonpos}]
Assume that $supp(H)\subset [0,1]\times M_{r'}$ for some $r'\in(0,1)$. For sufficiently small $\epsilon>0$, we pick a time independent Hamiltonian $F\in\cH_{<\tau}$ such that the restriction of $F$ to some Weinstein neighborhood of $L$ in $M$ is a lift of a $C^2$-small Morse function $f\in C^\infty(L)$ with $\|F\|_{C^2(M)}<\epsilon$. We also pick $K\in\cH_{<\tau}$ such that $K=H$ on $M_{r}$ for some $r\in (r',1)$, and $K|_{\widehat{M}\setminus M_r}=F|_{\widehat{M}\setminus M_r}$.
We ask further that $F$ has a \textbf{unique local (and global) minimum point} $q$ contained in $L\cap (M\setminus M_r)$. Let $H^s$ be a linear homotopy from $H^-=F$ to $H^+=K$ as in (\ref{eq:hmtp}). By Theorem~\ref{barricade} and Remark~\ref{rk:varmin}, there exists a $C^\infty$-small perturbation $h^s$ of $H^s$ with $h^-=H^-$ and $h^+|_{\widehat{M}\setminus M_r}=H^+|_{\widehat{M}\setminus M_r}$ and a family of almost complex structures $(J_t)_{t\in[0,1]}$ such that the pairs $(h,J)$ and $(h^\pm,J)$ are Floer-regular and have a barricade on $M_r\setminus M_{r'}^\circ$ for $r$ close enough to $r'$. We refer to Figure~\ref{fig:varfmin} for a schematic graph of these Hamiltonians.

Note that we have a chain-level isomorphism between the Morse complex $CM^*(L,f,g)$ and the wrapped Floer complex $CW^*(L,h^-)$  via the PSS-map $\psi_f^{h^-}$ by sending critical points of $f$ to the corresponding constant chords of $h^-$. So the point $q$ represents the fundamental class
$$[q]=\mathds{1}_L\in H^*(L)\xrightarrow[\Psi_{pss}^{h^-}]{\cong} HW^*(L,h^-).$$
By our construction, $h^+$ and $h^-$ have the same slope outside of\;$M$. It follows from Lemma~\ref{lem:slope} that
the continuation map $\Phi_{h^+h^-}:CW^*(L,h^+)\to CW^*(L,h^-)$ associated to $h$ induces an isomorphism in homology. Therefore, there exists a cycle representing $\Phi_{h^+h^-}^{-1}([q])$ in $CW^*(L,h^+)$ which is mapped by $\Phi_{h^+h^-}$ to
the constant chord $q$. Since the pair $(h,J)$ has a barricade for $L$ on $M_r\setminus M_{r'}^\circ$, the Floer strip $u$ of this continuation map starting at $x_-=q\in L\cap(M\setminus M_r)$ must end at some chords $x_+$ of $h^+$ outside of\; $M_r$, and hence the cycle $\Phi_{h^+h^-}^{-1}(q)$ equals to the sum $\sum_ip_i$ of some critical points $p_i$ of $h^+$ on $L\cap (M\setminus M_r)$. Since $h^+$ is $C^\infty$-close to $H^+=K$ and the values of $K$ are close to zero in the region $M\setminus M_r$, the action of every constant chord $p_i$ is close enough to zero. It follows from the continuity property of $\ell$ and $\Phi_{h^+h^-}\circ \Psi_{pss}^{h^+}(\mathds{1}_L)=\Psi_{pss}^{h^-}(\mathds{1}_L)$ (see~Lemma~\ref{lem:iso}) that
\[
\ell(H,\mathds{1}_L)\leq \ell(K,\mathds{1}_L)+\epsilon\leq \ell(h^+,\mathds{1}_L)+2\epsilon\leq \max_{p_i}\mathcal{A}_{L,h^+}(p_i)+2\epsilon\leq 3\epsilon.
\]
Since $\epsilon>0$ can be arbitrarily small, we obtain the desired inequality.
\end{proof}

\section{The proof of Theorem~\ref{barricade}}\label{sec:bar}

\subsection{Constructing Ganor-Tanny barricades}

%Let $Y$ denote the restriction of Liouville vector field $V_\theta$ on $\Omega$.
We consider a pair $(H,J)$ of a Hamiltonian homotopy $(H^s)_{s\in\R}$ and a family of almost complex structures $(J_t)_{t\in[0,1]}$. Fix $\mu\in(0,\tau)$ (where $\tau=\min\mathcal{R}(\partial L,\theta)$).

\begin{df}\label{df:bump}
We say that the pair $(H,J)$ has a \emph{cylindrical bump} of slope $\mu$ on $\Omega=W_1\setminus W_0^\circ$ if
\begin{itemize}
  \item[1.] $H=0$ on $\R\times [0,1]\times\partial \Omega$ and $\partial_s H^s\leq 0$ outside of\;$W_0$;
  \item[2.] $J$ is of contact type, i.e. $d\rho\circ J=-\theta$ (or equivalently $JY=R$), near the boundaries $\partial W_0,\partial W_1$;
  \item[3.] $\nabla_J H=\mu V_\theta$ near $\R\times[0,1]\times\partial W_0$ and $\nabla_J H=-\mu V_\theta$ near $\R\times[0,1]\times\partial W_1$  where $\nabla_J$ is the gradient with respect to the metric $d\lambda(\cdot,J\cdot)$;
  \item[4.] If $x_\pm\in \mathcal{C}(L,H^\pm)$ are not contained in $W_0$, then $x_\pm$ are critical points of $H^\pm$ on $\widehat{L}$ with values in the interval $(-\mu,\mu)$.
\end{itemize}
\end{df}

\begin{prop}\label{prop:bump}
Assume that $(H,J)$ is a pair with a cylindrical bump of slope $\mu\in(0,\tau)$ on $\Omega$. Then $(H,J)$ has a barricade on $\Omega$.
\end{prop}

For proving Proposition~\ref{prop:bump}, we use three lemmata as in~~\cite{GS} to exclude certain types of Floer trajectories. In what follows we denote $W:=W_0$ or $W_1$ for simplicity.

\begin{df}
Let $(H,J)$ be a pair of a Hamiltonian homotopy and a family of almost complex structures. We say that
$(H,J)$ is \emph{$\mu$-cylindrical} near $\partial W$ with $\mu\in\R\setminus\{0\}$ if
\begin{itemize}
  \item $J$ is of contact type, i.e. $d\rho\circ J=-\theta$ near $\partial W$;
  \item $H\in\cH_{<\tau}$ is independent of the $\R$-coordinate and the time coordinate near $\partial W$ and $\R\times[0,1]\times\partial W=\{H=a\}$ is a regular level set of $H$;
  \item $\nabla_J H=\mu V_\theta$ on a neighborhood of $\partial W$ and $H$ has no chords $x\in\hc$ intersecting this neighborhood.
\end{itemize}
\end{df}

The first lemma follows from an argument that has appeared first in~\cite[Lemma~7.2]{AS}. We include the proof for completeness because our setup differs slightly from the one there.
\begin{lem}[No escape lemma]\label{lem:maxi}
Let $(H,J)$ be a pair which is $\mu$-cylindrical near $W$. If $\partial_s H\leq 0$ on $W^c$, then every finite-energy solution $u$ with both asymptotes contained in $W$ is entirely contained in $W$.
\end{lem}
\begin{proof}
Suppose the contrary that $u$ is not entirely contained in $W$.
Note that $\Sigma:=u^{-1}(\widehat{M}\setminus W^\circ)$ is a compact surface with corners and the corners divide the boundary $\partial \Sigma$ into two pieces: the piece landing in the boundary $\partial W$ and the one landing in $\widehat{L}$. We write $\partial \Sigma=\partial_b \Sigma\cup \partial_l \Sigma$ according to these two pieces. Clearly, by our assumption $\partial_b \Sigma\neq\emptyset$.
We denote by $j$ the restriction of the complex structure from the strip $\R\times[0,1]$ to $\Sigma$.

The Floer equation for $u$ can be read as $(du-X_{H^s}(u)\otimes dt)^{0,1}=0$. Using $\partial_sH^s\leq 0$ outside of\;$W$ and Stokes' theorem, we have
\begin{eqnarray}\label{e:EE}
E(u|_\Sigma)&=&\frac{1}{2}\int_\Sigma \big|du-X_{H^s}\otimes dt\big|^2\hbox{Vol}_\Sigma\notag\\
&=&\int_\Sigma u^*d\theta-u^*(dH^s)\wedge dt\notag\\
&=&\int_\Sigma d\big(u^*\theta-(u^*H^s)dt\big)+(\partial_sH^s)ds\wedge dt\notag\\
&\leq&\int_\Sigma d\big(u^*\theta-(u^*H^s)dt\big)\notag\\
&=&\int_{\partial\Sigma} u^*\theta-(u^*H^s)dt\notag\\
\end{eqnarray}
For any connected component $\gamma$ of $\partial_l\Sigma$ we have $(u^*H^s)dt|_\gamma=0$, and since $u^*\theta|_L=u^*dk_L$ by Stokes' theorem we get $\int_\gamma u^*\theta=0$ for circles $\gamma$, while for intervals $\gamma$, $\int_\gamma u^*\theta=k_L(p)-k_L(q)$ for corners $p,q\in\partial_b\Sigma\cap \partial_l\Sigma$ and hence this also vanishes by the assumption that $k_L|_{\partial W\cap L}=0$. Therefore, from the last term in (\ref{e:EE}) we deduce that $E(u|_\Sigma)\leq\int_{\partial_b\Sigma} u^*\theta-(u^*H^s)dt$.
Since $\nabla_J H=\mu V_\theta$ and $d\rho\circ J=-\theta$ on a neighborhood of $\partial W$, we have $X_{H^s}=J\nabla_J H=\mu R$ and hence $\theta(\frac{a}{\eta}X_{H^s})=H$ along $\partial W$ where $\eta=\mu r$ or $\mu r'$ depending on $W=M_r$ or $M_{r'}$. Using this and the contact condition $d\rho\circ J=-\theta$ near $\partial W$,
\begin{eqnarray}\label{e:EE'}
E(u|_\Sigma)&\leq&\int_{\partial_b\Sigma} \theta\circ\bigg(du-\frac{a}{\eta}X_{H^s}(u)\otimes dt\bigg)\notag\\
&=&\int_{\partial_b\Sigma} (\theta\circ J)\circ\bigg(du-\frac{a}{\eta}X_{H^s}(u)\otimes dt\bigg)\circ(-j)\notag\\
&=&\int_{\partial_b\Sigma} d\rho\circ\bigg(du-\frac{a}{\eta}X_{H^s}(u)\otimes dt\bigg)\circ(-j)\notag\\
&=&\int_{\partial_b\Sigma} d\rho\circ du\circ(-j).
\end{eqnarray}
Let $\xi$ be a tangent vector to $\partial_b\Sigma$ which gives rise to the boundary orientation. Then $j\xi$ points into $\Sigma$, and thus $du(j\xi)$ does not point outwards along $\partial W$, so $d\rho\circ du(j\xi)\geq 0$. Integrating that it follows from (\ref{e:EE'}) that $E(u|_\Sigma)=0$. This implies that each connected component of $u|_\Sigma$ is contained in a single orbit of $X_{H^s}$. Since $\partial_b \Sigma\neq\emptyset$ and $X_{H^s}$ is tangent to $\partial W$ this orbit must be contained in $\partial W$. This contradicts our assumption that $H$ has no chords near $\partial W$.
\end{proof}

\begin{figure}[H]
	\centering
\qquad\includegraphics[scale=0.5]{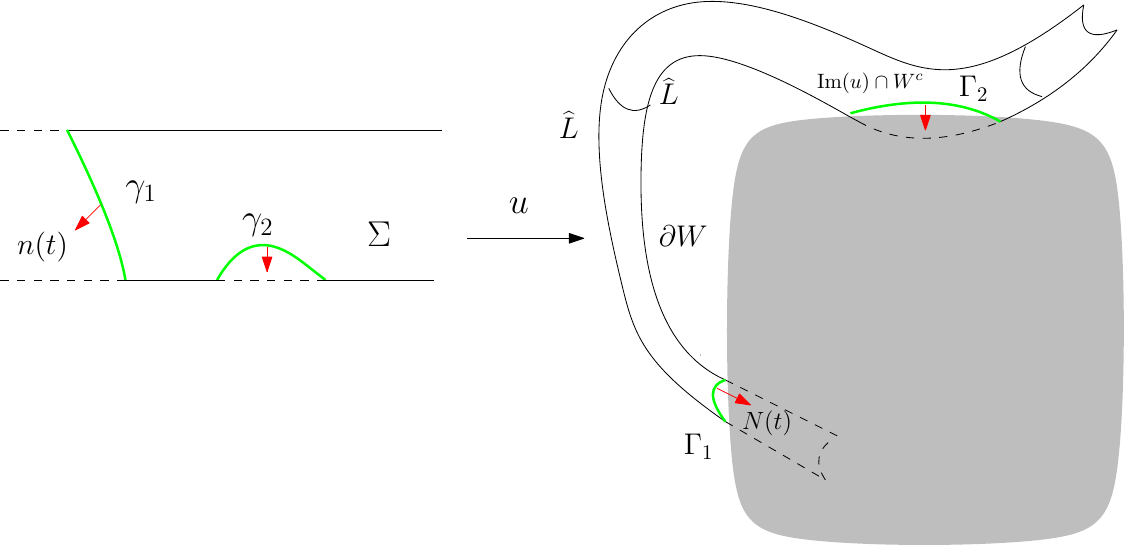}
	\caption{The image of a solution $u$ crossing the boundary $\partial W$.}\label{fig:inters}
\end{figure}

Next we give an upper bound for the integral of $\theta$ along the oriented curve $\gamma:=\partial (\im(u)\cap W^c)$ as illustrated in Figure~\ref{fig:inters},  where $u:\R\times [0,1]\to \widehat{M}$ is a solution of the $s$-dependent Floer equation~(\ref{e:feq}) with finite energy and the boundary condition (\ref{e:boundary}).

\begin{lem}\label{lem:int}
Let $u$ be as above with asymptotic chords $x_\pm\in \mathcal{C}(L,H^\pm)$. If $u$ intersects $\partial W$ transversely, then
\begin{equation}\notag
\int_\Gamma\theta\leq
\begin{cases}
-\mu,& \hbox{if}\;x_-\subset W,\;x_+\subset W^c,\\
\mu, &  \hbox{if}\;x_-\subset W^c,\;x_+\subset W,\\
0, & \hbox{if}\;x_\pm\subset W\;\hbox{or}\;x_\pm\subset W^c
\end{cases}
\end{equation}
where  $\Gamma:=\im(u)\cap \partial W$ is oriented as the boundary of $\im(u)\cap W^c$.

\end{lem}

The third lemma is an application of Lemma~\ref{lem:int} which is useful to bound the actions of the ends of Floer strips that cross the boundary of $W$ provided that the homotopy $H$ is non-increasing outside of\;$W$.
\begin{lem}\label{lem:action}
Let $(H,J)$ be a pair which is $\mu$-cylindrical near $W$. If $\partial_sH\leq 0$ on $W^c$, then every finite-energy solution $u$ with asymptotic chords $x_\pm\in \mathcal{C}(L,H^\pm)$ satisfies
 $\mathcal{A}_{L,H^+}(x_+)<a-\mu$ whenever $x_-\subset W$ and $x_+\subset W^c$, or
 $\mathcal{A}_{L,H^-}(x_-)>a-\mu$ whenever $x_-\subset W^c$ and $x_+\subset W$. Here $a$ is the value of $H$ on $\partial W$.

\end{lem}

The proofs of Lemma~3.2 and Lemma~3.4 in~\cite{GS} can be carried over to the above two lemmata respectively in a direct fashion. The only difference in the proofs is that in our case the portion $\Sigma:=\im(u)\cap W^c$ of a Floer strip $u$ outside of\;$W$ has an additional boundary $\partial_l\Sigma$ landing in $L$. But this would not cause new problems since in the proof of Lemma~\ref{lem:int} the extra term $\int_{\partial_l\Sigma}dt$ in the formula~(16) on page~664 in~\cite{GS} vanishes, and since in the proof of Lemma~\ref{lem:action} $\theta|_L=dk_L$ with $k_L$ vanishing near $\partial W\cap L$ the extra term $\int_{\partial_l\Sigma}\theta$ in the formula~(19) on page~666 in~\cite{GS} disappears.

As a consequence of Lemma~\ref{lem:action}, we have the following:
\begin{prop}\label{prop:bari}
If $(H,J)$ has a cylindrical bump of slope $\mu\in(0,\tau)$ on $\Omega$, then every solution $u:\R\times [0,1]\to \widehat{M}$ to (\ref{e:feq}) satisfying (\ref{e:boundary}) with asymptotic chords $x_\pm\in \mathcal{C}(L,H^\pm)$ has the following properties:
\begin{enumerate}
  \item if $x_-\subset W_0$ and $x_+\subset W_0^c:=\widehat{M}\setminus W_0$ then $\mathcal{A}_{L,H^+}(x_+)<-\mu$;
   \item if $x_+\subset W_1$ and $x_-\subset W_1^c:=\widehat{M}\setminus W_1$ then $\mathcal{A}_{L,H^-}(x_-)>\mu$;
   %\item if $x_-\subset M$ and $x_+\subset \widehat{M}\setminus M$ then $\mathcal{A}_{L,H^+}(x_+)>\mu$;
   %\item if $x_+\subset W$ and $x_-\subset \widehat{M}\setminus W$ then $\mathcal{A}_{L,H^-}(x_-)<-\mu$.
\end{enumerate}

\end{prop}

\begin{proof}[The proof of Proposition~\ref{prop:bump}]
Let $u:\R\times [0,1]\to \widehat{M}$ be a Floer strip for the pair $(H,J)$ with asymptotic chords $x_\pm\in \mathcal{C}(L,H^\pm)$ at $\pm \infty$. We only prove the first case in the definition of barricade since the proof of the second case is similar. Suppose that $x_-\subset W_0$. If $x_+\subset W_0^c$, then $x_+$ is a critical point of $H^+$ on $\widehat{L}$ with value in $(-\mu,\mu)$. On the other hand, it follows from the first statement of Proposition~\ref{prop:bari} that $\mathcal{A}_{L,H^+}(x_+)<-\mu$. So we get a contradiction and hence $x_+\subset W_0$. Then Lemma~\ref{lem:maxi} concludes that $\im(u)\subset W_0$.
\end{proof}

\subsection{Transversality}

To achieve transversality of moduli spaces we would like to perturb a given homotopy in suitable Banach spaces. 
The Floer $C_\varepsilon$-space introduced by Floer~\cite{Fl0} is a separable Banach space, so one can apply Sard-Smale theorem for genericity arguments.
For our purpose, given a compact interval $I\subset\R$ with non-empty interior, we consider the following perturbation space $C^\infty_{\varepsilon,I}(M)$ defined as follows.

For $h\in C^\infty(\R\times[0,1]\times \widehat{M})$ with $supp(h(s,t,\cdot))\subset M$ for all $s,t$, if $\varepsilon=(\varepsilon_k)_{k=0}^\infty$ is a sequence of positive numbers with $\varepsilon_k\to0$ we define the $C_\varepsilon$-norm $\|\cdot\|_\varepsilon$ by
\[
\|h\|_\varepsilon=\sum_{k=0}^\infty\varepsilon_k\sup_{\R\times[0,1]\times M}|d^kh|.
\]
Denote by $C^\infty_{\varepsilon,I}(M)$ the space of functions $h\in C^\infty(\R\times[0,1]\times \widehat{M})$ satisfying that
$h$ is supported in $I\times[0,1]\times M$, and the $C_\varepsilon$-norm of $h$ is finite, i.e. $\|h\|_\varepsilon<\infty$.
It can be shown that $C^\infty_{\varepsilon,I}(M)$ is a separable Banach space, see~\cite[Lemma~8.3.2]{AD} or \cite{GS}. Moreover, there is a sequence $\varepsilon$ such that the set $C^\infty_{\varepsilon,I}(M)$ is dense in the space $C^\infty_I(M)$ of smooth functions $h:\R\times[0,1]\times \widehat{M}\to\R$ with compact supports in $I\times[0,1]\times M$ with respect to $C^1$-topology,  see~\cite[Proposition~8.3.1]{AD} or~\cite{Fl0,GS}. In the following we fix such $\varepsilon$.

\begin{prop}\label{prop:reg}
Let $(H,J)$ be a pair of a stationary homotopy in $\cH_{<\tau}$ and an almost complex structure of contact type outside of\; $M$ such that $(H^\pm,J)$ are Floer-regular. Let $I\subset\R$ be a compact interval with non-empty interior. Then there is a residual subset $\mathcal{V}^{reg}_\varepsilon\subset C^\infty_{\varepsilon,I}(M)$ such that for every $h\in\mathcal{V}^{reg}_\varepsilon$ the pair $(H+h, J)$ is Floer-regular.
\end{prop}

For the pair $(H,J)$ given as in the above proposition, we consider the solutions $u:\mathbb{R}\times [0,1]\to \widehat{M}$ of the PDE
\begin{equation}\label{app:seq}
(\overline{\partial}_{H,J}(u))(s,t):=\partial_su(s,t)+J_t(u(s,t))\partial_tu(s,t)+\nabla_JH^s_t(u(s,t))=0
\end{equation}
with finite energy and subject to the boundary condition $u(\cdot,\{0,1\})\subset \widehat{L}$, and denote by $\mathcal{M}_{H,J}$ the set of all such solutions $u$. For two Hamiltonian chords $x_\pm\in \mathcal{C}(L,H^\pm)$, we denote by $\mathcal{M}_{H,J}(x_-,x_+)$ the set of the above solutions $u$ with asymptotic chords $x_\pm$ at $\pm\infty$. It can be shown that $\mathcal{M}_{H,J}=\bigcup_{x_\pm\in \mathcal{C}(L,H^\pm)}\mathcal{M}_{H,J}(x_-,x_+)$.

Let $C^\infty_{exp}(x_-,x_+)$ denote the space of smooth maps $u:\R\times[0,1]\to M$ converging to $x_\pm$ at the ends with exponentially decaying derivatives and satisfying $u(\cdot,\{0,1\})\subset L$. Since outside of\;$M$ the slopes $\tau_{H^s}$ of $H^s$ do not depend on $s\in\R$, it follows from Lemma~\ref{lem:bd} that every solution $u$ to (\ref{app:seq}) has image contained in $int(M)$. Using the condition that $H^\pm\in\cH_{<\tau}$ are non-degenerate one can further show that each $u\in \mathcal{M}_{H,J}(x_-,x_+)$ belongs to $C^\infty_{exp}(x_-,x_+)$.

For $k>2/p$ and $p>1$, we define
\[
\begin{split}
\mathcal{P}(x_-,x_+):=\big\{u:\R\times[0,1]\to M\big|&u(s,t)=\exp_{w(s,t)}\xi(s,t)\;\hbox{where}\;w\in C^\infty_{exp}(x_-,x_+),\\
&\xi\in W^{k,p}(w^*TM)\;\hbox{with}\;\xi(s,0)\in T_{u(s,0)}L\\&\hbox{and}\;\xi(s,1)\in T_{u(s,0)}L\big\}.
\end{split}
\]
Let $W^{k-1,p}(x_-,x_+)$ denote the Banach vector bundle over $\mathcal{P}(x_-,x_+)$ whose fiber at $u\in \mathcal{P}(x_-,x_+)$ is $W^{k-1,p}(u^*TM)$.
%consists of vector fields $\xi\in L^p(u^*TM)$ with $\xi(\cdot,i)\in T_{u(\cdot,i)}L, i=0,1$.

As in~\cite{Fl0}, to prove Proposition~\ref{prop:reg} it suffices to prove

\begin{lem}\label{lem:Op}
The section
\[
\begin{split}
\mathcal{F}:\mathcal{P}(x_-,x_+)\times C^\infty_{\varepsilon,I}(M)&\longrightarrow W^{k-1,p}(x_-,x_+)\\
(u,h)&\longmapsto\overline{\partial}_{H+h,J}(u)
\end{split}
\]
is smooth and its linearization is surjective on its zero set
$$\mathcal{Z}(x_-,x_+)=\big\{(u,h)\in\mathcal{P}(x_-,x_+)\times C^\infty_{\varepsilon,I}(M)\big|\overline{\partial}_{H+h,J}(u)=0\big\}.$$
\end{lem}

To prove the smoothness of $\mathcal{F}$, we choose an unitary trivialization (preserving the symplectic structure and the almost complex structure)
$$\Phi:u^*T\widehat{M}\to \R\times[0,1]\times \R^{2n}$$
such that $\Phi(u(\cdot,0)^*T\widehat{L})=\R\times\{0\}\times\R^n$ and $\Phi(u(\cdot,1)^*T\widehat{L})=\R\times\{1\}\times \R^n$. Under this trivialization, $\mathcal{P}(x_-,x_+)$ is modeled over the Banach space
\[
W^{k,p}_L(\R\times [0,1];\R^{2n})=\big\{\xi\in W^{k,p}(\R\times [0,1];\R^{2n})\big|\xi(\cdot,0),\xi(\cdot,1)\in\R^n\big\},
\]
and the linear operator $D\mathcal{F}$ has the form
\[
\begin{split}
\Upsilon:W^{k,p}_L(\R\times [0,1];\R^{2n})\times C^\infty_{\varepsilon,I}(M)&\longrightarrow W^{k-1,p}(\R\times [0,1];\R^{2n})\\
(\xi,\eta)&\longmapsto D(\overline{\partial}_{H+h,J})_u(\xi)+\nabla_u\eta.
\end{split}
\]
The smoothness follows from the above form immediately.

Note that the operator $E_u:=D(\overline{\partial}_{H+h,J})_u$ is of perturbed Cauchy-Riemann type, i.e.
\[
E_u=\overline{\partial} +T=\frac{\partial}{\partial s}+j\frac{\partial}{\partial t}+T
\]
where $T:\R\times [0,1]\to \hbox{End}(\R^{2n})$, and has the asymptotic limits of the form
$\overline{\partial} +T^\pm$ with $T^\pm:[0,1]\to \hbox{End}(\R^{2n})$.
Since $x_\pm$ are non-degenerate, $E_u$ are Fredholm operators for all $u\in\mathcal{M}_{H,J}(x_-,x_+)$ and have index
\[\hbox{Ind}(E_u)=\mu(x_-)-\mu(x_+),\]
see for instance~\cite{RS,Fl0}.

%\begin{lem}\label{reg} Let $u$ be a solution of (\ref{app:seq}) which is not constant and has limits $x_\pm\in \mathcal{C}(L,H^\pm)$. Then the subset $\mathcal{R}(u)\subset \R\times(0,1)$ of those $(s,t)$ which satisfies\[\partial_su(s,t)\neq0,\quad u(s,t)\notin u(\R\setminus\{s\},t)\cup\{x_+(t)\}\cup\{x_-(t)\}\]is open and dense in $\R\times[0,1]$.\end{lem}

Before proving the rest of the statement in Lemma~\ref{lem:Op}, we show how Proposition~\ref{prop:reg} follows from this lemma. Notice that the linearization $D\mathcal{F}$ equals to the sum of a Fredholm operator and a linear operator, it has right inverse, see~\cite[Lemma 8.5.6]{AD}. This, together with Lemma~\ref{lem:Op}, implies that $\mathcal{F}$ intersects the zero section transversally. By the implicit function theorem, the zero set $\mathcal{Z}(x_-,x_+)$ is a smooth Banach submanifold of $\mathcal{P}(x_-,x_+)\times C^\infty_{\varepsilon,I}(M)$. The fact that $C^\infty_{\varepsilon,I}(M)$ is separable implies that $\mathcal{Z}(x_-,x_+)$ is also separable. Clearly, the projection map
\[\pi:\mathcal{Z}(x_-,x_+)\longrightarrow C^\infty_{\varepsilon,I}(M)\quad (u,h)\longmapsto h\]
is smooth. Moreover, $\pi$ is a Fredholm map which has the same Fredholm index as $D(\overline{\partial}_{H+h,J})_u$'s. Since $\pi^{-1}(h)=\mathcal{M}_{H+h,J}(x_-,x_+)$, if $h\in C^\infty_{\varepsilon,I}(M)$ is a regular value of $\pi$ then $D(\overline{\partial}_{H+h,J})_u$ are surjective for all $u\in\mathcal{M}_{H+h,J}(x_-,x_+)$, i.e. $(H+h,J)$ is Floer-regular. By the Sard-Smale Theorem, the set of regular values of $\pi$ is of the second category in $C^\infty_{\varepsilon,I}(M)$ (a countable intersection of open and dense set). We define the set $\mathcal{V}^{reg}_\varepsilon\subset C^\infty_{\varepsilon,I}(M)$ as the intersection of the regular values of the projections for all pairs $(x_-,x_+)$ where $x_\pm\in \mathcal{C}(L,H^\pm)$.

\begin{proof}[Sketch of the proof of Lemma~\ref{lem:Op}]
Notice that adding any $h\in C^\infty_{\varepsilon,I}(M)$ to the homotopy $H$ does not change the ends $H^\pm$ of the resulting homotopy.
Consequently, $E_u$ has a closed range with finite-dimensional cokernel for every $u\in\mathcal{Z}(x_-,x_+)$, and hence the range of $\Upsilon$ is closed. Therefore, to prove $\Upsilon$ is surjective on $\mathcal{Z}(x_-,x_+)$ it suffices to prove that the image $\im(\Upsilon)$ is dense.

Suppose the contrary, then there exists a nonzero continuous linear functional $\Gamma$ on $W^{k-1,p}(\R\times [0,1];\R^{2n})$ such that
\[\Gamma(E_u\xi)=0\quad\forall \xi\in W^{k,p}_L(\R\times [0,1];\R^{2n})\]
and
\[
\Gamma(\nabla_u\eta)=0\quad\forall \eta\in C^\infty_{\varepsilon,I}(M).
\]
By the elliptic regularity theory, $\Gamma$ can be represented by some nonzero vector field $\gamma\in W^{l,q}$ with $1/p+1/q=1$ for some $l\in\N$ so that for every $\zeta\in W^{k-1,p}$,
\[
\Gamma(\zeta)=\langle\gamma,\zeta\rangle=\int_{\R\times[0,1]}\gamma(s,t)\cdot\zeta(s,t)dsdt.
\]
In particular, for every $\xi\in W^{k,p}_L(\R\times [0,1];\R^{2n})$ and every $\eta\in C^\infty_{\varepsilon,I}(M)$,
\begin{equation}\label{e:E}
\langle\gamma,E_u\xi\rangle=0,
\end{equation}
\begin{equation}\label{e:H}
\langle\gamma,\nabla_u\eta\rangle=0.
\end{equation}
By (\ref{e:E}), $\gamma$ is in the kernel of the $L^2$-adjoint operator $E_u^*$ associated to $E_u$ which is also a perturbed Cauchy-Riemannian operator. The unique continuation~\cite[Proposition~3.1]{FHS} implies that if $\gamma$ has an infinite-order zero, then it is identically zero. To arrive at a contradiction we shall show that $\gamma$ vanishes on $I\times[0,1]$ (hence $\gamma\equiv0$).

We define the map
\[
\begin{split}
\widetilde{u}:\R\times[0,1]&\longrightarrow \R\times[0,1]\times \widehat{M}\\
(s,t)&\longmapsto \big(s,t,u(s,t)\big).
\end{split}
\]
Clearly, this map is an embedding. We pull back $\gamma$ as the vector field along $\widetilde{u}$ to $\R\times[0,1]\times \widehat{M}$ that has no components in the directions $\partial/\partial s\in T\R$ and $\partial/\partial t\in T[0,1]$. We see that $\gamma$ is not tangent to $\im(\widetilde{u})$ at the points where it is not zero.

Now we suppose that there exists a point $(s_0,t_0)\in I\times[0,1]$ such that $\gamma(s_0,t_0)\neq0$. Without loss of generality, we may further assume that $(s_0,t_0)$ is an interior point in $I\times[0,1]$, then there exists a small open neighborhood $\Omega$ in $int(I\times[0,1])$ such that $\gamma(s,t)$ is nonzero on $\Omega$. Hence in this neighborhood $\gamma$ is transversal to $\widetilde{u}$. Pick a smooth function $\chi:\R\times[0,1]\to\R$ with support in $U$ which satisfies
$$\int_{\R\times[0,1]}\chi(s,t)dsdt\neq 0.$$
Let $U\subset I\times[0,1]\times int(M)$ be a tubular neighborhood of $\widetilde{u}(\Omega)$ such that $U\cap\im(\widetilde{u})=\widetilde{u}(\Omega)$. We pick a smooth function $h:\R\times[0,1]\times \widehat{M}\to\R$ with support in $U$ such that if $\phi_{s,t}(r)$ is a parameterized integral curve of $\gamma$ passing through $\widetilde{u}(s,t)$ at $r=0$, then
\[
h(s,t,\phi_{s,t}(r)):=\chi(s,t)r\quad \forall r\in(-\epsilon,\epsilon)
\]
where $\epsilon$ is sufficiently small. The condition that $\gamma$ is transversal to $\widetilde{u}(\Omega)$ guarantees that such $h$ can be well defined. Write $h_{s,t}=h(s,t,\cdot)$. Then we compute
\[
\begin{split}
\langle\gamma,\nabla_uh\rangle&=\int_{\R\times[0,1]}\gamma(s,t)\cdot \nabla_uh(s,t)dsdt\\
&=\int_{\R\times[0,1]}dh_{s,t}(\gamma(s,t))dsdt\\
&=\int_{\R\times[0,1]}dh_{s,t}\bigg(\frac{d}{dr}\bigg|_{r=0}\phi_{s,t}(r)\bigg)dsdt\\
&=\int_{\R\times[0,1]}\frac{d}{dr}\bigg|_{r=0}h_{s,t}\big(\phi_{s,t}(r)\big)dsdt\\
&=\int_{\R\times[0,1]}\chi(s,t)dsdt\neq0.
\end{split}
\]
From the construction we see that $h\in C^\infty_I(M)$. Using the fact that $C^\infty_{\varepsilon,I}(M)$ is dense in the space $C^\infty_I(M)$ with the $C^1$-topology, one can choose $\eta\in C^\infty_{\varepsilon,I}(M)$ which approximates arbitrarily to $h$ so that the equality~(\ref{e:H}) does not hold for this $\eta$. So we achieve a contradiction which means that $\Upsilon$ is surjective.
\end{proof}

\subsection{Gromov-Floer compactness and robustness of barricades under perturbations}

\begin{prop}\label{prop:survive}
Let $(H^s)_{s\in\R}\subset \cH_{<\tau}$ be a homotopy, stationary for $|s|\geq R$ with some $R>0$, from $H^-$ to $H^+$ such that
$H^\pm$ are non-degenerate and the pairs $(H,J)$ and $(H^\pm,J)$ have a barricade on $\Omega:=M_r\setminus M_{r'}^\circ$. Then, for every $C^\infty$-small perturbation $\widetilde{H}$ of $H$ which satisfies $supp (\partial_s\widetilde{H}^s-\partial_sH^s)\subset [-R,R]\times [0,1]\times M$ and $\widetilde{H}^\pm=H^\pm$, the pairs $(\widetilde{H},J)$ and $(\widetilde{H}^\pm,J)$ also have a barricade on $\Omega$.

\end{prop}

The proof of Proposition~\ref{prop:survive} is parallel to that of~\cite[Proposition 9.21]{GS},
we give a sketch of the proof for the sake of completeness.
The crucial ingredient of the proof is to apply the following Gromov-Floer compactness result.

\begin{thm}\label{prop:compact}
Let $(H^s)_{s\in\R}\subset \cH_{<\tau}$ be a homotopy, stationary for $|s|\geq R$ with some $R>0$, with non-degenerate ends $H^\pm$. Let $H_n$ be a sequence of homotopies in $\cH_{<\tau}$ such that the sequence $\{H_n-H\}_{n}\subset C^\infty_I(M)$ converges to $0$ in $C^\infty$-topology where $I=[-R,R]$.
Let $\{u_n\}\subset\mathcal{M}_{H_n,J}(x_-,x_+)$ be a sequence of solutions and $\{\sigma_n\}$ a sequence of real numbers. Then there exist subsequences of $\{u_n\}$ and $\{\sigma_n\}$ (still denoted by $\{u_n\}$ and $\{\sigma_n\}$ for simplicity), Hamiltonian chords $x_i\in\mathcal{C}(L,H^-),i=0,\ldots,k$ and $y_j\in\mathcal{C}(L,H^+),j=0,\ldots,l$  satisfying $x_0=x_-$ and $y_l=x_+$, and sequences of real numbers $\{\varsigma_n^i\}$ for $1\leq i\leq k$ and $\{\tau_n^j\}$ for $1\leq j\leq l$ such that
\[
\begin{split}
u_n&\longrightarrow w\in \mathcal{M}_{H,J}(x_k,y_0),
\\
u_n(\cdot+\varsigma_n^i,\cdot)&\longrightarrow v_i\in \mathcal{M}_{H^-,J}(x_{i-1},x_i),\\u_n(\cdot+\tau_n^j,\cdot)&\longrightarrow v_j'\in \mathcal{M}_{H^+,J}(y_{j-1},y_j)
\end{split}
\]
for $1\leq i\leq k$ and $1\leq j\leq l$ in $C^\infty_{loc}$-topology, and the sequence $u_n(\cdot+\sigma_n)$ converges to one of $v_i,w,v_j'$ in $C^\infty_{loc}$-topology up to a shift in the $s$-coordinate.

\end{thm}

\begin{figure}[H]
\centering
\qquad\includegraphics[scale=0.5]{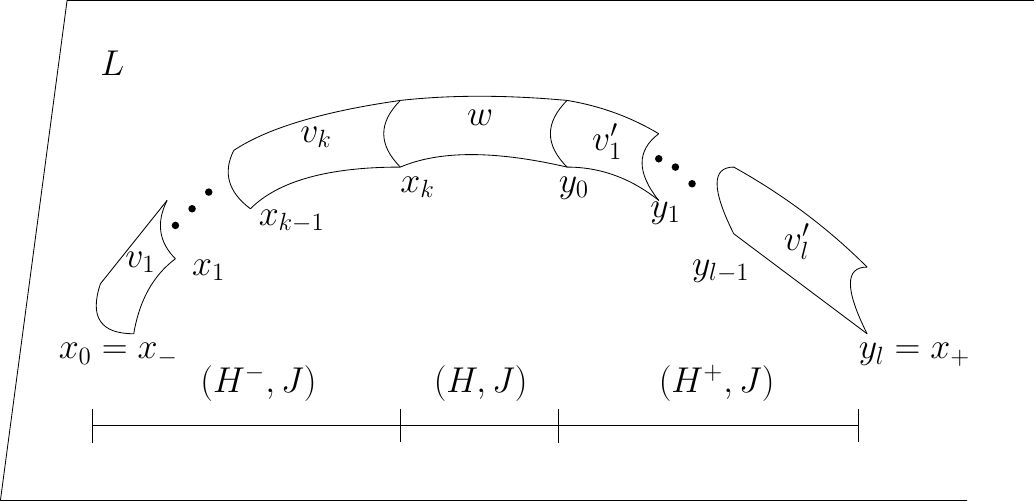}
\caption{Broken Floer trajectories.}\label{fig:GroF}
\end{figure}

For $x\in\mathcal{C}(L,H^-)$ and $y\in\mathcal{C}(L,H^+)$, we write
$$\mathcal{M}:=\cup_{n}\mathcal{M}_{H_n,J}\bigcup \mathcal{M}_{H,J},\quad \mathcal{M}(x,y):=\cup_{n}\mathcal{M}_{H_n,J}(x,y)\bigcup \mathcal{M}_{H,J}(x,y).$$

The finite sequence $(v_0,\ldots,v_k,w,v_0'\ldots,v_l')$ is called a \emph{broken trajectory} of $(H,J)$ as illustrated in Figure~\ref{fig:GroF}.
The proof of \cite[Theorem~11.1.10]{AD} or \cite[Proposition~9.14]{GS} can be carried over to Theorem~\ref{prop:compact} in a direct fashion. The key observation is that the energies $E(u)$ for all $u\in \mathcal{M}$ have a uniform bound. For simplicity we set $H_0:=H$. In our setting, even although the target manifold $\widehat{M}$ is non-compact, since $H_n\in \cH_{<\tau}$ for all $n\in\N\cup\{0\}$ and their slopes outside of\;$M$ do not depend on $s$ we deduce from Lemma~\ref{lem:bd} that $\im(u)\subset M$. And since $\partial_sH_n$ are supported in $[-R,R]\times[0,1]\times M$ it follows from the energy identity~(\ref{E}) that for every $u\in\mathcal{M}_{H_n,J}(x_-,x_+)$,
\[E(u)\leq\mathcal{A}_{L,H^-}(x_-)-\mathcal{A}_{L,H^+}(x_+)+2R\cdot\sup_n\sup_{[-R,R]\times [0,1]\times M}\big\{\partial_sH^s_n\big\}\]
where the fact that all $\partial_sH^s_n$ satisfy a uniform bound on $[-R,R]\times [0,1]\times M$ has been used due to the uniform convergence with derivatives of the sequence $\{H_n\}$ to $H$. Letting
\[
A:=\max_{x_\pm\in \mathcal{C}(L,H^\pm)}\big(\mathcal{A}_{L,H^-}(x_-)-\mathcal{A}_{L,H^+}(x_+)\big)++2R\cdot\sup_n\sup_{[-R,R]\times [0,1]\times M}\big\{\partial_sH^s_n\big\},
\]
we have that $E(u)\leq A$ for all $u\in\mathcal{M}$. This and the assumption that $L$ is exact (hence no bubbling disks or spheres occur) will result in a uniform bound for the $J$-gradient vector field $\nabla_J u$, that is,
\begin{lem}[{\cite[Proposition 11.3.8]{AD} or \cite[Lemma 9.15]{GS}}]
There exists a constant $C>0$ such that for each $u\in \mathcal{M}$ and each $(s,t)\in \R\times[0,1]$,
\[\|\partial_su(s,t)\|^2_J+\|\partial_tu(s,t)\|_J^2\leq C.\]
\end{lem}
The rest of the proof of Theorem~\ref{prop:compact} is a consequence of repeatedly applying the above lemma, Arzel\'{a}-Ascoli theorem and elliptic regularity,  we refer to~\cite[Theorem 11.1.10]{AD} for the completely parallel proof.

\begin{proof}[Sketch of the proof of Proposition~\ref{prop:survive}]
Let $(H^s)_{s\in\R}\subset \cH_{<\tau}$ be a homotopy, stationary for $|s|\geq R$ with some $R>0$, from $H^-$ to $H^+$ such that
$H^\pm$ are non-degenerate and the pairs $(H,J)$ and $(H^\pm,J)$ have a barricade on $\Omega:=M_r\setminus M_{r'}^\circ$. We will see that the similar restrictions of the barricade on Floer trajectories persist for broken trajectories of $(H,J)$. More precisely, we have
\begin{clm}\label{clm:broken}
Let $\overline{v}=(v_1,\ldots,v_k,w,v_1'\ldots,v_l')$ be a broken trajectory of $(H,J)$ connecting $x_\pm\in \mathcal{C}(L,H^\pm)$. Then it holds that
\begin{enumerate}
  \item[I.] if $\im(x_-)\subset W_0:=M_{r'}$ then $\overline{v}\subset W_0$.
  \item[II.] if $\im(x_+)\subset W_1:=M_r$ then $\overline{v}\subset W_1$.
\end{enumerate}

\end{clm}
\begin{proof}[Proof of Claim~\ref{clm:broken}]
We only prove statement I since statement II can be proved in almost exactly the same way. We first notice that $x_1:=\lim_{s\to+\infty}v_1(s,\cdot)$ has image in $W_0$ because $(H^-,J)$ has a barricade on $\Omega$ and $x_-=\lim_{s\to-\infty}v_1(s,\cdot)$ has image in $W_0$. Since by definition $x_1$ is also the negative end of $v_2$, i.e.  $x_1:=\lim_{s\to-\infty}v_2(s,\cdot)$ and since $(H^-,J)$ has a barricade on $\Omega$, we see that $\im(v_2)\subset W_0$. Repeatedly using the above argument we find that the images of $v_1,\ldots,v_k$ are all contained in $W_0$. Next, by definition we have $x_k:=\lim_{s\to+\infty}v_k(s,\cdot)=\lim_{s\to-\infty}w(s,\cdot)$. It follows from the assumption that $(H,J)$ has a barricade on $\Omega$ that $\im(w)\subset W_0$. Since $(H^+,J)$ has a barricade on $\Omega$ and $y_0:=\lim_{s\to-\infty}v_0'(s,\cdot)=\lim_{s\to+\infty}w(s,\cdot)$, we see that $v_0'$ has image in $W_0$. Finally, arguing in the same way we conclude that all $v_j',1\leq j\leq l$ have images in $W_0$. Therefore, the broken trajectory $\overline{v}$ is entirely contained in $W_0$.
This completes the proof of the claim.
\end{proof}

Now we are in position to finish the proof of Proposition~\ref{prop:survive}. Suppose the contrary that statement I in the definition of barricade does not hold. Then there exists a sequence $\{H_n\}$ of regular homotopies with the sequence $\{H_n-H\}_{n}\subset C^\infty_{[-R,R]}(M)$ converging to $0$ in $C^\infty$-topology such that for every $n$, the moduli space $\mathcal{M}_{H_n,J}$ admits an element $u_n$ which satisfies that the limit $x_-^n:=\lim_{s\to-\infty}u_n(s,\cdot)$ is contained in $W_0$ but $u_n$ is not. For every $n$ we pick a number $\sigma_n\in\R$ such that $u_n(\sigma_n,\cdot)$ is not contained in $W_0$.

Since $H^\pm_n=H^\pm$ and $H^\pm$ admit only finitely many Hamiltonian chords, we may assume that $x_\pm^n=x_\pm$ (after passing to a subsequence) for all $n$. In this case we have $u_n\in\mathcal{M}(x_-,x_+)$ for all $n$. By Theorem~\ref{prop:compact}, there exist subsequences of $\{u_n\}$ and $\{\sigma_n\}$ (still denoted by $\{u_n\}$ and $\{\sigma_n\}$) such that $\{u_n\}$ converges to a broken trajectory $\overline{v}$ of $(H,J)$ and $u_n(\cdot+\sigma_n,\cdot)$ converges to one of the solutions in $\overline{v}$ (perhaps up to a shift). Since $x_-=x_0\subset W_0$, the first statement of Claim~\ref{clm:broken} implies that $\overline{v}$ is completely contained in $W_0$, and hence $\lim_{n\to\infty}u_n(\cdot+\sigma_n,\cdot)\subset W_0$. Since the latter limit is taken in $C^\infty_{loc}$-topology, we have that
$\lim_{n\to\infty}u_n(\sigma_n,\cdot)=\lim_{n\to \infty}u_n(0+\sigma_n,\cdot)$ is also contained in $W_0$ --- a contradiction!

Using the second statement of Claim~\ref{clm:broken} and arguing in the same way, we see that if $n$ is sufficiently large, then every solution $u_n\in\mathcal{M}_{H_n,J}$ ending in $W_1$ is contained in $W_1$.
\end{proof}

\subsection{Finishing the proof of Theorem~\ref{barricade}}

\begin{proof}
Let $(H^s)_{s\in\R}$ be the homotopy as in the assertion. Let $(J_t)_{t\in[0,1]}\in \mathcal{J}_\theta$ be any family of almost complex structures (to be determined later) that
  are of contact type near the boundaries $\partial M_{r'},\partial M_r$. For sufficiently small $\mu\in(0,\tau)$ we make a $C^\infty$-small perturbation of $H$ into a homotopy $h$ such that $(h,J)$ admits a cylindrical bump of slope $\mu\in(0,\tau)$ on $\Omega:=M_r\setminus M_{r'}^\circ$ and $h^\pm\in \cH_{<\tau}$ are non-degenerate. This can be done by adding first a $C^\infty$-small radial bump function $\chi$ to $H$ on a small neighborhood of $\Omega$ in $M$, and then perturbing its ends $H^\pm+\chi$ out of a small neighborhood of $\partial \Omega$ in $M$ into non-degenerate Hamiltonians that are $C^2$-small Morse functions on $[r',1]\times\partial M$, finally making a generic perturbation of $H+\chi$ out of a small neighborhood of $\partial \Omega$ in $M$ so that the ends $h^\pm$ of the resulting homotopy $h$ agree with the non-degenerate perturbations of $H^\pm+\chi$ respectively. Clearly, $(h,J)$ is a pair with a cylindrical bump of slope $\mu\in(0,\tau)$ on $\Omega$. It follows from Proposition~\ref{prop:bump} that the pairs $(h,J)$ and $(h^\pm,J)$ have a barricade on $\Omega$. Moreover, we may require that $supp(\partial_sh-\partial_s H)\subset [-R,R]\times [0,1]\times M$ for large $R>0$ since $H$ is stationary for sufficiently large $|s|$.

The remaining problem is that the pairs $(h,J)$ and $(h^\pm,J)$ may be not regular. Suppose this case, we need to perturb the homotopy and its ends again to achieve regularity while keeping the barricade condition on $\Omega$. Indeed, this is possible by following an argument due to Ganor and Tanny~\cite[Section~9]{GS}. The crucial point is that if $(h^\pm,J)$ are regular, then for any compact interval $I$ with non-empty interior one can perturb $h$ on the set $I\times [0,1]\times M$ such that the resulting homotopy $h'$ is Floer-regular meanwhile the barricade property of $h$ on $\Omega$ survives for $h'$.
Now since $h^\pm$ are non-degenerate and has no chords near $\partial \Omega$, by the usual Floer theory one can choose generically a family of time-dependent almost complex structures $J=(J_t)_{t\in[0,1]}\in\mathcal{J}_\theta$ such that $(h^\pm,J)$ are regular and of contact type near the boundary $\partial \Omega$.
Fix such $J$. If the homotopy $h$ constructed above is a constant homotopy, we have already got the desired result. If $h$ is not a constant homotopy, by Proposition~\ref{prop:reg} one can choose a homotopy $h'$ close enough to $h$ with $supp(\partial_sh'-\partial_sh)\subset I\times [0,1]\times M$ for a fixed compact interval $I$ such that $(h',J)$ are regular. Then Proposition~\ref{prop:survive} gives rise to the desired pairs $(h',J)$ and $(h'^\pm,J)$.
\end{proof}

%I am thankful to Adrian P. Dawid for the insightful discussions and for sharing his master's thesis, supervised by Paul Biran, which employs a different method to demonstrate the unboundedness of the Lagrangian Hofer norm on Lagrangians isotopic to a fiber in the disk cotangent bundle of any sphere $S^n$ ($n\geq 1$) with the standard metric.

 \end{document}